\documentclass[sn-mathphys-num]{sn-jnl}

\usepackage{graphicx}%
\usepackage{multirow}%
\usepackage{amsmath,amssymb,amsfonts}%
\usepackage{amsthm}%
\usepackage{mathrsfs}%
\usepackage[title]{appendix}%
\usepackage{xcolor}%
\usepackage{textcomp}%
\usepackage{manyfoot}%
\usepackage{booktabs}%
\usepackage{listings}%
\usepackage{enumitem}

\usepackage{algorithm}
\usepackage{algorithmic}

\usepackage{url}
\usepackage{longtable}
\usepackage{caption}
\usepackage[labelformat=simple]{subcaption} % Suppresses parentheses around subfigure 
\usepackage{siunitx}
\usepackage{amsmath,amsthm,amssymb,dsfont,amsthm,upgreek,bbm}
\usepackage{graphicx,color,tikz,float,rotating,enumitem}
\usepackage{verbatim}
\hypersetup{%
  colorlinks=true,% hyperlinks will be coloured
  linkcolor={blue!66!black},% hyperlink text will be green
  linkbordercolor=red,% hyperlink border will be red
  citecolor={blue!66!black}
  }
\definecolor{indigo}{RGB}{75,0,130}
\usepackage{array}
\newcolumntype{L}[1]{>{\raggedright\let\newline\\\arraybackslash\hspace{0pt}}m{#1}}
\newcolumntype{C}[1]{>{\centering\let\newline\\\arraybackslash\hspace{0pt}}m{#1}}
\newcolumntype{R}[1]{>{\raggedleft\let\newline\\\arraybackslash\hspace{0pt}}m{#1}}

\usetikzlibrary{shapes.geometric, arrows}
\usepackage[title]{appendix}
\usepackage{arydshln}

\newcommand{\ignore}[1]{}

\usepackage{multirow}
\usepackage[symbol]{footmisc}

\usepackage{hyperref}

\begin{document}

\title{A Stochastic Programming Model for Anticipative Planning of Integrated Electricity and Gas Systems with Bidirectional Energy Flows under Fuel and CO2 Price Uncertainty}

\author*[1]{\fnm{Giovanni} \sur{Micheli}
\email{giovanni.micheli@unibg.it}}
\equalcont{These authors contributed equally to this work.}

\author[1]{\fnm{Maria Teresa} \sur{Vespucci}
\email{maria-teresa.vespucci@unibg.it}}
\equalcont{These authors contributed equally to this work.}

\author[2]{\fnm{Alessia} \sur{Cortazzi}
\email{alessia.cortazzi@cesi.it}}
\equalcont{These authors contributed equally to this work.}

\author[2]{\fnm{Cinzia} \sur{Puglisi}
\email{cinzia.puglisi@cesi.it}}
\equalcont{These authors contributed equally to this work.}

\affil*[1]{\orgdiv{Department of Management, Information and Production Engineering}, \orgname{University of Bergamo}, 
\orgaddress{\street{Via Salvecchio 19}, \city{Bergamo}, \postcode{24129}, \country{Italy}}}

\affil[2]{\orgname{CESI}, 
\orgaddress{\street{Via Rubattino 54}, \city{Milano}, \postcode{20134}, \country{Italy}}}

\abstract{A two-stage multi-period mixed-integer linear stochastic programming model is proposed to assist qualified operators in long-term generation and transmission expansion planning of electricity and gas systems to meet policy objectives.
The first-stage decisions concern investments in new plants, new connections in the electricity and gas sectors, and the decommissioning of existing thermal power plants; the second-stage variables represent operational decisions, with uncertainty about future fuel and CO$_2$ prices represented by scenarios.
The main features of the model are:
(i) the bidirectional conversion between electricity and gas enabled by Power-to-Gas and thermal power plants,
(ii) a detailed representation of short-term operation, crucial for addressing challenges associated with integrating large shares of renewables in the energy mix, and
(iii) an integrated planning framework to evaluate the operation of flexibility resources, their ability to manage non-programmable generation, and their economic viability.
Given the computational complexity of the proposed model, in this paper we also implement a solution algorithm based on Benders decomposition to compute near-optimal solutions.
A case study on the decarbonisation of the Italian integrated energy system demonstrates the effectiveness of the model. The numerical results show:
(i) the importance of including a detailed system representation for obtaining reliable results, and
(ii) the need to consider price uncertainty to design adequate systems and reduce overall costs.}


\maketitle

%\linenumbers
%\nolinenumbers

\section{ Introduction}
\label{sec:Intro}

To meet ambitious decarbonisation targets, such as a 55 \% reduction in greenhouse gas emissions by 2030 compared to 1990 levels \cite{EU}, energy systems will need to transition to large shares of non-programmable renewable generation capacity, particularly wind and solar photovoltaic (PV).
Electricity production from non-programmable renewable sources poses two distinct problems: (i) primary energy sources (wind, solar radiation) may not be available at times of high demand, and excess renewable energy production may occur at times of low demand: flexibility resources are therefore needed to shift energy in time, and (ii) renewable power plants are located in geographical areas with favourable climatic conditions, which are not necessarily the most industrialised and energy-intensive areas: flexibility resources are therefore needed to shift energy in space.
The energy system must be able to cope with these issues in order to make the best use of non-programmable renewable production. 
Various technologies can provide this flexibility. 
Programmable hydro and gas-fired thermal power plants are well suited for load balancing, as they can be started up within minutes when non-programmable renewable generation is unable to meet the load requirements.
Pumped hydro energy storage and battery energy storage shift energy in time \cite{luo2015overview}. 
Spatial shift is related to the transport of energy through electricity and gas transmission systems.
Another promising option is Power-to-Gas (PtG) technology \cite{jentsch2014optimal}, which uses electricity to produce gaseous fuels such as hydrogen and biomethane, which can be stored locally for later use (energy shift in time) or injected into the natural gas grid (energy shift in time and/or space).
As a result, the electricity and gas systems are becoming increasingly interconnected through bidirectional energy flows: gas-fired thermal power plants generate electricity, and excess electricity can be used to produce hydrogen and biomethane. Planning for the transition to a decarbonised energy system requires an integrated view of the electricity and gas systems.

Generation and Transmission Expansion Planning (GTEP) analysis provides the strategic framework for addressing this challenge. 
In energy system design, GTEP models determine how the system should evolve over the long term by defining investments in new generation and transmission. 
The GTEP problem can be approached in two main ways: from a centralised or decentralised perspective. 
If the goal of the analysis is to design the trajectory of energy systems, a centralised approach to GTEP is necessary. This approach allows for a comprehensive and unified plan where a single entity, such as a regulator or government agency, can define the strategic direction for system development. The centralised framework allows modeling the entire energy system in detail, taking into account the complex interactions between generation, transmission and storage components, while ensuring alignment with long-term decarbonisation goals.
The centralised approach is used to carry out \textit{anticipative planning}, i.e., to determine the most efficient configuration of the energy system in order to set policies and incentives that guide independent decision-makers to invest in a socially efficient manner \cite{pozo2013if}.
Once the centralised model has determined the optimal system configuration, decentralised models can be used to assess the behaviour of individual market participants and to validate policies in a competitive environment. 
Decentralised models reduce the technical detail of the analysis and focus on the interactions between market participants and their decision-making processes in a liberalised electricity sector \cite{hesamzadeh2020transmission}.
In this paper, we address the centralised GTEP problem with the aim of assisting qualified operators in defining a strategic plan for the transition to decarbonised energy systems. 
There is a large literature on the GTEP problem. 
We refer the reader to \cite{koltsaklis2018state,micheli2021survey,klatzer2022state} for comprehensive reviews on the topic.

\begin{sidewaystable}
\caption{Comparison of relevant contributions on expansion planning of integrated systems proposed in the literature.}\label{tab:review}
\begin{tabular*}{\textheight}{@{\extracolsep\fill}lllllllllll}
\toprule%
Ref. & Formulation & GEP & TEP & Storage & PtG & Stochasticity & Dynamism & Operations & UC & Policy goals \\
\midrule
%[28]
\cite{qiu2015linear} & LP & T & P/G & GSU & No & Deterministic & Dynamic & Annual & No & --  \\
%[57]
\cite{zhang2015long} & LP & T & P/G & -- & No & Deterministic & Dynamic & Load blocks & No & --  \\
%Mio [10]
\cite{khaligh2018multi} & MILP & T & P/G & -- & No & Deterministic & Static & Load blocks & No & -- \\
%[33]
\cite{zhang2015reliability} & MILP & -- & P/G & -- & No & Deterministic & Static & Load blocks & No & -- \\
%[58]
\cite{samsatli2018multi} & MILP & T/W & P/G & BESS/GSU & No & Deterministic & Dynamic & Rep. days & No & --  \\
%[59]
\cite{zhang2015optimal} & MILP & T & P & -- & No & Deterministic & Dynamic & Load blocks & No & --   \\
%Mio [9] - Loro [17]
\cite{unsihuay2010model} & MILP & T/W & P/G & PHS/GSU & No & Deterministic & Dynamic & Load blocks & No & -- \\
%Mio [8]
\cite{wang2017coordinated} & MINLP & T & P/G & GSU & No & Deterministic & Static & Single time & No & --  \\
%[60]
\cite{sanchez2016convex} & MINLP & T/W & P/G & -- & No & Deterministic & Static & Single time & No & --  \\
%Mio [11] - Loro [56]
\cite{chaudry2014combined} & MINLP & -- & P/G & GSU & No & Deterministic & Dynamic & Rep. days & No & --  \\
%[31]
\cite{barati2014multi} & MINLP & T & P/G & -- & No & Deterministic & Dynamic & Annual & No & -- \\
%[41]
\cite{qiu2014low} & MINLP & T & P/G & -- & No & Deterministic & Dynamic & Annual & No & Carbon emissions  \\
%Mio [13] - Loro [73]
\cite{zhao2017coordinated} & MILP & T & P/G & -- & No & Stochastic & Static & Load blocks & No & 
\\
%Mio [14] - Loro [44]
\cite{nunes2018stochastic} & MILP & T/W/S & P/G & -- & No & Stochastic & Static & Time blocks & No & RES generation \\
%[30]
\cite{yamchi2021multi} & MILP & T/S & P/G & -- & No & Stochastic & Static & Time blocks & No & -- \\
%Mio [15]
\cite{ding2017multi} & MILP & T & P/G & -- & No &  Stochastic & Dynamic & Load blocks & No & -- \\
%[54]
\cite{pantovs2013stochastic} & MILP & T & P/G & -- & No & Stochastic & Dynamic & Load blocks & No & -- \\
%[35]
\cite{khaligh2019stochastic} & MINLP & T/W & P/G & -- & No & Stochastic & Static & Load blocks & No & Bounds for wind capacity
\\
%Mio [12] - Loro [34]
\cite{odetayo2018chance} & MINLP & T & G & GSU & No & Stochastic & Dynamic & Annual & No & --  \\
%%% With P2G
%[82]
\cite{brown2018synergies} & LP & T/W/S & P & BESS/PHS & Yes & Deterministic & Static & Hourly & No & Cap on CO$_2$ emissions  \\
%[40]
\cite{tao2020carbon} & MILP & W/S & P/G & -- & Yes & Deterministic & Dynamic & Load blocks & No & Coal retirement \\
%[96]
\cite{zeng2017bi} & MINLP & T & -- & GSU & Yes & Deterministic & Dynamic & Rep. week & No & --  \\
%Nuovi
\cite{liang2020stochastic} & MILP & T & P/G & -- & Yes & Stochastic & Dynamic & Load blocks & No & -- \\
%He
\cite{he2017robust} & MILP & T & P/G & -- & Yes & Stochastic & Dynamic & Load blocks & No & -- \\
%%%
%This Paper
%%%
\hline \hline
& \multirow{4}{*}{MILP} & \multirow{4}{*}{T/W/S} & \multirow{4}{*}{P/G} & & \multirow{4}{*}{Yes} & \multirow{4}{*}{Stochastic} & \multirow{4}{*}{Dynamic} & \multirow{1}{*}{Rep. days with} & \multirow{4}{*}{Yes} & Coal retirement,  \\
This & & & & BESS/PHS/ & &&& constraints to & & RES penetration, \\
work & & & & GSU & &&& simulate the 
 & & cap on CO$_2$ emissions, \\
 & & & &  & & && full chronology & & wind and PV capacity \\
\botrule
\end{tabular*}
\footnotetext{LP, Linear Programming; MILP, Mixed Integer Linear Programming; MINLP, Mixed Integer Nonlinear Programming; GEP, Generation Expansion Planning; T, Thermal plants; W, Wind power plants; S, Solar PV power plants; TEP, Transmission Expansion Planning; P, Power lines; G, Gas pipelines; GSU, Gas Storage Unit; PHS, Pumped Hydro Storage; BESS, Battery Energy Storage System; Rep., Representative; UC, Unit Commitment.}
\end{sidewaystable}

Most of the work presented in the literature neglects the interactions between the electricity and gas systems. Only a few papers deal with planning the expansion of integrated systems. Table \ref{tab:review} provides the following information about the relevant papers on integrated system expansion planning: the type of mathematical model used (column 2); the power generation technologies considered for Generation Expansion Planning (GEP) (column 3); whether Transmission Expansion Planning (TEP) is considered for either electricity or gas (column 4); the storage technologies considered (column 5); whether PtG plants are considered (column 6); whether the model is deterministic or stochastic (column 7) and static or dynamic (column 8); the approach taken to represent system operation (column 9); whether Unit Commitment (UC) constraints are considered (column 10); the policy objectives (column 11).

Regarding the mathematical formulation, Mixed Integer Linear Programming (MILP) is the most common model chosen for the GTEP problem, as discrete investment decisions need to be represented. A few studies use Linear Programming (LP) models, mainly because of their computational efficiency and ease of implementation \cite{qiu2015linear,zhang2015long,brown2018synergies}.
Models that incorporate detailed representations of natural gas systems are often based on Mixed Integer Nonlinear Programming (MINLP), where nonlinearities typically arise from modelling the steady-state gas flow equations \cite{wang2017coordinated,sanchez2016convex,chaudry2014combined,barati2014multi,qiu2014low,khaligh2019stochastic,odetayo2018chance}.
While these formulations enhance the technical accuracy of gas flow representation, they also introduce significant computational challenges in finding optimal solutions \cite{klatzer2022state}.

Another key feature of capacity planning models is the expansion scope. 
Models are typically classified according to whether they consider GEP, TEP or both.
Table \ref{tab:review} details the generation technologies considered in the GEP expansion decisions, whether thermal (T) and/or wind (W) and/or solar photovoltaic (S). It also indicates whether the TEP concerns the electricity grid (P) or the natural gas grid (G) or both.
In most of the studies in the literature, non-programmable generation is an exogenous parameter of the model that determines the investment in thermal power plants, mainly fueled by natural gas, to satisfy the net load. However, this approach is not suitable for determining the expansion of energy systems if the achievement of decarbonisation targets requires investment mainly in non-programmable renewable generation technologies. 

%For this reason, references \cite{nunes2018stochastic,brown2018synergies} develop mathematical programs to plan investments in both programmable and non-programmable power generation.

In integrated energy systems with significant shares of non-programmable electricity generation, resources are needed to provide the flexibility required to compensate for the intermittency of renewable generation. However, many studies completely neglect the contribution of flexibility resources. Table \ref{tab:review} shows that only a few works analyse the role of specific flexibility resources, such as gas storage units (GSU) \cite{qiu2015linear,wang2017coordinated,chaudry2014combined,odetayo2018chance,zeng2017bi}, pumped hydro storage (PHS) \cite{unsihuay2010model,brown2018synergies}, battery energy storage systems (BESS) \cite{samsatli2018multi,brown2018synergies} and PtG plants \cite{brown2018synergies,tao2020carbon,zeng2017bi,liang2020stochastic,he2017robust}, in supporting the integration of large shares of non-programmable renewable generation. These resources are considered in our proposed model.

In the GTEP problem, the planning and operation of power systems are inherently affected by stochasticity. 
The existing literature considers various sources of uncertainty, including electricity demand \cite{odetayo2018chance,khaligh2019stochastic,nunes2018stochastic,pantovs2013stochastic,zhao2017coordinated,liang2020stochastic,ding2017multi}, natural gas demand \cite{zhao2017coordinated}, non-dispatchable renewable generation \cite{yamchi2021multi,khaligh2019stochastic,nunes2018stochastic,liang2020stochastic,he2017robust}, natural gas prices \cite{nunes2018stochastic,liang2020stochastic}, CO$_2$ prices \cite{liang2020stochastic}, and interest rates \cite{khaligh2019stochastic}. 
Approaches developed to deal with uncertainty include chance-constrained programming \cite{odetayo2018chance}, robust optimization \cite{he2017robust}, and stochastic programming \cite{yamchi2021multi,khaligh2019stochastic,nunes2018stochastic,pantovs2013stochastic,zhao2017coordinated,liang2020stochastic,ding2017multi}. 
Contributions in the area of stochastic programming can be further divided into two classes: two-stage and multi-stage stochastic models. 
In two-stage models \cite{yamchi2021multi,khaligh2019stochastic,nunes2018stochastic,zhao2017coordinated}, a single optimal investment plan is determined: expansion decisions (first-stage variables) are made before the uncertainty is revealed, and operational decisions (second-stage variables) depend on the realised values of the uncertain parameters. 
Multi-stage models \cite{pantovs2013stochastic,liang2020stochastic,ding2017multi}, on the other hand, define multiple investment plans with sequential investment decisions that depend on uncertainty that is gradually revealed over time. The computational burden of multi-stage problems requires a reduction in the level of detail of the analysis, typically representing short-term power system operation with low accuracy.

Table \ref{tab:review} also indicates, for each paper, whether the proposed approach is \textit{static} or \textit{dynamic}. 
A static model refers to a setting in which all investment decisions are taken at once, typically for a specific target year, without representing the system's evolution over time. 
In contrast, a dynamic model plans the evolution of the system over time, following the adopted time discretisation.
%By contrast, dynamic models allow for investment decisions at several points in time.

Another distinctive feature of the GTEP models reviewed in Table \ref{tab:review} is how short-term operations are represented. 
In principle, the most accurate approach would be to model all hours within the planning horizon. 
However, this is generally not computationally feasible due to the complexity and size of the resulting problem.
Only \cite{brown2018synergies}, where a deterministic linear program is used, considers the hourly resolution of the planning horizon, which consists of one year.
To address the computational limitations of considering longer planning horizons and multiple sources of stochasticity, four main modeling approaches can be identified in the literature, based on the level of detail used to represent operational conditions:
(i) models using annual averages for demand and non-programmable generation \cite{qiu2015linear,wang2017coordinated,sanchez2016convex,barati2014multi,qiu2014low,odetayo2018chance};
(ii) models that consider a limited number of values, based on duration curves \cite{zhang2015long,zhang2015optimal,unsihuay2010model,khaligh2018multi,khaligh2019stochastic,pantovs2013stochastic,tao2020carbon,he2017robust,zhang2015reliability,zhao2017coordinated,ding2017multi,liang2020stochastic};
(iii) models that divide each year into blocks of consecutive hours (e.g., days, weeks, or months) \cite{yamchi2021multi,nunes2018stochastic};
(iv) models that represent short-term operation with hourly resolution over a reduced set of representative periods \cite{samsatli2018multi,chaudry2014combined,zeng2017bi}, typically selected by clustering methods \cite{micheli2021selecting}. 
The first three approaches result in a highly simplified representation of system operation, fail to capture the intermittency and short-term variability of renewable generation, and are therefore not suitable for planning the transition to decarbonised integrated energy systems. 
In contrast, using representative periods with hourly resolution provides a more detailed representation of short-term dynamics, but limits the ability to account for seasonal effects, which are critical for accurately modeling long-term storage operation.

The increasing penetration of intermittent renewables needs to be supported by the increasing use of flexibility resources to manage the variability and uncertainty of generation from non-programmable sources. By describing the operation of thermal power plants through a detailed model that includes UC constraints, the thermoelectric capacity required to provide flexibility to the system is more accurately determined, as is the cost of the flexibility provided by such plants. Nevertheless, UC constraints are not considered in any of the papers \cite{qiu2015linear} to \cite{he2017robust}. As discussed in \cite{schwele2020unit}, this simplification is not appropriate for planning expansion with large shares of non-programmable generation.

Finally, with regard to decarbonisation targets, the last column of Table \ref{tab:review} shows that only a few papers consider the role of specific policy targets in defining the expansion of integrated systems, such as CO$_2$ emission caps \cite{qiu2014low,brown2018synergies}, minimum levels of non-programmable renewable energy \cite{nunes2018stochastic}, caps on installed wind power capacity \cite{khaligh2019stochastic} and coal phase-out \cite{tao2020carbon}.

A number of key features emerge from the literature review that need to be considered when designing an effective GTEP model that ensures adequate representation of operational dynamics and decarbonisation constraints:
\begin{itemize}
    \item develop an integrated planning framework that incorporates expansion decisions on programmable and non-programmable generation, electricity and gas networks, electricity and gas storage, and PtG plants to determine the combination of generation and flexibility resources that can meet decarbonisation targets at the lowest cost;
    \item capture the inherent stochasticity that affects the planning and operation of energy systems;
    \item address the dynamics of the problem by planning the evolution of the system over time;
    \item provide an accurate representation of short-term operation, representing hourly dispatch, while taking into account the seasonality of medium- and long-term storage and ensuring computational tractability;
    \item include UC constraints to effectively assess flexibility needs;
    \item include policy constraints in the mathematical formulation to assess the feasibility of transitioning to decarbonised energy systems. 
%    \item Develop solution methods to provide GTEP plans whose quality can be guaranteed with affordable effort.
\end{itemize}

Building on these key considerations, in this paper we develop a two-stage multi-period stochastic mixed-integer linear programming model for the anticipative planning of integrated electricity-gas systems with bidirectional conversion between electricity and gas, to help qualified operators define an optimal investment plan for the transition to decarbonised energy systems. 
We focus on the long-term uncertainty of future prices of fuels and of CO$_2$ emissions as they play a crucial role in expansion decisions. 
Indeed, fuel costs and CO$_2$ emission prices influence the operating costs of power plants and thus the merit order of power dispatch. As fuel prices fluctuate according to the geopolitical situation and CO$_2$ prices are influenced by political decisions, they introduce an element of uncertainty into long-term decisions.
Forecast scenarios developed by internationally recognised bodies are therefore considered by decision makers and analysts.
Such scenarios can be used to represent uncertainty in a stochastic programming approach.
Unlike the two-stage models proposed in the literature \cite{yamchi2021multi,khaligh2019stochastic,nunes2018stochastic,zhao2017coordinated}, where the expansion takes place in a single period before the uncertainty is revealed, the model we propose is dynamic, i.e., the first-level decisions (investments in new generation, storage, PtG units, power transmission lines and gas pipelines, as well as the decommissioning of existing thermal power plants) can be implemented in any of the years in which the planning period is discretised. 
The two-stage approach allows short-term system operation, modeled by second-stage decision variables, to be represented by hourly discretisation, taking into account unit commitment constraints for thermal power plants, reserve requirements and the medium-term seasonality of hydro and gas storage: the model thus determines how to effectively manage large shares of renewable generation, supplementing it with thermal, hydro or battery power at times and in system zones where it is insufficient, and storing it in electricity or gas storage for future use at times and in system zones where it exceeds demand.
The last row of Table \ref{tab:review} summarises the main features of the proposed method, highlighting its novelty compared to the existing literature.
Given the computational complexity of the proposed model, in this paper we also implement a solution algorithm based on Benders decomposition to compute near-optimal solutions.
The paper is organised as follows. 

The proposed stochastic programming model is presented in Section \ref{sec:Model}. Section \ref{sec:Algorithm} presents the solution algorithm for determining the GTEP plans. 
Section \ref{sec:Case_study} presents a case study related to the Italian energy system. The optimal expansion plan determined by the stochastic programming model (SPM plan) is compared with the one determined by the mean value problem (MVP plan), i.e. the corresponding deterministic model, where each uncertain parameter is assigned the expected value of its realisations in the scenarios used to represent the uncertainties in the stochastic programming model. In the SPM plan, investment in new renewable generation capacity is higher than in the MVP plan. The two expansion plans also differ in the choice of flexibility resources: the SPM plan installs more resources for shifting energy in time, while the MVP plan installs more resources for shifting energy in space. Furthermore, the MVP plan does not invest in the new PtG technology, whereas the SPM plan does.  The numerical results show that in scenarios with high prices for both gas and CO$_2$ emissions, the differences between the two plans affect the adequacy of the system.  
Section \ref{sec:Conclusions} concludes the paper.

\section{The two-stage stochastic programming model}
\label{sec:Model}
In this section we introduce the two-stage multi-period mixed-integer linear stochastic programming model for the anticipative planning of integrated electricity and gas systems. 
The sets, parameters and variables used in the model are listed in \ref{Notation}.  

\begin{align}
\label{eq:obj}
\min & \quad \sum_{y \in \mathcal{Y}} \frac{1}{(1+r)^{y-y_0}} \Bigg\{ 
\sum_{l \in \mathcal{L}_C} IC_l^\text{L} \  \delta_{l,y} + 
\sum_{j \in \mathcal{J}_C} IC_j^\text{J} \ \delta_{j,y} + 
\sum_{h \in \mathcal{H}_C} IC_h^\text{H} \ \overline{H}_h^{\text{OUT}} \ \delta_{h,y} + \nonumber \\
% 2
& + \sum_{z \in \mathcal{Z}} \Bigg[ 
\sum_{k \in \mathcal{K}_z} \big( IC_{k,y}^\text{K} \ \overline{P}_k \ N_{k,y}^{+} + DC_{k,y}^\text{K} \ \overline{P}_k \ N_{k,y}^{-} \big) + IC_{z,y}^\text{S} \ S_{z,y} \ + \nonumber \\
% 3
& + IC_{z,y}^\text{W} \ W_{z,y}  +  \sum_{b \in \mathcal{B}_z} IC_{b,y}^\text{B} \ B^\text{CAP}_{b,y}  +  \sum_{g \in \mathcal{G}_z} IC_{g}^\text{PtG} \ PtG^\text{CAP}_{g,y} \Bigg] \Bigg\} + \nonumber \\ 
% 4
& + \sum_{w \in \mathcal{W}} pb_w \Bigg\{ \sum_{y \in \mathcal{Y}} \sum_{c \in \mathcal{C}^y} \psi_c  \sum_{t=1}^{24} \Bigg[ \sum_{z \in \mathcal{Z}}  \Bigg(
\sum_{h \in \mathcal{H}_z}  C_h^\text{H}  \ H_{h,t,c,y,w}^{\text{OUT}} + \sum_{k \in \mathcal{K}_z} C_k^{\text{SU}}\ \alpha_{k,t,c,y,w} + \nonumber \\
\nonumber \\
% 5
& + \sum_{k \in \mathcal{K}_z} C_{k,y,w}^\text{M} \left( \! \underline{P}_k \gamma_{k,t,c,y,w}  +  p_{k,t,c,y,w} \right) + \sum_{b \in \mathcal{B}_z}  C_b^{B}  \ B_{b,t,c,y,w}^{\text{OUT}} + \nonumber \\
%& \quad \quad \quad \quad \  + \sum_{k \in \mathcal{K}_z} \left( OM_k + {CO_2}_k \ Pr^\text{CO$_2$}_{y,w} \right) \left( \underline{P}_k \ \gamma_{k,t,c,y,w} + p_{k,t,c,y,w} \right) + \nonumber \\
%&\quad \quad \quad \quad \  + \sum_{k \in \mathcal{K}_z \setminus \mathcal{K}_z^\text{G}} HR_k \ Pr^\text{FUEL}_{f(k),y,w} \left( \underline{P}_k \ \gamma_{k,t,c,y,w} + p_{k,t,c,y,w} \right)  + \nonumber \\
% 6
& + C^{\text{OG}} \ OG_{z,t,c,y,w} + C^{\text{ENP}} E^{\text{NP}}_{z,t,c,y,w} + C^{\text{RNP}} R^{\text{NP}}_{z,t,c,y,w} \Bigg) + \nonumber \\
% 7
&  +   \sum_{n \in \mathcal{N}}  \Bigg(  
C^{\text{G}}_{n,y,w} G_{n,t,c,y,w}  +
 \sum_{g \in \mathcal{G}_n} C^\text{PtG}_g G^{\text{PtG}}_{g,t,c,y,w}+ 
 C^{\text{GC}}  G^{\text{CURT}}_{n,t,c,y,w} 
\Bigg) \Bigg]  \Bigg\} 
\end{align}    
subject to
{\allowdisplaybreaks
\medskip
\begin{itemize}
    \item for $y \in \mathcal{Y}$:
\end{itemize}
\begin{alignat}{2}
% VINCOLO 2
& \theta_{l,y} = \sum_{i=1}^y \delta_{l,i}, && l \in \mathcal{L}_C, \label{eq:inv_line} \\
% VINCOLO 3
& \delta_{l,y} \ \theta_{l,y} \in \{0,1\}, && l \in \mathcal{L}_C, \label{eq:bin_line} \\ 
% VINCOLO 4
& \theta_{j,y} = \sum_{i=1}^y \delta_{j,i}, && j \in \mathcal{J}_C, \label{eq:inv_pipe}\\
% VINCOLO 5
& \delta_{j,y}, \ \theta_{j,y} \in \{0,1\}, \qquad && j \in \mathcal{J}_C, \label{eq:bin_pipe} \\
% VINCOLO 6
& \theta_{h,y} = \sum_{i=1}^y \delta_{h,i}, && h \in \mathcal{H}_C, \label{eq:inv_hydro}\\
% VINCOLO 7
& \delta_{h,y}, \ \theta_{h,y} \in \{0,1\}, && h \in \mathcal{H}_C. \label{eq:bin_hydro}
\end{alignat}

\bigskip
\begin{itemize}
    \item for $z \in \mathcal{Z}, y \in \mathcal{Y}$:
\end{itemize}
\begin{alignat}{2}
%\text{subject to} \\
% VINCOLO 8
& N_{k,y} = N_{k,y-1} + N_{k,y}^{+} - N_{k,y}^{-}, \quad && k \in \mathcal{K}_z, \label{eq:inv_the}\\
% VINCOLO 9
& \underline{N}_{k,y} \leq N_{k,y} \leq \overline{N}_{k,y}, && k \in \mathcal{K}_z, \label{eq:bound_the}\\
% VINCOLO 10
& N_{k,y}, \ N_{k,y}^{+}, \ N_{k,y}^{-} \in \mathbb{N}, && k \in \mathcal{K}_z, \label{eq:int_the}\\
% VINCOLO 11
& \underline{S}_{z,y} \leq S_{z,0} + \sum_{i=1}^y S_{z,i} \leq \overline{S}_{z,y}, \label{eq:bound_sol}\\
% VINCOLO 12
& \underline{W}_{z,y} \leq W_{z,0} + \sum_{i=1}^y W_{z,i} \leq \overline{W}_{z,y}, \qquad && \label{eq:bound_wind}\\
% VINCOLO 13
& S_{z,y}, \ W_{z,y} \geq 0, &&  \label{eq:pos_sw}\\
% VINCOLO 15
& B^\text{CAP}_{b,y} \geq 0, && b \in \mathcal{B}_z, \label{eq:pos_bat}\\
% VINCOLO 17
& PtG^\text{CAP}_{g,y} \geq 0, && g \in \mathcal{G}_z.  \label{eq:pos_ptg}
\end{alignat}

\bigskip
\begin{itemize}
    \item for $z \in \mathcal{Z}$:
\end{itemize}
\begin{alignat}{2}
% VINCOLO 14
& B^\text{CAP}_{b,0} + \sum_{y \in \mathcal{Y}} B^\text{CAP}_{b,y} \leq \overline{B}^\text{CAP}_{b}, && b \in \mathcal{B}_z, \label{eq:inv_bat} \\
% VINCOLO 16
& PtG^\text{CAP}_{g,0} + \sum_{y \in \mathcal{Y}} PtG^\text{CAP}_{g,y} \leq \overline{PtG}^\text{CAP}_{g}, \qquad && g \in \mathcal{G}_z.  \label{eq:inv_ptg}
\end{alignat}
}
% VINCOLO 18
\bigskip
\begin{itemize}
    \item for $z \in \mathcal{Z}, 1 \leq t \leq 24, c \in \mathcal{C}^y, y \in \mathcal{Y}$:
\end{itemize}
\begin{align}
% VINCOLO 18
& RES_{z,t,c,y}
= \Big( S_{z,0} + \sum_{i=1}^y S_{z,i} \Big) \ \mu_{z,t,c} 
+ \Big( W_{z,0} + \sum_{i=1}^y W_{z,i} \Big) \ \rho_{z,t,c}.  \label{eq:RES}
\end{align}

\bigskip
\begin{itemize}
    \item for $z \in \mathcal{Z}, 1 \leq t \leq 24, c \in \mathcal{C}^y, y \in \mathcal{Y}, w \in \mathcal{W}$:
\end{itemize}
\begin{align}
% VINCOLO 19
& RES_{z,t,c,y} 
+ \sum_{k \in \mathcal{K}_z} (\underline{P}_k \ \gamma_{k,t,c,y,w} +
p_{k,t,c,y,w}) 
+ \sum_{h \in \mathcal{H}_z} H_{h,t,c,y,w}^{\text{OUT}} 
+ \sum_{b \in \mathcal{B}_z} B_{b,t,c,y,w}^{\text{OUT}} 
+ \nonumber \\
& + \sum_{l \in BS^{\text{P}}_z} F^\text{L}_{l,t,c,y,w} + E^{\text{NP}}_{z,t,c,y,w} 
= D^{\text{P}}_{z,t,c,y}
+ \sum_{g \in \mathcal{G}_z} \frac{G^{\text{PtG}}_{g,t,c,y,w}}{\eta_g}
+ \nonumber \\
& + \sum_{h \in \mathcal{H}_z} H_{h,t,c,y,w}^{\text{IN}}  
+ \sum_{b \in \mathcal{B}_z} B_{b,t,c,y,w}^{\text{IN}}
+ \sum_{l \in FS^{\text{P}}_z} F^\text{L}_{l,t,c,y,w} 
+ OG_{z,t,c,y,w},  \label{eq:balance} \\
& E^{\text{NP}}_{z,t,c,y,w} \geq 0, \label{eq:non_neg_ENP}\\
& OG_{z,t,c,y,w} \geq 0. \label{eq:non_neg_OG}
\end{align}

\bigskip
\begin{itemize}
    \item for $1 \leq t \leq 24, c \in \mathcal{C}^y, y \in \mathcal{Y}, w \in \mathcal{W}$:
\end{itemize}
\begin{alignat}{2}
% VINCOLO 20
&\underline{F}^\text{L}_l \leq F^\text{L}_{l,t,c,y,w} \leq \overline{F}^\text{L}_l, && l \in \mathcal{L}_E, \label{eq:flow_el} \\
% VINCOLO 21
&\underline{F}^\text{L}_l \ \theta_{l,y} \leq F^\text{L}_{l,t,c,y,w} \leq \underline{F}^\text{L}_l \ \theta_{l,y}, \qquad && l \in \mathcal{L}_C.  \label{eq:flow_el_NEW}
\end{alignat}

\bigskip
\begin{itemize}
\item for $1 \leq t \leq 24, c \in \mathcal{C}^y,y \in \mathcal{Y}, w \in \mathcal{W}$:
\end{itemize}
\begin{alignat}{2}
% VINCOLO 22
& 0 \leq H_{h,t,c,y,w}^{\text{OUT}} \leq \overline{H}_h^{\text{OUT}}, && h \in \mathcal{H}_E,  \label{eq:hydro_OUT}\\
% VINCOLO 23
& 0 \leq H_{h,t,c,y,w}^{\text{OUT}} \leq \overline{H}_h^{\text{OUT}} \ \theta_{h,y}, \qquad && h \in \mathcal{H}_C,  \label{eq:hydro_OUT_new} \\
% VINCOLO 24
& 0 \leq H_{h,t,c,y,w}^{\text{OUT}} \leq F_{h,t,c,y}, && h \in \mathcal{H}_{NP},  \label{eq:hydro_OUT_NP}\\
% VINCOLO 25
& 0 \leq H_{h,t,c,y,w}^{\text{IN}} \leq \overline{H}_h^{\text{IN}}, && h \in \mathcal{H}_E \cap \mathcal{H}_{P}, \label{eq:hydro_IN}\\
% VINCOLO 26
& 0 \leq H_{h,t,c,y,w}^{\text{IN}} \leq \overline{H}_h^{\text{IN}} \ \theta_{h,y}, && h \in \mathcal{H}_C \cap \mathcal{H}_{P}, \label{eq:hydro_IN_new}\\
% VINCOLO 27
& 0 \leq H_{h,t,c,y,w}^{\text{SPILL}} \leq \overline{H}_h^{\text{SPILL}}, && h \in \mathcal{H}_E \cap \mathcal{H}_{P},  \label{eq:hydro_SPILL}\\
% VINCOLO 28
& 0 \leq H_{h,t,c,y,w}^{\text{SPILL}} \leq \overline{H}_h^{\text{SPILL}} \ \theta_{h,y}, && h \in \mathcal{H}_C \cap \mathcal{H}_{P}.  \label{eq:hydro_SPILL_new}
\end{alignat}

% VINCOLI da 29 a 32
\bigskip
\begin{itemize}
    \item for $h \in (\mathcal{H}_E \cup \mathcal{H}_C) \cap \mathcal{H}_P,y \in \mathcal{Y}, w \in \mathcal{W}$:
\end{itemize}
\begin{alignat}{2}
% VINCOLO 29
& 
H^{\text{LT}}_{h,M,y,w} 
= H^{\text{LT}}_{h,0} 
+ \sum_{d=1}^{M} \sum_{t=1}^{24} \Big( F_{h,t,c(d),y} \ \theta_{h,y}  
+ \lambda_h^{\text{IN}} \ H^{\text{IN}}_{h,t,c(d),y,w} + 
&& \nonumber \\
& \quad 
- \lambda_h^{\text{OUT}} \ H^{\text{OUT}}_{h,t,c(d),y,w}
- H^{\text{SPILL}}_{h,t,c(d),y,w} \Big), && \label{eq:hydro_LT_0} \\
% VINCOLO 30
& H^{\text{LT}}_{h,\xi M,y,w} 
= H^{\text{LT}}_{h,(\xi-1)M,y,w} 
+ \sum_{d=(\xi-1)M}^{\xi M} \sum_{t=1}^{24} \Big( F_{h,t,c(d),y} \ \theta_{h,y} +
&& \nonumber \\
& \quad
+ \lambda_h^{\text{IN}} \ H^{\text{IN}}_{h,t,c(d),y,w} 
- \lambda_h^{\text{OUT}} \ H^{\text{OUT}}_{h,t,c(d),y,w}
- H^{\text{SPILL}}_{h,t,c(d),y,w} \Big), && 2 \leq \xi \leq \overline{\xi}, \label{eq:hydro_LT_M} \\
% VINCOLO 31
& H^{\text{LT}}_{h,0} 
= H^{\text{LT}}_{h,\overline{\xi}M,y,w} 
+ \sum_{d=\overline{\xi}M+1}^{365} \sum_{t=1}^{24} \Big( F_{h,t,c(d),y} \ \theta_{h,y} 
+ \qquad && \nonumber \\
& + \lambda_h^{\text{IN}} \ H^{\text{IN}}_{h,t,c(d),y,w} 
- \lambda_h^{\text{OUT}} \ H^{\text{OUT}}_{h,t,c(d),y,w}
- H^{\text{SPILL}}_{h,t,c(d),y,w} \Big), &&  \label{eq:hydro_LT_final} \\
% VINCOLO 32
& 0 \leq H^{\text{LT}}_{h,\xi M,y,w} \leq 
EPR_h \ \overline{H}_h^{\text{OUT}}, && 1 \leq \xi \leq \overline{\xi}.  \label{eq:hydro_LT_pos}
\end{alignat}

% VINCOLI da 33 a 38
\bigskip
\begin{itemize}
    \item for $k \in \mathcal{K}_z, z \in \mathcal{Z},c \in \mathcal{C}^y,y \in \mathcal{Y}, w \in \mathcal{W}$:
\end{itemize}
\begin{alignat}{2}
% VINCOLO 33
&\gamma_{k,t,y,c,w} \leq N_{k,y}, && 1 \leq t \leq 24,  \label{eq:the_max_ON} \\
% VINCOLO 34
&\gamma_{k,t,y,c,w} = \gamma_{k,t-1,y,c,w} + \alpha_{k,t,y,c,w} - \beta_{k,t,y,c,w}, \qquad && 1 \leq t \leq 24, \label{eq:the_status}  \\ 
% VINCOLO 35
&\sum_{\tau = t - (MUT_k-1)}^t \alpha_{k,\tau,y,c,w} \leq \gamma_{k,t,y,c,w}, && MUT_k \leq t \leq 24,  \label{eq:the_MUT} \\
% VINCOLO 36
&\sum_{\tau = t - (MDT_k-1)}^t \beta_{k,\tau,y,c,w} \leq N_{k,y}-\gamma_{k,t,y,c,w}, && MDT_k \leq t \leq 24,  \label{eq:the_MDT} \\
% VINCOLO 37
&0 \leq p_{k,t,c,y,w} \leq (\overline{P}_k - \underline{P}_k) \ \gamma_{k,t,c,y,w}, && 1 \leq t \leq 24,  \label{eq:the_max_power} \\
% VINCOLO 38
&\alpha_{k,t,c,y,w}, \ \beta_{k,t,c,y,w}, \ \gamma_{k,t,c,y,w} \in \mathbb{N}, && 1 \leq t \leq 24.  \label{eq:the_int} 
\end{alignat}

% VINCOLI da 39 a 40
\bigskip
\begin{itemize}
    \item for $z \in \mathcal{Z}, 1 \leq t \leq 24, c \in \mathcal{C}^y, y \in \mathcal{Y}, w \in \mathcal{W}$:
\end{itemize}
\begin{align}
% VINCOLO 39
& \sum_{k \in \mathcal{K}_z} \Big[ (\overline{P}_k-\underline{P}_k) \ \gamma_{k,t,c,y,w} - p_{k,t,c,y,w} \Big] + R^{\text{NP}}_{z,t,c,y,w} \geq R_{z,t,c,y},  \label{eq:reserve} \\
% VINCOLO 40
& R^{\text{NP}}_{z,t,c,y,w} \geq 0. \label{eq:pos_reserve} 
\end{align}

% VINCOLI da 41 a 42
\bigskip
\begin{itemize}
    \item for $a \in \mathcal{A},y \in \mathcal{Y}$:
\end{itemize}
\begin{align}
% VINCOLO 41
& \sum_{z \in \mathcal{Z}_a} \sum_{c \in \mathcal{C}^y} \psi_c \sum_{t=1}^{24} RES_{z,t,c,y} 
\geq 
\phi_{a,y} \ \Big( \sum_{z \in \mathcal{Z}_a} \sum_{c \in \mathcal{C}^y} \psi_c \sum_{t=1}^{24} D_{z,t,c,y}^\text{P} \Big),  \label{eq:penetration} \\ 
% VINCOLO 42
& \sum_{z \in \mathcal{Z}_a} 
\sum_{k \in \mathcal{K}_z} 
\sum_{c \in \mathcal{C}^y} 
\psi_c \sum_{t=1}^{24} 
\Big[ 
%HR_k \ 
{CO_2}_k \
(\underline{P}_k \ \gamma_{k,t,c,y,w} + p_{k,t,c,y,w} ) \Big] 
\leq {\overline{CO}_2}_{a,y}, \qquad w \in \mathcal{W}.  \label{eq:CO2} 
\end{align}

% VINCOLI da 43 a 47
\bigskip
\begin{itemize}
    \item for $b \in \mathcal{B}_z, z \in \mathcal{Z}, y \in \mathcal{Y}, w \in \mathcal{W}$:
\end{itemize}
\begin{alignat}{2}
% VINCOLO 43
& B_{b,t,c,y,w} 
= (1-\lambda_b) \ B_{b,t-1,c,y,w} 
+ \lambda_b^{\text{IN}} \ B_{b,t,c,y,w}^{\text{IN}} 
- \lambda_b^{\text{OUT}} \ B_{b,t,c,y,w}^{\text{OUT}}, \, \, \, && 1 \leq t \leq 24,  \label{eq:bat_level} \\
% VINCOLO 44
& 0 \leq B_{b,t,c,y,w} \leq EPR_b \ \Big( B^\text{CAP}_{b,0} + \sum_{i=1}^y B^\text{CAP}_{b,i} \Big), && 1 \leq t \leq 24,  \label{eq:bat_inst}\\
% VINCOLO 45
& 0 \leq B^\text{IN}_{b,t,c,y,w} \leq B^\text{CAP}_{b,0} + \sum_{i=1}^y B^\text{CAP}_{b,i}, && 1 \leq t \leq 24,  \label{eq:bat_IN}\\
% VINCOLO 46
& 0 \leq B^\text{OUT}_{b,t,c,y,w} \leq B^\text{CAP}_{b,0} + \sum_{i=1}^y B^\text{CAP}_{b,i}, && 1 \leq t \leq 24,  \label{eq:bat_OUT}\\
% VINCOLO 47
& B_{b,24,c,y,w} = B_{b,0,c,y}.  \label{eq:bat_final}
\end{alignat}

% VINCOLO 48
\bigskip
\begin{itemize}
    \item for $n \in \mathcal{N},1 \leq t \leq 24, c \in \mathcal{C}^y, y \in \mathcal{Y}, w \in \mathcal{W}$:
\end{itemize}
\begin{align}
& G_{n,t,c,y,w} 
+ \sum_{g \in \mathcal{G}_n} G_{g,t,c,y,w}^{\text{PtG}} 
+ G^\text{OUT}_{n,t,c,y,w} 
+ \sum_{j \in BS^{\text{G}}_n} F^J_{j,t,c,y,w} 
+ G^\text{CURT}_{n,t,c,y,w} = D^\text{G}_{n,t,c,y} + \nonumber \\
&   
+ \sum_{k \in \mathcal{K}^\text{G}_n} HR_k \ 
(\underline{P}_k \ \gamma_{k,t,c,y,w} + p_{k,t,c,y,w} )
+ G^\text{IN}_{n,t,c,y,w} 
+ \sum_{j \in FS^{\text{G}}_z} F^J_{j,t,c,y,w}, \label{eq:gas_balance} \\
& G^\text{CURT}_{n,t,c,y,w} \geq 0. \label{non_neg_GCURT}
\end{align}

% VINCOLI da 49 a 50
\bigskip
\begin{itemize}
    \item for $1 \leq t \leq 24, c \in \mathcal{C}^y, y \in \mathcal{Y}, w \in \mathcal{W}$:
\end{itemize}
\begin{alignat}{2}
% VINCOLO 49
&\underline{F}^\text{J}_j \leq F^\text{J}_{j,t,c,y,w} \leq \overline{F}^\text{J}_j, && j \in \mathcal{J}_E,  \label{eq:flow_gas} \\
% VINCOLO 50
&\underline{F}^\text{J}_j \ \theta_{j,y} \leq F^\text{J}_{j,t,c,y,w} \leq \underline{F}^\text{J}_j \ \theta_{j,y}, \qquad && j \in \mathcal{J}_C.  \label{eq:flow_gas_NEW}
\end{alignat}

% VINCOLO 51
\bigskip
\begin{itemize}
    \item for $1 \leq t \leq 24, g \in \mathcal{G}_n, n \in \mathcal{N}, c \in \mathcal{C}^y, y \in \mathcal{Y}, w \in \mathcal{W}$:
\end{itemize}
\begin{align}
& 0 \leq G_{g,t,c,y,w}^{\text{PtG}} \leq PtG^\text{CAP}_{g,0} + \sum_{i=1}^y PtG^\text{CAP}_{g,i}.  \label{eq:PtG_prod}
\end{align}

% VINCOLI da 52 a 54
\bigskip
\begin{itemize}
    \item for $1 \leq t \leq 24, c \in \mathcal{C}^y, y \in \mathcal{Y}, w \in \mathcal{W}$:
\end{itemize}
\begin{alignat}{2}
% VINCOLO 52
& \underline{G}_n \leq G_{n,t,c,y,w} \leq \overline{G}_n, && n \in \mathcal{N},  \label{eq:gas_prod} \\
% VINCOLO 53
& 0 \leq G^\text{IN}_{n,t,c,y,w} \leq \overline{G}_n^{\text{IN}}, && n \in \mathcal{N},  \label{eq:gas_IN} \\
% VINCOLO 54
& 0 \leq G^\text{OUT}_{n,t,c,y,w} \leq \overline{G}_n^{\text{OUT}}, \qquad && n \in \mathcal{N}. \label{eq:gas_OUT}
\end{alignat}

% VINCOLI da 55 a 58
\bigskip
\begin{itemize}
    \item for $n \in \mathcal{N},y \in \mathcal{Y}, w \in \mathcal{W}$:
\end{itemize}
\begin{alignat}{2}
% VINCOLO 55
& G^{\text{LT}}_{n,M,y,w} 
= G^{\text{LT}}_{n,0} 
+ \sum_{d=1}^{M} \sum_{t=1}^{24} 
\Big( \lambda_n^{\text{IN}} \ G^{\text{IN}}_{n,t,c(d),y,w} 
- \lambda_n^{\text{OUT}} \ G^{\text{OUT}}_{n,t,c(d),y,w} \Big), && \label{eq:gas_LT_0} \\
% VINCOLO 56
& G^{\text{LT}}_{n,\xi M,y,w}   =   G^{\text{LT}}_{n,(\xi-1)M,y,w}   +     \sum_{d = ( \xi- 1 ) M}^{\xi M}   \sum_{t=1}^{24}   \Big(   \lambda_n^{\text{IN}}   G^{\text{IN}}_{n,t,c(d),y,w}  + \nonumber \\
& \quad -   \lambda_n^{\text{OUT}}    G^{\text{OUT}}_{n,t,c(d),y,w}   \Big), \   \quad 2   \leq   \xi \leq   \overline{\xi},  \label{eq:gas_LT_M}\\
% VINCOLO 57
& G^{\text{LT}}_{n,0} 
= G^{\text{LT}}_{n,\overline{\xi}M,y,w} 
+ \sum_{d=\overline{\xi}M+1}^{365} \sum_{t=1}^{24} 
\Big( \lambda_n^{\text{IN}} \ G^{\text{IN}}_{n,t,c(d),y,w} - \lambda_n^{\text{OUT}} \ G^{\text{OUT}}_{n,t,c(d),y,w} \Big), &&  \label{eq:gas_LT_final} \\
% VINCOLO 58
& 0 \leq G^{\text{LT}}_{n,\xi M,y,w} \leq \overline{G}_n^{\text{LT}}, \qquad 1   \leq   \xi \leq  \overline{\xi}.  \label{eq:gas_LT_pos}
\end{alignat}

\setcounter{footnote}{0}

The power system is modeled as a set $\mathcal{Z}$ of zones connected by power transmission lines. 
The set of transmission lines existing in the initial configuration of the power system 
(i.e., at the beginning of the planning period) 
and the set of candidate lines for power transmission capacity expansion are denoted by $\mathcal{L}_E$ and $\mathcal{L}_C$, respectively. 
The set of lines entering zone $z$ is denoted by ${BS}^{\text{P}}_z \subseteq (\mathcal{L}_E \cup \mathcal{L}_C)$ and 
the set of lines leaving zone $z$ is denoted by ${FS}^{\text{P}}_z \subseteq (\mathcal{L}_E \cup \mathcal{L}_C)$.
Constraints \eqref{eq:inv_line} and \eqref{eq:bin_line} model the investment decisions regarding each candidate line $l \in \mathcal{L}_C$: either candidate line $l$ is not built ($\delta_{l,y}=0, \forall y \in \mathcal{Y}$) and therefore is not available in any year of the planning period ($\theta_{l,y}=0, \forall y \in \mathcal{Y}$), or it is built in year $\hat{y}$ ($\delta_{l,\hat{y}}=1, \delta_{l,y}=0$, for $y \ne \hat{y}$) and therefore is available from year $\hat{y}$ until the end of the planning period ($\theta_{l,y}=0$, for $y<\hat{y}$; $\theta_{l,y}=1$, for $y \geq \hat{y}$)\footnote{For the sake of simplicity, the construction times of generation and transmission facilities are not considered in this paper. However, the model can be readily extended to include construction times, as done in \cite{vespucci2014stochastic}.}. 
Similarly, the gas system is modeled as a set $\mathcal{N}$ of zones connected by gas pipelines: the set of pipelines existing in the initial configuration and the set of candidate pipelines for gas network expansion are denoted by $\mathcal{J}_E$ and $\mathcal{J}_C$, respectively. 
The set of pipelines entering zone $n$ is denoted by ${BS}^{\text{G}}_n \subseteq (\mathcal{J}_E \cup \mathcal{J}_C)$ and 
the set of lines leaving zone $n$ is denoted by ${FS}^{\text{G}}_n \subseteq (\mathcal{J}_E \cup \mathcal{J}_C)$.
Constraints \eqref{eq:inv_pipe} and \eqref{eq:bin_pipe} model the investment decisions regarding each candidate pipeline $j \in \mathcal{J}_C$.

The power generation technologies represented in the model are hydro, thermal, solar and wind. The set of existing hydropower plants in the initial configuration and the set of candidate hydropower plants for generation capacity expansion are denoted by $\mathcal{H}_E$ and $\mathcal{H}_C$, respectively. The following types of hydropower plants are represented in the model: (i) run-of-river hydropower plants without water storage, which operate as intermittent (non-programmable) energy sources, (ii) pumped-storage hydropower plants, and (iii) hydroelectric valleys, which are modelled as equivalent hydro-power plants with reservoir. The sets of programmable and non-programmable hydropower plants are denoted by $\mathcal{H}_P$ and $\mathcal{H}_{NP}$, respectively, with $\mathcal{H}_P \subseteq (\mathcal{H}_E \cup \mathcal{H}_C)$ and $\mathcal{H}_{NP} \subseteq (\mathcal{H}_E \cup \mathcal{H}_C)$.  
Constraints \eqref{eq:inv_hydro} and \eqref{eq:bin_hydro} model the investment decisions with respect to each candidate hydropower plant $h \in \mathcal{H}_C$. 

The thermoelectric generation system existing in zone $z$ at the beginning of the planning period is represented by the values of the integer parameters $N_{k,0}, k \in \mathcal{K}_z$, where $\mathcal{K}_z$ is the set of thermal power plant types in zone $z \in \mathcal{Z}$. Plant types $k_1$ and $k_2$, with $k_1,k_2 \in \mathcal{K}_z$, differ in their technical parameters (e.g., minimum power output, capacity, minimum uptime and downtime, heat rate). If $N_{k,0}=0$, thermal power plants of type $k$ are candidates for investment in thermal power generation, but are not present in the initial generation mix of zone $z$. The number $N_{k,y}$ of units of type $k$ available for production in year $y$ is determined by constraint \eqref{eq:inv_the} as the number $N_{k,y-1}$ of units available in year $y-1$ plus the number $N_{k,y}^+$ of units built in year $y$ minus the number $N_{k,y}^-$ of units decommissioned in year $y$, where $N_{k,y}$, $N_{k,y}^+$ and $N_{k,y}^-$ are integer variables. In constraint \eqref{eq:bound_the}, the lower bound $\underline{N}_{k,y}$ and the upper bound $\overline{N}_{k,y}$ allow the imposition of mandatory decommissioning and investment\footnote{For example, suppose that four thermal power plants of type $k1$ and three thermal power plants of type $k2$ are in the generation mix of zone $z$ at the beginning of the planning period, i.e. $N_{k1,0}=4$ and $N_{k2,0}=3$, with $k1 \in K_z$ and $k2 \in K_z$. The mandatory decommissioning of one plant of type $k1$ within the first 5 years of the planning period is modelled by assigning $\overline{N}_{k1,y}=4, 1 \leq y \leq 4$ and $\overline{N}_{k1,5}=3$, so that constraint \eqref{eq:bound_the} ensures that at least one decommissioning takes place within year 5. In a similar way, if in the period between year 2 and year 5, one new plant of type $k2$ has to be mandatorily built and two new plants of type $k2$ can be optionally built, then the following assignment have to be made: $\underline{N}_{k2,y}=3, 1 \leq y \leq 4$, $\underline{N}_{k2,5}=4$, $\overline{N}_{k2,1}=3$, and $\overline{N}_{k2,y}=6, 2 \leq y \leq 5$.}. 

Constraints \eqref{eq:bound_sol}$-$\eqref{eq:pos_sw} relate to the expansion of non-programmable renewable generation capacity. Given the solar power capacity $S_{z,0}$ existing in zone $z$ at the beginning of the planning period, constraint \eqref{eq:bound_sol} determines the solar power capacity $S_{z,y}$ to be installed in zone $z$ in year $y$, such that the total available capacity, resulting from all investments up to year $y$, lies between the lower bound $\underline{S}_{z,y}$ and the upper bound $\overline{S}_{z,y}$. Similarly, given the wind power capacity $W_{z,0}$ existing in zone $z$ at the beginning of the planning period, constraint \eqref{eq:bound_wind} determines the new capacity $W_{z,y}$ to be installed in zone $z$ in year $y$, such that the total available capacity, resulting from all investments up to year $y$, lies between the lower bound $\underline{W}_{z,y}$ and the upper bound $\overline{W}_{z,y}$. Parameters $\underline{S}_{z,y}, \overline{S}_{z,y}, \underline{W}_{z,y}, \overline{W}_{z,y}$ may be defined by policy targets.

Further flexibility can be provided to the system by batteries and PtG plants.
%In addition to pumped storage hydro plants, both batteries and PtG plants can contribute to increasing the share of renewable generation. 
%In addition to pumped storage in hydropower plants, other storage systems, such as batteries and PtG plants, can contribute to increasing the share of RES generation. 
%
The storage capacity of batteries of technology $b$ available in zone $z$ at the beginning of the planning period is represented by the real-valued parameters $B_{b,0}^{\text{CAP}},b \in \mathcal{B}_z$, where $\mathcal{B}_z$ is the set of battery technologies either existing or candidate in zone $z$. 
Constraints \eqref{eq:pos_bat} and \eqref{eq:inv_bat}  determine the new capacity $B_{b,y}^{\text{CAP}}$ of batteries of technology $b \in \mathcal{B}_z$ to be installed in zone $z$ in year $y$, with the total available capacity at the end of the planning period bounded above by $\overline{B}_{b}^{\text{CAP}}$. 
Similarly, the capacity of PtG plants of technology $g$ available in zone $z$ at the beginning of the planning period is represented by the real-valued parameters ${PtG}_{g,0}^{\text{CAP}}, g \in \mathcal{G}_z$, where $\mathcal{G}_z$ is the set of PtG technologies either existing or candidate in zone $z$.
% SHORT TERM
Constraints \eqref{eq:pos_ptg} and \eqref{eq:inv_ptg} determine the new capacity $PtG_{g,y}^{\text{CAP}}$ of PtG plants of technology $g \in \mathcal{G}_z$ to be installed in zone $z$ in year $y$, with the total available capacity at the end of the planning period bounded above by $\overline{PtG}_{g}^{\text{CAP}}$.

Constraints \eqref{eq:RES}$-$\eqref{eq:gas_LT_pos} model the short-term system operation.
In the presence of large shares of renewable generation, the system operation needs to be modelled with hourly discretisation in order to capture the intermittency of non-programmable renewable generation and to provide an accurate estimate of thermal generation costs and storage dynamics. 
Since considering all days of the planning period with hourly discretisation would be computationally infeasible, in our approach we obtain an accurate representation of short-term operation by considering a small set of representative days selected from historical data by a clustering procedure. Representative days allow the short-term variability of primary energy sources (wind speed, solar radiation and water inflows), electric load and gas demand in each system zone to be included in the GTEP analysis.
For more details on the selection of representative days, we refer the reader to \cite{micheli2021selecting}.
We denote with $c \in \mathcal{C}^y$ the set of representative days in year $y \in \mathcal{Y}$. 
Each representative day is discretised in hours, with $\mathcal{T}:=\{1,2,\ldots,24\}$ being the set of the hours of the day, and is characterized by specific hourly profiles of: 
(i) wind power production coefficients $\rho_{z,t,c}$, 
(ii) solar power production coefficients $\mu_{z,t,c}$, 
(iii) natural inflows for hydro reservoirs $F_{h,t,c,y}$, 
(iv) electricity demand $D^{\text{P}}_{z,t,c,y}$, and 
(v) natural gas demand $D^{\text{G}}_{n,t,c,y}$. 
Furthermore, in each year $y \in \mathcal{Y}$, 
$\mathcal{D}^y$ denotes the corresponding set of calendar days, and 
$c(d)$ is the injective mapping that determines the cluster (i.e., the representative day) $c \in \mathcal{C}^y$ 
to which the calendar day $d \in \mathcal{D}^y$ belongs. 
The $c(d)$ mapping is used in the mass balance constraints of hydro and gas reservoirs, which capture the medium-term seasonality.

Constraints \eqref{eq:RES}$-$\eqref{eq:gas_LT_pos} model the system operation in each year $y \in \mathcal{Y}$ and each scenario $w \in \mathcal{W}$.
In particular: 
(i) the operation of the integrated electricity and gas systems at each hour $t$ of each representative day $c \in \mathcal{C}^y$ is defined by constraints \eqref{eq:RES}$-$\eqref{eq:hydro_SPILL_new}, \eqref{eq:the_max_ON}$-$\eqref{eq:pos_reserve} and \eqref{eq:bat_level}$-$\eqref{eq:gas_OUT}; 
(ii) constraints \eqref{eq:hydro_LT_0}$-$\eqref{eq:hydro_LT_pos} control the reservoirs of the programmable hydropower plants; 
(iii) constraints \eqref{eq:gas_LT_0}$-$\eqref{eq:gas_LT_pos}  control the gas storage; 
(iv) constraints \eqref{eq:penetration} and \eqref{eq:CO2} impose policy targets on renewables penetration and on CO$_2$ emissions.

% VINCOLO 18
%Constraints \eqref{eq:RES}$-$\eqref{eq:CO2} pertain to the electricity system. 
Equation \eqref{eq:RES} defines the non-programmable renewable generation $RES_{z,t,c,y}$ in zone $z$ at hour $t$ of representative day $c$ of year $y$ as 
the solar power capacity available in year $y$ multiplied by the solar power production coefficient $\mu_{z,t,c}$ in zone $z$ at hour $t$ of representative day $c$ 
plus the wind power capacity available in year $y$ multiplied by the wind power production coefficient $\rho_{z,t,c}$ in zone $z$ at hour $t$ of representative day $c$. 

% VINCOLO 19
Constraints \eqref{eq:balance} enforce the energy balance in each zone $z$ at each hour $t$:
the electricity available is the sum of 
(i) the output $RES_{z,t,c,y}$ of solar and wind power plants,
(ii) the output $\sum_{k \in \mathcal{K}_z} (\underline{P}_k \ \gamma_{k,t,c,y,w} + p_{k,t,c,y,w})$ of thermal power plants,
(iii) the output $\sum_{h \in \mathcal{H}_z} H_{h,t,c,y,w}^{\text{OUT}}$ of hydro power plants, 
(iv) the battery discharge $\sum_{b \in \mathcal{B}_z} B_{b,t,c,y,w}^{\text{OUT}}$, and 
(v) the import $\sum_{l \in BS^{\text{P}}_z} F^\text{L}_{l,t,c,y,w}$ from other zones; 
the electricity used is the sum of 
(i) the zonal load $D^{\text{P}}_{z,t,c,y}$,
(ii) the electricity $\sum_{g \in \mathcal{G}_z} \frac{G^{\text{PtG}}_{g,t,c,y,w}}{\eta_g}$ used by PtG plants,
(iii) the electricity $\sum_{h \in \mathcal{H}_z} H_{h,t,c,y,w}^{\text{IN}}$ stored in hydro reservoirs,
(iv) the battery charge $\sum_{b \in \mathcal{B}_z} B_{b,t,c,y,w}^{\text{IN}}$, and
(v) the export $\sum_{l \in FS^{\text{P}}_z} F^\text{L}_{l,t,c,y,w}$ to other zones. 
Possible imbalances are taken into account either by the nonnegative variable $E^{\text{NP}}_{z,t,c,y,w}$ \eqref{eq:non_neg_ENP}, which represents the energy not supplied when consumption exceeds the available energy, or by the nonnegative variable $OG_{z,t,c,y,w}$ \eqref{eq:non_neg_OG}, which represents the overgeneration when the available energy exceeds consumption.

% VINCOLI 20 e 21
Constraints \eqref{eq:flow_el} and \eqref{eq:flow_el_NEW} enforce lower and upper bounds on power flows at each hour $t$ through existing transmission lines $l \in \mathcal{L}_E$ and through candidate transmission lines $l \in \mathcal{L}_C$ that have been built by year $y$ (i.e., $\theta_{l,y}=1$). Constraints \eqref{eq:flow_el_NEW} enforce zero flows at hour $t$ through candidate transmission lines $l \in \mathcal{L}_C$ that have not been built (i.e., $\theta_{l,y}=0$).

%%%%%%%%%%%%%%%%%%%%%%%%%%%%%%%%% VINCOLI da 22 a 32 %%%%%%%%%%%%%%%%%%%%%%%%%%%%%%%%%
Constraints \eqref{eq:hydro_OUT}$-$\eqref{eq:hydro_LT_pos} model the operation of the hydropower subsystem.  
% VINCOLI 22 e 23
Constraints \eqref{eq:hydro_OUT} and \eqref{eq:hydro_OUT_new} apply to both programmable and non-programmable hydro power plants. They 
%enforce the upper bound $\overline{H}_h^\text{OUT}$ on the power output of both existing plants $h \in \mathcal{H}_E$ and candidate plants $l \in \mathcal{H}_C$ that have been built by year $y$. Zero power output is enforced by constraints \eqref{eq:hydro_OUT_new} for candidate plants $l \in \mathcal{H}_C$ that have not been built. Constraints \eqref{eq:hydro_OUT} and \eqref{eq:hydro_OUT_new} 
state that the hourly power output $H_{h,t,c,y,w}^\text{OUT}, 1 \leq t \leq 24, c \in \mathcal{C}^y, y \in \mathcal{Y}, w \in \mathcal{W}$, of any hydropower plant available in year $y$ is bounded above by the maximum output $\overline{H}_h^\text{OUT}$, while power output is zero for candidate hydropower plants not available in year $y$. 
%%%%
% VINCOLI 24
The output of non-programmable run-of-river plants is also bounded above by the water flow $F_{h,t,c,y}$ in hour $t$, as stated by constraints \eqref{eq:hydro_OUT_NP}.
Programmable plants are modelled as equivalent plants consisting of one reservoir and one turbine. Three types of programmable plants are considered: (a) plants with natural inflows and no pumped storage; (b) plants with pumped storage without natural inflows; (c) plants with pumped storage and natural inflows.  
% VINCOLI 25 e 26
Constraints \eqref{eq:hydro_IN} and \eqref{eq:hydro_IN_new} impose that the hourly pumping power $H_{h,t,c,y,w}^\text{IN}$ of any hydropower plant available in year $y$ is bounded above by the maximum pumping power $\overline{H}_h^\text{IN}$, being zero for candidate hydropower plants not available in year $y$. 
% VINCOLI 27 e 28
Similarly, constraints \eqref{eq:hydro_SPILL} and \eqref{eq:hydro_SPILL_new} impose bounds to the spillage $H_{h,t,c,y,w}^\text{SPILL}$.
In order to limit the number of variables and constraints, the reservoir levels of programmable plants $h$ available in year $y$ are not checked at the end of every hour, but only every $M$ days: therefore, in each year $y$ the level of any reservoir is checked $M \cdot \xi$ times, $1 \leq \xi \leq \overline{\xi}$, where  $\overline{\xi}=\lfloor \frac{365}{M} \rfloor$ is the maximum integer preceding $\frac{365}{M}$. 
% VINCOLI 29 e 30
For each year $y$ and each scenario $w$, constraint \eqref{eq:hydro_LT_0} performs the first check in the year, i.e., determines the reservoir level at the end of day $M$, while constraints \eqref{eq:hydro_LT_M} determine the reservoir level at the end of days $M \cdot \xi$, for $2 \leq \xi \leq \overline{\xi}$. 
% VINCOLI 31
Since hydroelectric reservoirs are storage reservoirs with an annual cycle, constraint \eqref{eq:hydro_LT_final} states that the level of the reservoir at the end of the year must be equal to the level $H_{h,0}^\text{LT}$ of the reservoir at the beginning of the year. 
% VINCOLI 32
Finally, constraints \eqref{eq:hydro_LT_pos} impose that the reservoirs levels $H_{h,\xi M,y,w}^\text{LT}, 1 \leq \xi \leq \overline{\xi}$, do not exceed the energy level $EPR_h \ \overline{H}_h^\text{OUT}$, where $EPR_h$ denotes the energy-to-power ratio of reservoir $h$. 
%\textcolor{red}{Since only representative days are considered in the operational problem, in constraints \eqref{eq:hydro_LT_0}, \eqref{eq:hydro_LT_M} and \eqref{eq:hydro_LT_final} each original day $d \in \mathcal{D}^y$ of year $y$ is replaced by the correspondent representative day $c(d) \in \mathcal{C}^y$.}

% VINCOLI da 33 a 38
%For $k \in \mathcal{K}_z, z \in \mathcal{Z}, y \in \mathcal{Y}, w \in \mathcal{W}$ 
Constraints \eqref{eq:the_max_ON}$-$\eqref{eq:the_int} model the daily operation of thermal power plants of type $k \in \mathcal{K}_z$ in zone $z \in \mathcal{Z}$, on representative day $c$ of year $y$ under scenario $w$. 
% VINCOLO 33
Constraints \eqref{eq:the_max_ON} state that at each hour $t$ the number $\gamma_{k,t,c,y,w}$ of online units is bounded above by the number of units available in the corresponding year $y$. 
% VINCOLO 34
Constraints \eqref{eq:the_status} state that the number of online units at hour $t$ equals the number of online units at hour $t-1$ plus the number of units started up at hour $t$ minus the number of units shut down at hour $t$. 
% VINCOLO 35
Any unit can be started-up at most once in any interval of $MUT_k$ consecutive hours, as stated by minimum uptime constraint \eqref{eq:the_MUT}. 
% VINCOLO 36
Any unit can be shut down at most once in any interval of $MDT_k$ consecutive hours, as imposed by minimum downtime constraint \eqref{eq:the_MDT}. 
Since the representative days $c \in \mathcal{C}^y$ are disconnected from each other, the minimum uptime constraints are not enforced for hours 1 to $MUT_k-1$ and the minimum downtime constraints are not enforced for hours 1 to $MDT_k-1$. 
% VINCOLO 37
Constraint \eqref{eq:the_max_power} states that the total power output of all online units of type $k$ in hour $t$ must be between $\underline{P}_k \ \gamma_{k,t,c,y,w}$ and $\overline{P}_k \ \gamma_{k,t,c,y,w}$  with $p_{k,t,c,y,w}=0$ if no units of type $k$ are online at hour $t$.
% VINCOLI 39 e 40
Constraints \eqref{eq:reserve} and \eqref{eq:pos_reserve} determine whether or not the unused capacity of online thermal power plants provides the spinning reserve $R_{z,t,c,y}$ requested in zone $z$ at hour $t$. Positive values of $R_{z,t,c,y,w}^\text{NP}$, the reserve not provided, are penalized in the objective function.

% VINCOLI 41 e 42
Constraints \eqref{eq:penetration} and \eqref{eq:CO2} impose policy goals that must be met by aggregates of zones $a \in \mathcal{A}$ (e.g., nations). Constraints \eqref{eq:penetration} represent targets on renewable penetration: for each $a \in \mathcal{A}$ the total non-programmable renewable generation in year $y$ must cover at least ratio $\phi_{a,y}$ of the total yearly load. Constraints \eqref{eq:CO2} limit the total CO$_2$ emissions in each zone aggregate $a$ and each year $y$ under each scenario $w$.

% VINCOLI 43 e 47
Constraints \eqref{eq:bat_level}$-$\eqref{eq:bat_final} model the daily operation of batteries of type $b \in \mathcal{B}_z$ in zone $z \in \mathcal{Z}$ on representative day $c \in \mathcal{C}^y$ of year $y \in \mathcal{Y}$ under scenario $w \in \mathcal{W}$. The balance equation \eqref{eq:bat_level} defines the electricity stored at the end of hour $t$ as the energy stored at the end of the previous hour (reduced by the loss coefficient $\lambda_b$, 
$0 \leq \lambda_b \leq 1$), plus the charge (reduced by the loss coefficient $\lambda_b^\text{IN} \leq 1$), minus the discharge (reduced by the coefficient $\lambda_b^\text{OUT} \geq 1$). In each hour $t$ the energy level is bounded above by the product of the available capacity times the energy-to-power ratio $EPR_b$, as stated by constraint \eqref{eq:bat_inst}, while charge and discharge are bounded above by the capacity available in the corresponding year, see constraints \eqref{eq:bat_IN} and \eqref{eq:bat_OUT} respectively. 
Constraint \eqref{eq:bat_final} imposes a 24-hour cycle, i.e. the energy level $B_{b,24,c,y,w}$ at the end of the day equals the energy level $B_{b,0,c,y,w}$ at the beginning of the day.

% VINCOLO 48
Constraints \eqref{eq:gas_balance}$-$\eqref{eq:gas_LT_pos} describe the short-term operation of the natural gas system.
Constraints \eqref{eq:gas_balance}$-$\eqref{eq:gas_OUT} are imposed for all hours of each representative day $c \in \mathcal{C}^y$ of each year $y \in \mathcal{Y}$ under each scenario $w \in \mathcal{W}$.
The zonal gas balance equation \eqref{eq:gas_balance} states that at any hour $t$ the amount of gas available must be equal to the amount of gas required. 
In particular, the amount of gas available in zone $n$ is the sum of 
(i) the zonal gas supply $G_{n,t,c,y,w}$;
(ii) the amount $G_{n,t,c,y,w}^\text{OUT}$ of gas withdrawn from the zonal storage;
(iii) the amount $\sum_{g \in \mathcal{G}_n} G_{g,t,c,y,w}^\text{PtG}$ of synthetic biomethane produced by PtG plants located in the zone; 
(iv) the amount $\sum_{j \in BS^{\text{G}}_n} F^\text{J}_{j,t,c,y,w}$ of gas imported from other zones. 
The amount of gas required in zone $n$ is the sum of 
(i) the exogenous industrial and tertiary (i.e., residential and office) gas demand $D^{\text{G}}_{n,t,c,y}$;
(ii) the endogenous gas demand $\sum_{k \in \mathcal{K}^\text{G}_n} HR_k (\underline{P}_k \gamma_{k,t,c,y,w} + p_{k,t,c,y,w})$ for thermal power generation, where $HR_k$ is the heat rate of thermal power plants of type $k$;
(iii) the quantity $G_{n,t,c,y,w}^\text{IN}$ of gas injected into the zonal storage; 
(iv) the amount $\sum_{j \in FS^{\text{G}}_z} F^J_{j,t,c,y,w}$ of gas exported to other zones. 
The inability to meet the gas demand is taken into account by the nonnegative variable $G^{\text{CURT}}_{n,t,c,y,w}$ \eqref{non_neg_GCURT}, which represents the gas demand curtailment.
% VINCOLI 49 e 50
Constraints \eqref{eq:flow_gas} and \eqref{eq:flow_gas_NEW} enforce lower and upper bounds on gas flows at each hour $t$ through existing pipelines $j \in \mathcal{J}_E$ and through candidate pipelines $j \in \mathcal{J}_C$ that have been built by year $y$ (i.e., $\theta_{j,y}=1$). Constraints \eqref{eq:flow_gas_NEW} enforce zero flows at hour $t$ through candidate pipelines $j \in \mathcal{J}_C$ that have not been built (i.e., $\theta_{j,y}=0$).
% VINCOLI 51 e 52
The hourly production of synthetic biomethane from PtG plants of technology $g \in \mathcal{G}_n$ in zone $n$ is bounded above by the total installed capacity of PtG plants in year $y$, as imposed by constraints \eqref{eq:PtG_prod}, and the hourly gas supply $G_{n,t,c,y,w}$ in each zone $n$ is bounded by lower and upper supply limits, as imposed by constraints \eqref{eq:gas_prod}.
% VINCOLI 53 e 54
All gas storages in zone $n$ are represented by a single equivalent storage associated with the zone.
Constraint \eqref{eq:gas_IN} requires that the gas injected into the storage at hour $t$ does not exceed the maximum hourly injection rate $\overline{G}^\text{IN}_{n}$; similarly, constraint \eqref{eq:gas_OUT} requires that the gas withdrawn from the storage at hour $t$ does not exceed the maximum hourly withdrawal rate $\overline{G}^\text{OUT}_{n}$.

% VINCOLI 55 e 58
In each year and scenario, the capacity constraint check for each gas storage is every $M$ days rather than at the end of every hour.
Constraint \eqref{eq:gas_LT_0} determines the storage content at the end of day $M$ and constraints \eqref{eq:gas_LT_M} determine the gas storage level at the end of days $M \cdot \xi$, for $2 \leq \xi \leq \overline{\xi}$. 
% VINCOLI 31
Gas storages have an annual cycle, therefore constraint \eqref{eq:gas_LT_final} states that the storage level at the end of the year must be equal to the level $G^\text{LT}_{n,0}$ at the beginning of the year. 
% VINCOLI 32
Finally, constraints \eqref{eq:gas_LT_pos} impose that the gas storage levels 
$G^\text{LT}_{n,\xi M,y,w}, 1 \leq \xi \leq \overline{\xi}$, do not exceed the capacity $\overline{G}^\text{LT}_n$ of gas storage $n$.

%\color{red}
The objective, expressed in \eqref{eq:obj}, is to minimise the sum of the following two terms: (i) total investment and decommissioning costs over the planning period, discounting costs paid in year $y$ to reference year $y_0$ using discount rate $r$; (ii) the expected cost of operating the electricity and gas systems over the planning period, where the uncertainty in fuel and CO$_2$ prices is represented by the set $\mathcal{W}$ of scenarios.
In (ii), the cost at hour $t$ of representative day $c$ of year $y$ in scenario $w$ is the sum of the following terms:
\begin{enumerate}
    \item cost of hydro power production $\sum_{z \in \mathcal{Z}} \sum_{h \in \mathcal{H}_z} C^\text{H}_h \ H^{\text{OUT}}_{h,t,c,y,w}$;
    \item start-up costs of thermal plants $\sum_{z \in \mathcal{Z}} \sum_{k \in \mathcal{K}_z} C_k^{\text{SU}}\ \alpha_{k,t,c,y,w}$;
    \item operating costs of thermal plants $\sum_{z \in \mathcal{Z}} \sum_{k \in \mathcal{K}_z}  C^\text{M}_{k,y,w} \left( \underline{P}_k \gamma_{k,t,c,y,w}  +  p_{k,t,c,y,w} \right)$;
    \item cost of electricity from batteries $\sum_{z \in \mathcal{Z}} \sum_{b \in \mathcal{B}_z}  C^{B}_b  \ B^{\text{OUT}}_{b,t,c,y,w}$;
    \item penalties for overgeneration $\sum_{z \in \mathcal{Z}} C^{\text{OG}} \ OG_{z,t,c,y,w}$;
    \item penalties for not supplied electricity $\sum_{z \in \mathcal{Z}} C^{\text{ENP}} E^{\text{NP}}_{z,t,c,y,w}$;
    \item penalties for not supplied reserve $\sum_{z \in \mathcal{Z}} C^{\text{RNP}} R^{\text{NP}}_{z,t,c,y,w}$;
    \item cost of natural gas supply $\sum_{n \in \mathcal{N}} C^{\text{G}}_{n,y,w} \ G_{n,t,c,y,w}$;
    \item operating costs of PtG plants $\sum_{n \in \mathcal{N}} \sum_{g \in \mathcal{G}_n} \ C^\text{PtG}_g \ G^{\text{PtG}}_{g,t,c,y,w}$;
    \item penalties for gas demand curtailment $\sum_{n \in \mathcal{N}} C^{\text{GC}} \ G^{\text{CURT}}_{n,t,c,y,w}$.
\end{enumerate}
Since the sum of terms 5 and 6 expresses the cost of the gas used to satisfy both the exogenous gas demand, $D^\text{G}_{n,t,c,y}$, and the endogenous gas demand of gas-fired thermal power plants, $\sum_{k \in \mathcal{K}^\text{G}_n} HR_k \ 
(\underline{P}_k \ \gamma_{k,t,c,y,w} + p_{k,t,c,y,w})$, 
%expressed by the term  
% in \eqref{eq:gas_balance}
the marginal costs $C^\text{M}_{k,y,w}$ of power plants $k \in \mathcal{K}^\text{G}_n$ in term 4 are defined as
\begin{align}
& C^\text{M}_{k,y,w} = OM_k + {CO_2}_k \ Pr^\text{CO$_2$}_{y,w}.
\label{eq:cost_k1}
\end{align}
In \eqref{eq:cost_k1}, $OM_k$ and ${CO_2}_k$ are the operating and maintenance costs and the $CO_2$ emission rate, respectively, 
of thermal power plants of type $k$,   
and $Pr^\text{CO$_2$}_{y,w}$ is the price of $CO_2$ emissions.
For non-gas-fired thermal power plants, the marginal cost is defined as
\begin{align}
& C^\text{M}_{k,y,w} 
= OM_k 
+ {CO_2}_k \ Pr^\text{CO$_2$}_{y,w}
+ HR_k \ Pr^\text{FUEL}_{f(k),y,w},
\label{eq:cost_k2}
\end{align}
where $Pr^\text{FUEL}_{f(k),y,w}$ is the price of fuel $f \in \mathcal{F}$ used by plants of type $k$ and $\mathcal{F}$ is the set of all fuels except gas (e.g. oil, coal).

\section{The solution algorithm}
%  based on Benders decomposition
\label{sec:Algorithm}

Solving the proposed two-stage stochastic programming model for the integrated planning of electricity and natural gas systems poses significant computational challenges due to the high level of temporal and technical detail over a long-term planning horizon. 
To address the computational burden, we propose a multi-cut Benders decomposition algorithm, which is obtained by suitably modifying the algorithm developed in \cite{micheli2020two} for planning the expansion of electricity systems.
Benders decomposition algorithms have been widely used in power system planning \cite{RAHMANIANI2017801}, with notable applications in large-scale transmission expansion planning \cite{LUMBRERAS2013transmission}, economic dispatch of combined heat and power \cite{ABDOLMOHAMMADI2013benders}, strategic planning of energy hubs \cite{mansouri2020stochastic}, day-ahead scheduling of energy communities operating under peer-to-peer energy trading schemes\cite{garcia2023benders}, and operational planning for aggregators in smart grids \cite{soares2017two}.
To describe the proposed solution algorithm, we write the objective function \eqref{eq:obj} as follows

%\begin{align*}
%    \min \sum_{y\in\mathcal{Y}} \Big( \mathbf{CI}_y^\top \, \mathbf{x}_y + \sum_{w\in\mathcal{W}} pb_w  \
%    \mathbf{CO}^\top_{y,w} \ \mathbf{s}_{y,w}  \Big),
%\end{align*}
\begin{align*}
    \min \sum_{y\in\mathcal{Y}} \Big( \mathbf{CI}_y^{\top} \, \mathbf{x}_y + \sum_{w\in\mathcal{W}} pb_w  \
    \mathbf{CO}^{\top}_{y,w} \ \mathbf{s}_{y,w}  \Big).
\end{align*}
The term 
\begin{align*}
\mathbf{CI}_y^{\top}\mathbf{x}_{y} =
&
\frac{1}{(1+r)^{y-y_0}} \Bigg\{ 
\sum_{l \in \mathcal{L}_C} IC_l^\text{L} \  \delta_{l,y} + 
\sum_{j \in \mathcal{J}_C} IC_j^\text{J} \ \delta_{j,y} + 
\sum_{h \in \mathcal{H}_C} IC_h^\text{H} \ \overline{H}_h^{\text{OUT}} \ \delta_{h,y} + \nonumber \\
% 2
& + \sum_{z \in \mathcal{Z}} \Bigg[ 
\sum_{k \in \mathcal{K}_z} \big( IC_{k,y}^\text{K} \ \overline{P}_k \ N_{k,y}^{+} + DC_{k,y}^\text{K} \ \overline{P}_k \ N_{k,y}^{-} \big) + IC_{z,y}^\text{S} \ S_{z,y} \ + \nonumber \\
% 3
& + IC_{z,y}^\text{W} \ W_{z,y}  +  \sum_{b \in \mathcal{B}_z} IC_{b,y}^\text{B} \ B^\text{CAP}_{b,y}  +  \sum_{g \in \mathcal{G}_z} IC_{g}^\text{PtG} \ PtG^\text{CAP}_{g,y} \Bigg] \Bigg\}
\end{align*}    
represents the discounted investment cost at year $y$, where vector $\mathbf{x}_y$ is
%and $\mathbf{s}_{y,w}$ are 
defined as 
\begin{alignat*}{2}
\mathbf{x}_y = \big[ 
& RES_{z,t,c,y}, \theta_{l,y}, \theta_{h,y}, \theta_{j,y}, N_{k,y}, B^{\text{CAP}}_{b,y}, PtG^{\text{CAP}}_{g,y},  && \\    
& S_{z,y}, W_{z,y}, \delta_{l,y}, \delta_{h,y}, \delta_{j,y}, N_{k,y}^+, N_{k,y}^- \big], \qquad && y \in \mathcal{Y},  
\end{alignat*} 
and $\mathbf{CI}_y$ 
%and $\mathbf{CO}_{y,w}$ 
is the vector of coefficients of $\mathbf{x}_y$.
%of coefficients of components of vectors $\mathbf{x}_{y}$ and $\mathbf{s}_{y,w}$, respectively.
The term
\begin{align*}
\mathbf{CO}_{y,w}^{\top} \, \mathbf{s}_{y,w} =
& \sum_{c \in \mathcal{C}^y} \psi_c  \sum_{t=1}^{24} \Bigg[ \sum_{z \in \mathcal{Z}}  \Bigg(
\sum_{h \in \mathcal{H}_z}  C_h^\text{H}  \ H_{h,t,c,y,w}^{\text{OUT}} + \sum_{k \in \mathcal{K}_z} C_k^{\text{SU}}\ \alpha_{k,t,c,y,w} + \nonumber \\
\nonumber \\
% 5
& + \sum_{k \in \mathcal{K}_z} C_{k,y,w}^\text{M} \left( \! \underline{P}_k \gamma_{k,t,c,y,w}  +  p_{k,t,c,y,w} \right) + \sum_{b \in \mathcal{B}_z}  C_b^{B}  \ B_{b,t,c,y,w}^{\text{OUT}} + \nonumber \\
% 6
& + C^{\text{OG}} \ OG_{z,t,c,y,w} + C^{\text{ENP}} E^{\text{NP}}_{z,t,c,y,w} + C^{\text{RNP}} R^{\text{NP}}_{z,t,c,y,w} \Bigg) + \nonumber \\
% 7
&  +   \sum_{n \in \mathcal{N}}  \Bigg(  
C^{\text{G}}_{n,y,w} G_{n,t,c,y,w}  +
 \sum_{g \in \mathcal{G}_n} C^\text{PtG}_g G^{\text{PtG}}_{g,t,c,y,w}+ 
 C^{\text{GC}}  G^{\text{CURT}}_{n,t,c,y,w} 
\Bigg) \Bigg]
\end{align*} 
represents the operational cost at year $y$ in scenario $w$, where vector $\mathbf{s}_{y,w}$ is 
defined as 
\begin{alignat*}{2}
\mathbf{s}_{y,w} = \big[ &    
\alpha_{k,t,c,y,w},
\beta_{k,t,c,y,w},
\gamma_{k,t,c,y,w},
p_{k,t,c,y,w},
H_{h,t,c,y,w}^{\text{IN}},
H_{h,t,c,y,w}^{\text{OUT}},
&& \quad \\
& 
H_{h,t,c,y,w}^{\text{SPILL}},
B_{b,t,c,y,w},
B_{b,t,c,y,w}^\text{IN},
B_{b,t,c,y,w}^\text{OUT},
F_{l,t,c,y,w}^\text{L},
OG_{z,t,c,y,w},
&& \\
& 
E^\text{NP}_{z,t,c,y,w},
R^\text{NP}_{z,t,c,y,w},
G_{g,t,c,y,w}^{\text{PtG}},
G_{n,t,c,y,w},
G^\text{IN}_{n,t,c,y,w},
G^\text{OUT}_{n,t,c,y,w},
 \\
&
F_{j,t,c,y,w}^\text{J},
G^\text{CURT}_{n,t,c,y,w}
\big] && y \in \mathcal{Y}, \ w \in \mathcal{W}, 
\end{alignat*}
and $\mathbf{CO}_{y,w}$ is the vector of coefficients of $\mathbf{s}_{y,w}$. 
Furthermore, let $\underline{GTEP}$ be the MILP problem obtained from the GTEP problem by defining the unit commitment variables $\alpha_{k,t,c,y,w}, \ \beta_{k,t,c,y,w}, \ \gamma_{k,t,c,y,w}$ as real nonnegative variables.
The algorithm is as follows:

\medskip
\noindent
\textbf{Initialization}: assign the convergence tolerance $\varepsilon$ and the maximum number of iterations $i^\text{max}$. 
%    Set $z_\text{LB}^{(0)}=0$ and $z_\text{UP}^{(0)}=\infty$.
    \\

    \noindent
    \textbf{Iterative cycle}: for $i = 1, \ldots, i^\text{max}$
 \medskip
    
    \begin{itemize}
    \item[A)]
        Solve the MILP problem $GTEP_\text{LB}^{(i)}$ defined as
  \begin{align}
     \min_{\substack{\mathbf{x}_{y},y \in \mathcal{Y} \\ 
     \theta_w,w \in \mathcal{W}}} z_\text{LB} = \ 
     & \sum_{y\in\mathcal{Y}} \mathbf{CI}_{y}^\top \mathbf{x}_{y} + \sum_{w \in \mathcal{W}} pb_w \theta_w
    \label{eq:obj_MP} \\
    \text{s.t. } 
    & \eqref{eq:inv_line}-\eqref{eq:RES}, \eqref{eq:penetration}, 
    \label{eq:prec_MP} \\
    & \theta_w \geq \sum_{y \in \mathcal{Y}} \Big[ z_{y,w}^{(\nu)} + \Big( \mathbf{x}_{y} - \mathbf{x}_{y}^{(\nu)}
    \Big)^\top \mathbf{\lambda}_{y,w}^{(\nu)} \; \Big], \quad 1\leq\nu\leq i-1, 
    \nonumber \\
    & \qquad \quad \quad \,\,  w \in\mathcal{W}, \quad \text{if $i>1$,} 
    \label{eq:cut} \\
    &  \theta_w = 0, \quad w \in \mathcal{W}, \quad \text{if $i=1$}. 
    \label{eq:min_theta} 
  \end{align}
  \item[] Set $\mathbf{x}_{y}^{(i)}=\mathbf{x}_{y}^{*}, \ y \in \mathcal{Y}$, and $z_\text{LB}^{(i)}=z_\text{LB}^{*}$.
  
 \medskip
 \item[B)]
  For each year $y \in \mathcal{Y}$ in each scenario $w\in\mathcal{W}$, solve the LP problem $OPE_{y,w}^{(i)}$ defined as
  \begin{align}
    \min_{\mathbf{x}_{y},\mathbf{s}_{y,w}}  \ z_{y,w} = \ & 
    \mathbf{CO}_{y,w}^\top \mathbf{s}_{y,w} 
    \label{eq:obj_SP} \\
    \text{s.t. } 
    & \eqref{eq:balance}-\eqref{eq:the_max_power}, \eqref{eq:reserve}, \eqref{eq:pos_reserve}, \eqref{eq:CO2}-\eqref{eq:gas_LT_pos}, 
    \label{eq:prec_SP} \\
    & \alpha_{k,t,c,y,w}, \ \beta_{k,t,c,y,w}, \ \gamma_{k,t,c,y,w} \geq 0, \quad k\in \mathcal{K}, \nonumber \\
    & \ 1 \leq t \leq 24, \ c\in\mathcal{C}^y,  
    \label{eq:linear} \\ 
    &  \mathbf{x}_y = \mathbf{x}_y^{(i)} \ \ \ \ \ : \ \ \ \ \ \mathbf{\lambda}_{y,w}.
    \label{eq:dual}
\end{align}  
   \item[] Set $z^{(i)}_{y,w}=z^{*}_{y,w}$, %\textcolor{blue}{$\mathbf{s}^{(i)}_{y,w}=\mathbf{s}^{*}_{y,w}$} 
   and $\mathbf{\lambda}^{(i)}_{y,w}=\mathbf{\lambda}^{*}_{y,w}$.
  \medskip
  \item[C)]
   Compute 
  \begin{equation}
   z_\text{GTEP}^{(i)} = \sum_{y\in\mathcal{Y}}  \Big( \mathbf{CI}_y^\top \mathbf{x}_y^{(i)} + \sum_{w\in\mathcal{W}} pb_w \ z_{y,w}^{(i)} \Big).
   \label{eq:z_i_GTEP}
  \end{equation}
    \item[] If $i=1$, set $z_\text{UB}^{(1)}=z_\text{GTEP}^{(1)}$ and $\mathbf{x}^\text{BEST}_y = \mathbf{x}_y^{(1)}, \ y\in\mathcal{Y}$,
 
    \item[] else if $z_\text{GTEP}^{(i)} < z_\text{UB}^{(i-1)}$,
    \item[] $\quad$ set $z_\text{UB}^{(i)} = z_\text{GTEP}^{(i)}$ and $\mathbf{x}^\text{BEST}_y = \mathbf{x}_y^{(i)}, \ y\in\mathcal{Y}$,

    \item[] else
    \item[] $\quad$ set $z_\text{UB}^{(i)} = z_\text{UB}^{(i-1)}$.

    \medskip
    \item[D)]
     If $\frac{z_\text{UB}^{(i) } - z_\text{LB}^{(i)}}{z_\text{UB}^{(i)}} > \varepsilon$,
    go to A;
    \item[] else
%    \begin{itemize}
     for each year $y \in \mathcal{Y}$ in each scenario $w\in\mathcal{W}$, solve the MILP problem $OPE_{y,w}$ and STOP. Problem $OPE_{y,w}$ is defined as follows
\begin{align}
    \min_{\mathbf{s}_{y,w}}  \ \ & 
    \mathbf{CO}_{y,w}^\top \mathbf{s}_{y,w} \label{eq:obj_SP_final}\\
    \text{s. t. } & \eqref{eq:balance}-\eqref{eq:pos_reserve}, 
    \eqref{eq:CO2}-\eqref{eq:gas_LT_pos}, \label{eq:prec_SP_final} \\
    &  \mathbf{x}_y = \mathbf{x}_y^{\text{BEST}}. \label{eq:fix_final}
\end{align}  

    \end{itemize}

\medskip

%    

%The procedure starts by assigning the maximum number of iterations $i^\text{max}$ and the tolerance $\varepsilon$ for the convergence test. The lower and upper bounds of the optimal value of the GTEP objective function are initialized to $z_\text{LB}^{(0)}=0$ and $z_\text{UB}^{(0)}=\infty$, respectively.
%
At step A of iteration $i$, the MILP problem $GTEP_{LB}^{(i)}$ is solved to determine the investment decisions $\mathbf{x}_{y}$ in each year $y \in \mathcal{Y}$ and the approximated operating cost $\theta_w$ in each scenario $w \in \mathcal{W}$ that minimise the objective function \eqref{eq:obj_MP}, i.e. the sum of the investment cost and the expected value of the approximated operating costs over the scenarios.
The constraints of $GTEP_{LB}^{(i)}$ are 
(i) the first stage constraints \eqref{eq:inv_line}-\eqref{eq:RES} and \eqref{eq:penetration} of the GTEP problem, and 
(ii) the definition \eqref{eq:cut} of the approximated operating cost $\theta_w$ in each scenario $w$. 
At iteration $i=1$, the approximated operating cost $\theta_w$ is zero in all scenarios. 
At iteration $i > 1$, the approximated operating cost $\theta_w$ is determined by using information computed at step B of previous iterations $\nu$, $1 \leq \nu \leq i-1$.
To this aim, in constraint \eqref{eq:cut} the approximated operating cost in year $y$ under scenario $w$ is expressed as the following linear function ${\theta}^{(\nu)}_{y,w} (\mathbf{x}_{y})$ of the investment variables $\mathbf{x}_{y}$   
\begin{equation}
    {\theta}^{(\nu)}_{y,w} (\mathbf{x}_{y}) = z^{(\nu)}_{y,w} + \Big(\mathbf{x}_{y} - \mathbf{x}_{y}^{(\nu)} \Big)^\top \mathbf{\lambda}_{y,w}^{(\nu)},
\label{approxtheta}
\end{equation}
where $z^{(\nu)}_{y,w}$, $\mathbf{x}^{(\nu)}_{y}$ and $\mathbf{\lambda}^{(\nu)}_{y,w}$ have been determined at a previous iteration $\nu < i$, and the approximated operating cost $\theta_w$ is computed as ${\theta}_{w} =\max_{1 \leq \nu \leq i-1} \{ \sum_{y\in\mathcal{Y}} {\theta}^{(\nu)}_{y,w} \}$.
%
%\textcolor{blue}{Since variable $\theta_w$ approximates by means of linear functions the true operating cost from below \cite{birge2011introduction}, the optimal value $z^{*}_{LB}$ of the objective function \eqref{eq:obj_MP} is a lower bound of the minimum value of the objective function \eqref{eq:obj} of the GTEP problem} and it is stored in $z^{(i)}_{LB}$. 
Based on \cite{birge2011introduction}, Chapter 5.1, the optimal value $z^{*}_{LB}$ of the objective function \eqref{eq:obj_MP} is a lower bound of the minimum value of the objective function of problem $\underline{GTEP}$ and it is stored in $z^{(i)}_{LB}$.
The optimal investments determined by $GTEP_{LB}^{(i)}$ are stored in $\mathbf{x}^{(i)}_y$, $y \in \mathcal{Y}$, to be used at step B of the current iteration and at step A of subsequent iterations $\nu > i$.

At step B of iteration $i$, for each year $y \in \mathcal{Y}$ under each scenario $w \in \mathcal{W}$ the LP problem $OPE^{(i)}_{y,w}$ is solved to determine the operational decisions $\mathbf{s}_{y,w}$ that minimize the operational cost \eqref{eq:obj_SP}. 
The constraints of $OPE^{(i)}_{y,w}$ are
(i) the second stage constraints \eqref{eq:balance}-\eqref{eq:the_max_power}, \eqref{eq:reserve}, \eqref{eq:pos_reserve} and \eqref{eq:CO2}-\eqref{eq:gas_LT_pos} of the GTEP problem (constraint \eqref{eq:prec_SP}),  
(ii) the relaxation of the integer unit commitment variables $\alpha_{k,t,c,y,w}$, $\beta_{k,t,c,y,w}$ and $\gamma_{k,t,c,y,w}$, which in $OPE^{(i)}_{y,w}$ are real nonnegative variables (constraint \eqref{eq:linear}), and 
(iii) the assignment of the optimal values $\mathbf{x}^{(i)}_{y}$, computed at step A, to the investment decisions $\mathbf{x}_{y}$ (constraint \eqref{eq:dual}), where $\mathbf{\lambda}^{(i)}_{y,w}$ denotes the vector of the dual variables of the assignment constraints. 
The optimal value $z^{*}_{y,w}$ of the objective function and 
the optimal values $\mathbf{\lambda}^{*}_{y,w}$ of the dual variables 
are stored in $z^{(i)}_{y,w}$ and $\mathbf{\lambda}^{(i)}_{y,w}$, respectively, 
to be used at step A of subsequent iterations $\nu > i$.

At step C of iteration $i$, the cost $z_\text{GTEP}^{(i)}$ computed by \eqref{eq:z_i_GTEP} is an upper bound of the optimal value of problem $\underline{GTEP}$, because the values of $\mathbf{x}_y^{(i)}$ and $z_{y,w}^{(i)}$ used in \eqref{eq:z_i_GTEP} have been obtained under constraints \eqref{eq:prec_MP}, \eqref{eq:prec_SP} and \eqref{eq:linear}. 
Since $z^{(\nu)}_{\text{GTEP}}$ may not necessarily decrease over iterations, we define $z^{(i)}_{\text{UB}}$ as 
$z^{(i)}_{\text{UB}} = \min_{1 \leq \nu \leq i}  \{ z^{(\nu)}_{\text{GTEP}} \}$.
Vectors $\mathbf{x}^\text{BEST}_{y}$, $y \in \mathcal{Y}$, contain the values of the investment variables in the solution corresponding to $z^{(i)}_{\text{UB}}$.

%When the algorithms stops at step D of iteration $i$, vectors $\mathbf{x}^{BEST}_{y}$, $y \in \mathcal{Y}$, contain the the values of the investment variables in the solution corresponding to the best computed upper bound.
%***********************************************************************************
At step D of iteration $i$, if the relative gap between upper and lower bounds is not greater than $\varepsilon$, the MILP problem $OPE_{y,w}$ is solved to determine the optimal operational variables associated to the investment decisions $\mathbf{x}_y^\text{BEST}$ and the solution procedure terminates.

%***********************************************************************************
The proposed algorithm is a Benders-type algorithm, where $GTEP_\text{LB}^{(i)}$ is the master problem and constraints \eqref{eq:cut} are the optimality cuts.
Feasibility cuts are not required: indeed, each problem $OPE_{y,w}^{(i)}$ is feasible since any set of values feasible for constraints
\eqref{eq:non_neg_ENP}-\eqref{eq:the_max_power}, \eqref{eq:pos_reserve}, \eqref{eq:CO2}-\eqref{eq:bat_final} and \eqref{non_neg_GCURT}-\eqref{eq:gas_LT_pos} also satisfies constraints \eqref{eq:balance}, \eqref{eq:reserve} and \eqref{eq:gas_balance}  for suitable non-negative values of variables $E^\text{NP}_{z,t,c,y,w}$, $R^\text{NP}_{z,t,c,y,w}$, $\text{OG}_{z,t,c,y,w}$, and $G_{n,z,t,c,y,w}^\text{CURT}$.

\section{Case study}
\label{sec:Case_study}
%This section discusses the analysis of challenges for the Italian energy system exploring the decarbonization pathways. The configuration of the Italian integrated system and the main input to the analysis are presented in Section \ref{sec:Scenario}, while the relevant results are discussed in Section \ref{sec:Results}. 
The model presented in Section \ref{sec:Model} can provide useful guidance on how an energy system should evolve to meet policy targets, while exploiting the potential for integration between the electricity and gas systems and taking into account the uncertainty of future fuel and CO$_2$ prices. 
To illustrate this, we show the analysis carried out on a case study related to the Italian energy system, based on publicly available data from 
ENTSO-E (the European Network of Transmission System Operators for Electricity), 
ENTSOG (The European Network of Transmission System Operators for Gas), 
TERNA (the Italian Transmission System Operators for Electricity), 
SNAM (the Italian Transmission System Operators for Gas)
and GME (the Italian Electricity Market Operator) \cite{TYND}, \cite{Terna}, \cite{cole2019cost}, \cite{ralon2017electricity}, \cite{EERA}, \cite{SnamTerna}. 
The planning period considered is 2023 to 2040. 
By means of the proposed model, a development plan for the Italian electricity and gas system is to be defined, indicating the interventions to be carried out year by year in order to meet the expected demand, taking into account policy objectives. 
The total electricity load of the Italian system is 334 TWh in 2023, with an assumed growth rate of 1\% per year in the following years.
In accordance with the Italian Implementation Plan \cite{ItalyMarketReform2020} of the European Commission, a penalty of 3,000 $\frac{\text{\texteuro}}{\text{MWh}}$ is imposed for electricity demand curtailment. 
The costs of unsupplied reserve and over-generation are set at 3,000 $\frac{\text{\texteuro}}{\text{MWh}}$ and 100 $\frac{\text{\texteuro}}{\text{MWh}}$, respectively. 
According to data published on the SNAM website \cite{SNAM}, the exogenous gas demand for the industrial and tertiary (residential and commercial) sectors is 501 TWh$_{th}$ in 2023. 
In subsequent years, an annual reduction of 0.6\% is assumed, based on the SNAM-TERNA scenarios of natural gas demand for final consumption \cite{SnamTerna}. 
A penalty of 3,000 $\frac{\text{\texteuro}}{\text{MWh$_{th}$}}$ has been set for natural gas demand curtailment.

For the period covered by the plan, the following policy objectives have been set for the Italian energy system:
\begin{enumerate}
    \item decommissioning of coal-fired and oil-fired power plants (with the exception of the oil-fired plant in central-south Italy) by 2024, according to \cite{TYND}; 
    \item $45\%$ lower limit for the penetration of renewable energies from 2030, according to \cite{SnamTerna}, i.e. in constraint \eqref{eq:penetration}, $\phi_{\textit{Italy},y} = 0.45$ for years $y = 2030, \ldots, 2040$;
    \item a maximum limit of 50 million tonnes of $CO_2$ emissions in 2030, decreasing by 1 million tonnes each year thereafter, according to the elaboration by CESI in \cite{Ministero}, \cite{ISPRA} of the long-term objectives for reducing the greenhouse gas impact of the Italian electricity system, i.e. in constraint \eqref{eq:CO2}, ${\overline{CO}_2}_{\textit{Italy},y} = 50-(y-2030)$ for years $y = 2030,\ldots,2040$.  
\end{enumerate} 
In order to achieve the 45\% target for the penetration of energy from renewable sources, it is necessary to identify an appropriate mix of non-programmable generation, programmable generation and flexibility resources that are capable of meeting demand under all conditions of primary energy sources (solar and wind). 
Due to the planned decommissioning of coal- and oil-fired power plants, potential investments in new hydro and thermal power plants must be considered in order to ensure an adequate level of programmable generation.
Given the lack of additional sites for hydroelectric valleys in Italy, pumped storage represents a viable option for the expansion of hydroelectric capacity. 
With regard to the expansion of thermoelectric capacity, the decision on the installation of new combined-cycle gas turbine (CCGT) plants and new gas turbine (GT) plants must take into account the future evolution of fuel costs and the reduction of CO$_2$ emissions.
We also want to evaluate candidate investments in new flexibility resources for energy shifting in time (e.g. batteries and PtG, in addition to the pumping stations mentioned above), new electricity transmission lines for energy shifting in space, and new gas interconnections to increase gas imports from regions with lower gas prices.
\color{black}
\begin{table}[ht!]
    \centering
    \begin{tabular}{l|rrr|rrr}
    \hline 
       \multirow{2}{*}{Fuel} & \multicolumn{3}{c|}{Scenario L} & \multicolumn{3}{c}{Scenario H} \\
                  &  2023 &  2030 &  2040 &  2023 &  2030 &   2040 \\
        \hline 
         Oil   & 90.50 & 91.25 & 71.58 & 90.50 & 85.81 & 102.14 \\
         IT gas   & 25.07 & 28.88 & 27.63 & 25.07 & 36.84 &  30.56 \\
         EU gas   & 23.07 & 26.88 & 25.63 & 23.07 & 34.84 &  28.56 \\
         GR gas   & 23.07 & 26.88 & 25.63 & 23.07 & 34.84 &  28.56 \\
         NAfr gas & 21.07 & 24.88 & 23.63 & 21.07 & 32.84 &  26.56 \\ 
         Coal     &  9.79 & 10.47 & 10.47 &  9.79 & 13.40 &  10.47 \\      
         \hline
    \end{tabular}
    \caption{Scenarios of fuel prices [$\frac{\text{\texteuro}}{\text{Gcal}}$] developed in \cite{TYND}.}
    \label{tab:Fuel_Price}
\end{table}

\begin{table}[ht!]
    \centering
    \begin{tabular}{c|ccc}
    \hline 
        Scenario & 2023 & 2030 & 2040 \\
        \hline 
               A &  33  &  35  &  75  \\
               B &  33  &  40  &  80  \\
               C &  33  &  53  & 100  \\
         \hline
    \end{tabular}
    \caption{Scenarios of CO$_2$ prices [$\frac{\text{\texteuro}}{\text{ton}}$]  developed in \cite{TYND}.}
    \label{tab:CO2_Price}
\end{table}

\begin{comment}
\begin{figure}[ht!]
    \centering
    % Panel 1
    \begin{subfigure}{\textwidth}
        \centering
    \includegraphics[width=0.9\linewidth]{Scenarios_Fuel.png}
    \caption{Scenarios of fuel prices [$\frac{\text{\texteuro}}{\text{Gcal}}$].}
    \label{fig:fuels}
    \end{subfigure}
%    \vspace{1 cm}
    %Panel 2
    \begin{subfigure} {\textwidth}
    \centering
    \includegraphics[width=0.9\linewidth]{Scenarios_CO2.png}
    \caption{Scenarios of CO$_2$ prices [$\frac{\text{\texteuro}}{\text{ton}}$].}
    \label{fig:CO2}
    \end{subfigure}
    \caption{Scenarios of fuel and CO$_2$ prices developed in \cite{TYND}.}
    \label{fig:scenarios}
\end{figure}
\end{comment}
In anticipative planning, scenarios provided by institutional bodies are typically considered to account for long-term uncertainties. In line with this approach, we represent the uncertainty of future fuel and CO$_2$ prices by six scenarios, obtained by combining scenarios developed by ENTSO-E and ENTSOG \cite{TYND}: the two scenarios (L and H) for future fuel prices reported in Table \ref{tab:Fuel_Price} and the three scenarios (A, B and C) for future CO$_2$ prices reported in Table \ref{tab:CO2_Price}.
%
%Specifically, Figure \ref{fig:fuels} shows the price trends for gas in scenarios L and H. 
%In addition, projections for oil fuel prices, based on estimates from \cite{TYND}, are also considered. 
%These projections, in $\frac{\text{\texteuro}}{\text{Gcal}}$, are as follows for the years 2023, 2030, and 2040: in scenario L, 90.50, 91.25, and 71.58; and in scenario H, 90.50, 85.81, and 102.14. 
Prices for the remaining years of the planning period are interpolated linearly.
%For coal, the price in 2023 is set to 9.79 $\frac{\text{\texteuro}}{\text{Gcal}}$, with estimates for 2024 of 9.87 $\frac{\text{\texteuro}}{\text{Gcal}}$ in scenario L and 10.30 $\frac{\text{\texteuro}}{\text{Gcal}}$ in scenario H. Notice that future coal price estimates are not relevant, as Italian energy policy mandates the decommissioning of all coal-fired power plants by 2024.
%
\medskip

In the following subsections we present (i) the configuration of the Italian electricity and gas systems in 2022, along with the candidate projects for the expansion of programmable and non-programmable generation, storage capacity, and transmission capacity; (ii) a brief outline of the solution procedure; (iii) the optimal plan determined by the stochastic programming model, with a comparison with the optimal plan obtained by solving the mean value problem.
%Decidere l'espansione degli impianti solari ed eolici, privilegiando le zone in cui sono più alti i coefficienti di producibilità. 

\subsection{The Italian electricity and gas systems in 2022 and candidate projects for generation and transmission expansion}
\label{sec:Scenario}
%%%%
%%%% 
%%%%
%\bigskip
To model the Italian electricity and gas systems, $11$ power zones and $5$ gas zones are defined:
\begin{itemize}
    \item seven power zones are the zones into which the Italian electricity market is divided, namely North (ITn), Central-North (ITcn), Central-South (ITcs), South (ITs), Calabria (ITcal), Sicily (ITsic) and Sardinia (ITsar);   
    \item four power zones represent the power systems to which the Italian system is either connected, namely Europe (EU), Montenegro (ME) and Greece (GR), or could potentially be connected, namely Tunisia (TN): these zones allow modelling imports and exports of the Italian system;
    \item two gas zones represent the Italian gas system, namely continental Italy+Sicily (ITA) and Sardinia (Sar): the ITA gas zone includes all Italian power zones except Sardinia, since the Italian gas transmission network does not have any bottlenecks that could limit the gas flow between the power zones included in the ITA gas zone;
    \item three gas zones, Europe (EU), Greece (GR) and North Africa (NAfr), represent the countries from which the Italian system can import gas.
\end{itemize}

The power generation capacity mix at the beginning of the planning period (year 2022) is shown in Table \ref{tab:Installed_cap}: column 1 shows the generation technologies, column 2 the total installed capacity by technology, and columns 3 to 9 the zonal details. In 2022 the installed capacity is 30\% non-programmable renewables, 24\% hydro and 46\% thermal. 
\begin{table}[ht!]
    \caption{Installed capacity [GW] per technology and zone in 2022.}
    \centering
    \normalsize
    \begin{tabular}{l|r|rrrrrrrr}
    \hline 
                 & TOT    & ITn    & ITcn  & ITcs  & ITs   & ITcal & ITsic & ITsar \\
        \hline 
%       Thermal & 49.527 & 22.573 & 1.619 & 8.510 & 6.894 & 3.506 & 4.599 & 1.826 \\
        Solar   & 21.08 &  9.40 & 2.51 & 2.94 & 3.32 & 0.55 & 1.52 & 0.84 \\
        Wind    & 11.21 &  0.19 & 0.21 & 2.04 & 4.61 & 1.19 & 1.90 & 1.07 \\
        Hydro   & 25.86 & 18.91 & 1.00 & 3.60 & 0.35 & 0.81 & 0.69 & 0.51 \\
        CCGT    & 37.08 & 19.79 & 1.59 & 5.16 &	3.54 & 3.28 & 3.08 & 0.63 \\
        GT      &  2.88 &  0.33 & 0.03 & 1.36 &	0.43 & 0.22 & 0.52 & 0.00 \\
        OIL     &  1.42 &  0.18 & 0.00 & 0.02 &	0.00 & 0.00 & 0.99 & 0.23 \\
        COAL    &  8.15 &  2.29 & 0.00 & 1.98 &	2.92 & 0.00 & 0.00 & 0.97 \\
        \hline
    \end{tabular}
    \label{tab:Installed_cap}
\end{table}

\begin{figure}[ht!]
    \centering
    \includegraphics[width=1\linewidth]{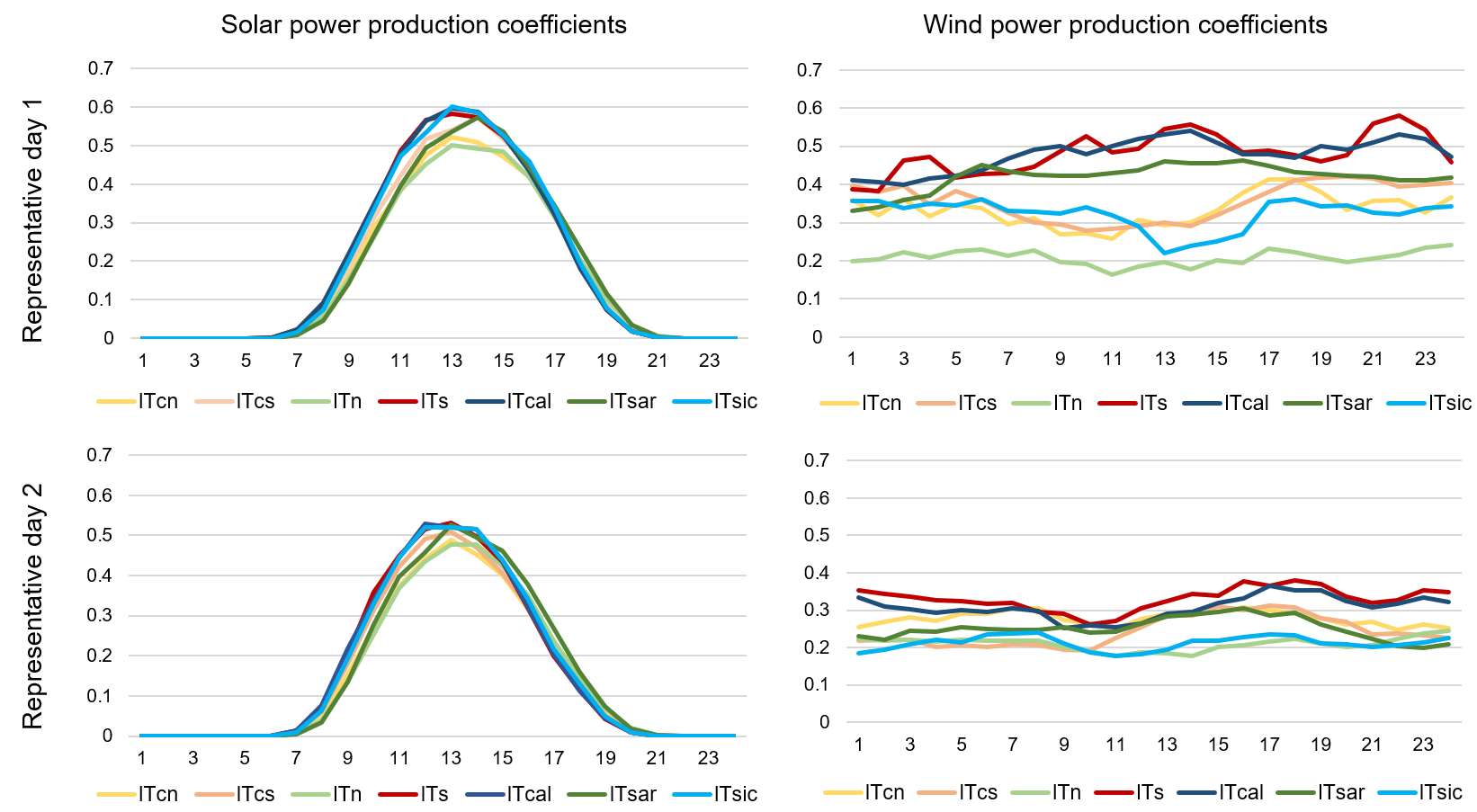}
    \caption{Hourly solar (left) and wind (right) power production coefficients in Italian power system zones in representative days 1 (top) and 2 (bottom).}
    \label{fig:rep_days}
\end{figure}

Investment decisions in solar and wind power capacity are affected by the hourly zonal power production coefficients $\mu_{z,t,c}$ and $\rho_{z,t,c}$.
Figure \ref{fig:rep_days} shows the coefficients $\mu_{z,t,c}$ and $\rho_{z,t,c}$ of the Italian power zones for two of the five representative days used to model short-term system operation. 
According to \cite{TYND}, \cite{cole2019cost} and \cite{SnamTerna}, investment costs for new solar and wind power capacity decrease over the planning period. 
In the first subperiod, for $y=2023, \ldots ,2030$, the investment costs are  
$IC^\text{S}_{z,y} = 880 - 40 (y-2023)$ for solar power capacity and  
$IC^\text{W}_{z,y} = 1180 - 40 (y-2023)$ for wind power capacity;
in the second subperiod, for $y=2030, \ldots, 2040$, the investment costs are  
$IC^\text{S}_{z,y} = 600 - 15 (y-2030)$ for solar power capacity and  
$IC^\text{W}_{z,y} = 900 - 21 (y-2030)$ for wind power capacity.
These costs apply to all zones $z \in \mathcal{Z}_\textit{Italy}$.

%%%%%%%%%%%%%%%%%%%%%%%%%%%%%%%%%%%%%%%%%%%%%
%%%%%%%%%%%%%%% THERMAL POWER PLANTS
%%%%%%%%%%%%%%%%%%%%%%%%%%%%%%%%%%%%%%%%%%%%%
%DA DIRE NELLA PARTE DI DEFINIZIONE DEL CASO STUDIO: non vi sono impianti candidati per le tecnologie OIL e COAL; di tutti gli impianti esistenti (both OIL and COAL) è imposta la dismissione, tranne che per l'impianto OIL in zona ITcs, che può essere o dismesso o mantenuto): 
The Italian thermal power generation system in 2022 consists of 28 power plant clusters: 
12 clusters of combined cycle power plants (CCGT1 to CCGT12), 6 clusters of gas turbine power plants (GT1 to GT6), 5 clusters of oil-fired power plants (OIL1 to OIL5) and 5 clusters of coal-fired power plants (COAL1 to COAL5).
The number of existing plants in 2022 and the number of plants to be decommissioned by 2024 are shown in columns 3 and 4 of Table \ref{tab:K}, respectively.
% ORDINAMENTO (decomm e max_new): = 0 e = 0; > 0 e = 0; > 0 e > 0; = 0 e > 0;
%         From & To & Code of pipeline $j \in \mathcal{J}_E \cup \mathcal{J}_C$ &    $\underline{F}_j^\text{J}$ [MW] &  $\overline{F}_j^\text{J}$ [MW] & $IC_j^\text{J}$ [$\text{M\texteuro}$] \\
%(Gcal/MWh)
\begin{table}[ht!]
\caption{Thermal power plants: data of clusters $k \in \mathcal{K}_z$, $z \in \mathcal{Z}_\textit{Italy}$.}
\centering
\begin{tabular}{llcccrrrc}
%    \begin{tabular} { p{1.2 cm} p{1.2 cm} R{1 cm} R{1 cm} R{1cm} R{1cm} R{0.8cm} R{2.0cm} R{2.0 cm} }
\hline 
Zone & Cluster & Number of & Number of & Number of & $\underline{P}_k$ & $\overline{P}_k$ & $C^{\text{SU}}_k$ & $HR_k$                             \\
     & code    & plants    & plants to & candidate & [MW]              & [MW]             & [k\texteuro]       & [$\frac{\text{Gcal}}{\text{MWh}}$] \\
     &         & in 2022   & be closed & plants    &                   &                  &                   &                                    \\ 
    \hline
    \multirow{8}{*}{ITn} 
    & COAL1  &  9 & 9 & 0 & 156.00 & 254.17 & 54& 2.937 \\
    & OIL1   &  2 & 2 & 0 &  24.45 &  78.35 & 20& 2.457 \\
    & OIL2   &  1 & 1 & 0 &  16.00 &  19.00 & 20& 2.457 \\
    & CCGT1  & 10 & 4 & 0 &  34.62 & 160.11 & 40& 2.199 \\
    & CCGT2  & 27 & 2 & 0 & 153.44 & 361.84 & 40& 1.670 \\
    & GT1    &  2 & 0 & 7 &  69.97 & 157.79 &  8& 2.114 \\
    & CCGT3  & 11 & 0 & 1 & 232.42 & 765.72 & 40& 1.680 \\    
    & CCGT13 &  0 & 0 & 7 & 139.14 & 359.29 & 40& 1.536 \\ 
    \hline
    \multirow{4}{*}{ITcn} 
    & CCGT4  &  2 & 0 & 0 & 190.50 & 385.55 & 40& 1.643 \\
    & CCGT5  &  4 & 1 & 0 &  66.68 & 204.00 & 40& 2.072 \\
    & GT2    &  2 & 1 & 0 &   5.00 &  16.00 &  8& 2.614 \\
    & GT7    &  0 & 0 & 2 &  25.20 &  63.00 &  8& 2.155 \\
    \hline 
    \multirow{7}{*}{ITcs}  
    & OIL3   &  1 & 0 & 0 &   0.00 &  15.60 & 20& 2.523 \\
    & COAL2  &  3 & 3 & 0 & 260.00 & 615.00 & 54& 2.211 \\
    & COAL3  &  2 & 2 & 0 &  27.30 &  65.00 & 54& 3.598 \\
    & CCGT6  &  5 & 3 & 0 &  63.70 & 184.66 & 40& 2.081 \\
    & CCGT7  & 10 & 1 & 3 & 178.77 & 423.85 & 40& 1.580 \\
    & GT3    & 11 & 8 & 2 &  36.46 & 123.38 &  8& 2.500 \\
    & CCGT14 &  0 & 0 & 3 & 286.67 & 566.67 & 40& 1.536 \\
    \hline    
    \multirow{5}{*}{ITs} 
    & CCGT8  &  7 & 0 & 0 & 195.43 & 505.93 & 40& 1.571 \\
    & COAL4  &  6 & 6 & 0 & 186.33 & 486.67 & 54& 3.024 \\
    & GT4    &  2 & 2 & 1 &  26.33 & 216.00 &  8& 2.443 \\
    & CCGT15 &  0 & 0 & 2 & 125.00 & 250.00 & 40& 1.536 \\
    & CCGT16 &  0 & 0 & 4 & 205.00 & 375.00 & 40& 1.536 \\
    \hline    
    \multirow{3}{*}{ITcal} 
    & CCGT9  & 6 & 0 & 0 & 229.67 & 547.07 & 40& 1.548 \\
    & GT5& 2 & 1 & 0     &  29.00 & 111.70 &  8& 2.614 \\
    & CCGT17 & 0 & 0 & 1 & 150.00 & 400.00 & 40& 1.536 \\
    \hline    
    \multirow{6}{*}{ITsic}      
    & GT6    & 5 & 0 & 0 &  71.84 & 104.73 &  8& 2.614 \\
    & CCGT10 & 3 & 0 & 0 & 179.33 & 493.00 & 40& 1.538 \\
    & OIL4   & 6 & 6 & 0 &  73.83 & 165.67 & 20& 2.457 \\
    & CCGT11 & 7 & 2 & 0 &  65.07 & 228.86 & 40& 1.786 \\
    & CCGT18 & 0 & 0 & 2 & 135.00 & 288.00 & 40& 1.536 \\
    & GT8    & 0 & 0 & 3 &  46.33 & 116.00 &  8& 2.102 \\
    \hline    
    \multirow{5}{*}{ITsar} 
    & COAL5  & 4 & 4 & 0 & 175.75 & 241.50 & 54& 3.283 \\
    & OIL5   & 3 & 3 & 0 &  38.33 &  76.67 & 20& 2.457 \\
    & CCGT12 & 1 & 1 & 0 & 115.00 & 630.00 & 40& 2.207 \\
    & CCGT19 & 0 & 0 & 1 & 115.00 & 630.00 & 40& 2.102 \\
    & GT9    & 0 & 0 & 2 &  80.00 & 200.00 &  8& 2.102 \\
    \hline
    \end{tabular}
    \label{tab:K}
\end{table}
In addition, 7 new CCGT clusters (CCGT13 to CCGT19) and 3 new GT clusters (GT7 to GT9) are candidates for investment.
For all CCGT and GT clusters, the number of candidate plants is given in column 5 of Table \ref{tab:K}.
The clusters differ in their technical characteristics: minimum power $\underline{P}_k$, capacity $\overline{P}_k$, start-up cost $C^{\text{SU}}_k$ and heat rate (i.e., fuel consumption per unit of generation) $HR_k$ are reported in columns 6 to 9.
%
%corresponding to 9.86 GW for CCGT and 2.44 GW for GT.
% NO: la produzione estera è rappresentata dagli impianti termoelettrici delle zone estere. C'è poi la domanda (inelastica) della zona estera e la curva di offerta della zona estera è formata dalle coppie quantità-prezzo degli impianti esteri e dall'importazione dall'Italia che il modello determina in modo da minimizzare il costo complessivo 
Similarly to hydro power, thermal power is a mature technology and investment costs $IC^\text{CCGT}=800$ $\frac{\text{\texteuro}}{\text{kW}}$ and $IC^\text{GT}=400$ $\frac{\text{\texteuro}}{\text{kW}}$, constant over the planning period, are considered, based on \cite{TYND} and \cite{SnamTerna}.
In each foreign zone $z \in \mathcal{Z}\setminus\mathcal{Z}_\textit{Italy}$, dummy thermal power plants are used to model the zonal stepwise aggregate supply curve: the capacity and the marginal cost of a dummy power plant represent the quantity/price pair that defines a step of the supply curve.

%%%%%%%%%%%%%%%%%%%%%%%%%%%%%%%%%%%%%%%%%%%%%
%%%%%%%%%%%%%%% HYDRO POWER PLANTS
%%%%%%%%%%%%%%%%%%%%%%%%%%%%%%%%%%%%%%%%%%%%%
%In the Italian hydro power system, there are three types of hydro power plants: (i) run-of-river, (ii) pumped-storage, and (iii) hydroelectric valleys. All hydro power plants of the same type in each zone $z \in \mathcal{Z}_\textit{Italy}$ are represented by a single equivalent plant. 
The Italian hydropower system in 2022 is represented by 18 equivalent plants of three types: 7 run-of-river plants, 7 hydroelectric valleys, and 4 pumped storage plants.
An equivalent plant represents all hydropower plants of the same type located in the same zone.
A new equivalent pumped storage hydropower plant is candidate in each Italian zone, except for zones ITn and ITcn.
%All hydropower plants considered (both existing and candidate) are programmable, with the exception of run-of-river plants.
%Table \ref{tab:H} reports data on both the equivalent plants existing at the beginning of the planning period and the candidate equivalent plants.
In Table \ref{tab:H}, for both existing and candidate equivalent plants $h \in \mathcal{H}_z, z \in \mathcal{Z}_\textit{Italy}$, column 1 reports the zone in which the plant is located, column 2 gives the plant code, 
and columns 3 to 5 give maximum power output $\overline{H}_h^\text{OUT}$, maximum pumping power $\overline{H}_h^\text{IN}$, and energy to power ration $EPR_h$, respectively.
Based on \cite{TYND} and \cite{SnamTerna}, the investment cost $IC_h^\text{H}=60$ $\frac{\text{\texteuro}}{\text{kW}}$, constant over the planning period, has been considered, since hydropower is a mature technology.

\begin{table}[ht!]
\caption{Data of equivalent hydro power plants $h \in \mathcal{H}_z$, $z \in \mathcal{Z}_\textit{Italy}$.}
    \centering
    \normalsize
    \begin{tabular}{llrrrr}
    \hline
%       $z \in \mathcal{Z}$ & $h \in \mathcal{H}_z$  & $\overline{H}_h^\text{OUT}$ (MW) & $\overline{H}_h^\text{IN}$ (MW) & $EPR_h$ (hours) & $IC_h^\text{H}$ (\texteuro/kW) \\
    \multirow{2}{*}{Zone}        
     & \multirow{2}{*}{Equivalent plant}
     & $\overline{H}_h^\text{OUT}$ 
     & $\overline{H}_h^\text{IN}$ 
     & $EPR_h$ \\
     &                  & [MW]                      & [MW]                   & [hours] \\
%       Zone & Equivalent plant  & $\overline{H}_h^\text{OUT}$ (MW) & $\overline{H}_h^\text{IN}$ (MW) & $EPR_h$ (hours) & $IC_h^\text{H}$ (\texteuro/kW) \\
         \hline
    \multirow{3}{*}{ITn}        
                           & Run-of-river 1                        & 7630 &    0 &  0    \\
                           & Hydroelectric valley 1                & 7755 &  406 & 500    \\
                           & Pumped storage hydro 1                & 3525 & 3154 &  4.4 \\
    \hline
    \multirow{3}{*}{ITcn}
                           & Run-of-river 2                        &  232 &    0 &  0    \\
                           & Hydroelectric valley 2                &  464 &    0 & 520    \\
                           & Pumped storage hydro 2                &  300 &  271 &  6.4 \\
    \hline
    \multirow{4}{*}{ITcs}
                           & Run-of-river 3                        & 1029 &    0 &  0    \\
                           & Hydroelectric valley 3                & 1562 &  201 & 540    \\
                           & Pumped storage hydro 3                & 1006 &  986 &  3.4 \\
                           & New pumped storage hydro 1 (NPSH1)    & 1000 & 1000 & 14    \\
    \hline
    \multirow{3}{*}{ITs} 
                           & Run-of-river 4                        &  117 &    0 &  0    \\
                           & Hydroelectric valley 4                &  228 &    0 & 950    \\
                           & New pumped storage hydro 2 (NPSH2)    &  450 &  450 & 14    \\
    \hline
    \multirow{3}{*}{ITcal}  
                           & Run-of-river 5                        &  273 &    0 &  0    \\
                           & Hydroelectric valley 5                &  534 &    0 & 950    \\
                           & New pumped storage hydro 3 (NPSH3)    & 1250 & 1250 & 14    \\
    \hline
    \multirow{4}{*}{ITsic}
                           & Run-of-river 6                        &   32 &    0 &  0    \\
                           & Hydroelectric valley 6                &  208 &   52 & 490    \\
                           & Pumped storage hydro 4                &  450 &  445 &  3.6 \\
                           & New pumped storage hydro 4 (NPSH4)    &  480 &  480 & 14    \\
    \hline
    \multirow{3}{*}{ITsar}  
                           & Run-of-river 7                        &   87 &    0 &  0    \\
                           & Hydroelectric valley 7                &  423 &  218 &  5.8 \\
                           & New pumped storage hydro 5 (NPSH5)    &  800 &  800 & 14    \\
    \hline
    \end{tabular}
    \label{tab:H}
\end{table}

%%%%%%%%%%%%%%%%%%%%%%%%%%%%%%%%%%%%%%%%%%%%%
%%%%%%%%%%%%%%% STORAGE
%%%%%%%%%%%%%%%%%%%%%%%%%%%%%%%%%%%%%%%%%%%%%
At the beginning of the planning period, storage capacity is provided only by the hydropower system.
%, as detailed in Table \ref{tab:tot_cap_store_initial}.
In addition to building new hydroelectric pumping stations, storage capacity can be increased by installing batteries.
Candidate battery technologies are 
lithium-ion batteries (Li-Batt), with energy to power ratio $EPR_\text{Li-Batt}=4$ hours and operating cost $C^\text{B}_\text{Li-Batt}= 20$ $\frac{\text{\texteuro}}{\text{MWh}}$, and
sodium-ion batteries (Na-Batt), with energy to power ratio $EPR_\text{Na-Batt}=6$ hours and operating cost $C^\text{B}_\text{Na-Batt}= 30$ $\frac{\text{\texteuro}}{\text{MWh}}$. 
For each battery technology, the maximum installable capacity is $\overline{B}^\text{CAP}=600$ MW in each zone.
\begin{comment}
\begin{table}[ht!]
    \centering
    \begin{tabular}{l|r|rrrrrrrr}
    \hline 
         Technology   & TOT    & ITn    & ITcn   & ITcs   & ITs    & ITcal  & ITsic  & ITsar \\
        \hline 
Hydroelectric valleys & 579.05 &	387.75 & 24.13 & 84.35 & 21.66 & 50.73 & 10.19 & 0.25 \\
Hydro pumped storage  &   2.25 &      1.55 & 0.19 &  0.34 & 0.00 & 0.00 & 0.162 & 0.00 \\	
%Lithium-ion batteries &  0.000 &  0.000 &  0.000 &  0.000 &  0.000 &  0.000 &  0.000 & 0.000 \\ 
%Sodium-ion batteries  &  0.000 &  0.000 &  0.000 &  0.000 &  0.000 &  0.000 &  0.000 & 0.000 \\ 
    \hline
    \end{tabular}
    \caption{Storage capacity [GWh] per technology and zone in 2022.}
    \label{tab:tot_cap_store_initial}
\end{table}
\end{comment}
According to \cite{TYND}, \cite{cole2019cost} and \cite{SnamTerna}, investment costs for batteries decrease over the planning period: 
in the first subperiod, for $y=2023, \ldots, 2030$, the investment costs are
$IC^\text{B}_{\text{Li-Batt},y} = 288.5 - 20.5 (y-2023)$ for lithium-ion batteries
and  
$IC^\text{B}_{\text{Na-Batt},y} = 334 - 22 (y-2023)$ for sodium-ion batteries;
in the second subperiod, for $y=2030, \ldots, 2040$, the investment costs are  
$IC^\text{B}_{\text{Li-Batt},y} = 145 - 2.4 (y-2030)$ for lithium-ion batteries
and  
$IC^\text{B}_{\text{Na-Batt},y} = 180 - 2.4  (y-2030)$ for sodium-ion batteries. 
These costs apply to all zones $z \in \mathcal{Z}_\textit{Italy}$.
Additional flexibility can be provided to the system by the PtG technology, which produces synthetic gas from excess electricity production.
In 2022 there is no PtG capacity. 
A PtG plant with capacity $\overline{PtG}^{\text{CAP}} = 10,000$ MW$_{th}$, efficiency $\eta=68\%$ and investment cost $IC^\text{PtG}=300$ 
$\frac{\text{\texteuro}}{\text{kW$_{th}$}}$
%\texteuro/kWh$_{th}$ 
is candidate in the ITcn zone.

\begin{table} [ht!]
\caption{Data of power transmission lines $l \in \mathcal{L_E} \cup \mathcal{L}_C$ and gas pipelines $j \in \mathcal{J_E} \cup \mathcal{J}_C$.}
    \centering
    \normalsize
    \begin{tabular}{llcrrr}
    \hline
    \multicolumn{3}{l}{Power transmission lines $l \in \mathcal{L_E} \cup \mathcal{L}_C$} \\
%    \hline 
         From & To & Code of line                             & $\underline{F}_l^\text{L}$ [MW] &  $\overline{F}_l^\text{L}$ [MW] & $ IC_l^\text{L} [$\text{M\texteuro}$]$      \\
    \hline
         ITn   & ITcn  & L-E1  & $-2500$ &  4000 &    - \\
         ITcn  & ITcs  & L-E2  & $-2700$ &  2100 &    - \\
         ITcs  & ITs   & L-E3  & $-4600$ & 10000 &    - \\
         ITs   & Tcal  & L-E4  & $-3350$ & 10000 &    - \\
         ITsar & ITcn  & L-E5  & $- 400$ &   400 &    - \\
         ITsar & ITcs  & L-E6  & $- 720$ &   900 &    - \\
         ITsic & ITcal & L-E7  & $-1500$ &  1500 &    - \\
         UE    & ITn   & L-E8  & $-5365$ & 10425 &    - \\
         ME    & ITcs  & L-E9  & $- 600$ &   600 &    - \\
         GR    & ITs   & L-E10 & $- 500$ &   500 &    - \\
\hdashline
         ITn   & ITcn  & L-R1  & $- 400$ &   400 &  240 \\
         ITn   & ITcn  & L-R2  & $-1000$ &   600 &  600 \\
         ITcn  & ITcs  & L-R3  & $- 150$ &   150 &   90 \\
         ITcn  & ITcs  & L-R4  & $-1000$ &  1000 & 1000 \\
         ITcs  & ITs   & L-R5  & $- 900$ &     0 &  270 \\
         ITcs  & ITs   & L-R6  & $- 200$ &     0 &  600 \\
         EU    & ITn   & L-R7  & $- 660$ &   660 &  198 \\
         EU    & ITn   & L-R8  & $- 100$ &   100 &   40 \\
         EU    & ITn   & L-R9  & $- 630$ &   630 &  400 \\
         EU    & ITn   & L-R10 & $- 600$ &   600 &  252 \\
         EU    & ITn   & L-R11 & $- 600$ &  1000 & 1000 \\
         ME    & ITcs  & L-R12 & $- 600$ &   600 &  600 \\
         ITcs  & ITsic & L-N1  & $-1000$ &  1000 & 1300 \\
         ITsar & ITsic & L-N2  & $-1000$ &  1000 & 1300 \\
         TN    & ITsic & L-N3  & $- 600$ &   600 &  360 \\
 \hline
 \hline
 \multicolumn{6}{l}{Gas pipelines $j \in \mathcal{J_E} \cup \mathcal{J}_C$} \\
% \hline
           From & To & Code of pipeline                         & $\underline{F}_j^\text{J}$ [MW] &  $\overline{F}_j^\text{J}$ [MW] & $IC_j^\text{J}$ [$\text{M\texteuro}$]     \\
\hline
         EU   & ITA & J-E1 &  $-38661$ & 38661 &     - \\
         NAfr & ITA & J-E2 &  $-36620$ & 36620 &     - \\
         GR   & ITA & J-E3 &  $-10873$ & 10873 &     - \\
\hdashline
         GR   & ITA & J-R1 &  $-10873$ & 10873 &  1011 \\
         NAfr & Sar & J-N1 &  $ -9750$ &  9750 &  1502 \\
         Sar  & ITA & J-N2 &  $ -9750$ &  9750 &  1004 \\
 \hline
    \end{tabular}
    \label{tab:L}
\end{table}
With regard to the spatial shift of energy, in 2022 the power zones are connected by 10 power transmission lines: ITn-ITcn, ITcn-ITcs, ITcs-ITs, ITs-ITcal, ITsar-ITcn, ITsar-ITcs, ITsic-ITcal, EU-ITn, ME-ITcs, and GR-ITs.
Candidate projects for the development of the power transmission system are
\begin{itemize}
    \item 12 projects to reinforce existing lines: 2 for the ITn-ITcn line, 2 for the ITcn-ITcs line, 2 for the ITcs-ITs line, 5 for the EU-ITn line and 1 for the ME-ITcs line;
    \item 3 new connections: ITcs-ITsic, ITsar-ITsic and TN-ITsic.
\end{itemize}
Data on both existing and candidate power lines are reported in Table \ref{tab:L}.
A line $l \in \mathcal{L}_E \cup \mathcal{L}_C$ is represented as an arc from the zone given in column 1 (tail of the arc) to the zone given in column 2 (head of the arc): flows in the conventional direction (tail to head) are positive and flows in the opposite direction are negative.
Column 3 gives the line code (E-Existing line; R-Reinforcement project; N-New line). 
Columns 4 and 5 show the maximum flow in the direction opposite to the conventional one and in the conventional direction, respectively. 
The investment costs for the candidate lines $l \in \mathcal{L}_C$ are given in column 6.

%

\begin{comment}
\begin{table}[ht!]
    \centering
    \begin{tabular}{l|ccccc}
    \hline 
        Gas zone $n$        &  ITA &  UE  &  GR  & NAfr \\
        \hline 
         $\overline{G}_{n}$ & 22.2 & 38.7 & 21.7 & 42.4 \\
        \hline
    \end{tabular}
    \caption{Maximum hourly natural gas supply [GW$_{th}$] in each gas zone.}
    \label{tab:Gas_supply}
\end{table}
\end{comment}
%%%%%%%%%%%%%%% Description of gas production
Regarding gas supply, 
%Table \ref{tab:Gas_supply} shows 
according to the SNAM-Terna scenario elaborations \cite{SnamTerna},
the maximum hourly gas supply $\overline{G}_{n}$ for foreign gas zones UE, GR and NAfr are set to 38.7 GW$_{th}$, 21.7 GW$_{th}$ and 42.4 GW$_{th}$, respectively. 
The maximum hourly gas supply of the ITA zone, $\overline{G}_{ITA}$, is set to 22.2 GW$_{th}$ and includes both domestic production and LNG imports.
In the ITA gas zone there are also storage facilities with total capacity $\overline{G}_{\text{ITA}}^{\text{LT}} =$ 190 TWh$_{th}$ and hourly injection and withdrawal rates $\overline{G}_{\text{ITA}}^{\text{IN}}=\overline{G}_{\text{ITA}}^{\text{OUT}} =$ 148 GW$_{th}$.
In 2022, the ITA gas zone is connected to the three import zones by the EU-ITA, GR-ITA and NAfr-ITA pipelines. 
The flow through the EU-ITA pipeline represents the supply of gas from Russia and the North Sea through the entry points Passo Gries, Tarvisio, Gorizia. 
The flow through the GR-ITA pipeline (also known as the Trans Adriatic Pipeline - TAP) represents the supply of gas from Azerbaijan through Turkey, Greece and Albania to Melendugno, Puglia in southern Italy. 
The flow through the NAfr-ITA pipeline represents the supply of gas from Algeria and Tunisia through the entry points of Mazara Del Vallo and Gela in Sicily, respectively.
The Sardinia (Sar) gas zone is not connected to other gas zones.
%\textcolor{red}{Zona non connessa tramite pipeline con altre zone, ma riceve una fornitura tramite metaniere} %(dire la provenienza del gas che usa: import con metaniere e/o produzione interna).
Candidate projects for the development of the gas transmission system are the doubling of the GR-ITA pipeline and the two branches of the \textit{Gasdotto Algeria Sardegna Italia} (GALSI), namely the pipeline from Algeria to Sardinia (NAfr-Sar) and the pipeline from Sardinia to the Italian mainland (Sar-ITA).  
Gas pipelines are represented in the model in a similar way to power transmission lines: data on existing and candidate gas pipelines are given in Table \ref{tab:L}.

Note that both hydro and gas storage levels are checked at the end of every week, i.e. $M = 7$ in constraints \eqref{eq:hydro_LT_0}$-$\eqref{eq:hydro_LT_pos} and \eqref{eq:gas_LT_0}$-$\eqref{eq:gas_LT_pos}. 

%As described in detail in the previous section, the natural gas demand is composed of two parts: (i) the demand related to the thermal power generation, which is endogenously determined, and (ii) the industrial and tertiary (residential and offices) demand, which is an input to the model. 

%\bigskip

\subsection{Solving the GTEP stochastic programming model}
\label{sec:algo}
The model was coded in GAMS 24.7.4 and solved on a computer with two 2.10 GHz Intel Xeon Platinum 8160 CPU processors and 128 GB of RAM, using Gurobi with GUSS \cite{bussieck2012guss}, a GAMS tool specifically designed to efficiently handle the batch processing of optimization problems with identical structure but varying data inputs. 
Due to its large dimension (3'226'806 real variables, 1'869'660 integer variables and 6'961'460 constraints), solving the stochastic model as a monolithic program is computationally infeasible, requiring the application of the solution algorithm introduced in Section \ref{sec:Algorithm}.
%Specifically, we set the tolerance for convergence test at $\varepsilon=10^{-3}$, while solving at each iteration the master problem up to optimality.
%Regarding the operational subproblems, 
At iteration $i$ of the proposed solution algorithm, we solve by means of language extension GUSS 108 independent second-stage problems $OPE_{y,w}^{(i)}$, each of which determines the operation of the system in year $y$ and scenario $w$.
More precisely, GUSS solves the first second-stage problem $OPE_{1,1}^{(i)}$ (referred to as the 'base case') and stores the factors of its optimal basis. 
This information is used to obtain an advanced initial basis for solving the following second-stage problems (called the 'updated subproblems').
%At the end of each outer iteration, it is checked whether $(UB-LB)/UB$ is less than a given tolerance, where $LB$ is the optimal value of the objective function of the master problem, and $UB$ is the sum of the optimal investment cost of the master problem and the optimal expected cost of the second-stage problems. 
%If $(UB-LB)/UB$ is less than the tolerance, the iteration stops, otherwise the master problem is updated by adding new Benders cuts generated from the dual information provided by the second-stage problems, and a new outer iteration is performed.
With the termination tolerance set to $\varepsilon=10^{-3}$, the algorithm converged in 24 iterations, always solving the master problem $GTEP_{LB}^{(i)}$ to optimality.
In the most computationally expensive iteration, iteration 24, the master problem (comprising $4,795$ variables and $5,441$ constraints) was solved in $10.02$ seconds, while the time taken to solve the second-stage problems (each one comprising $306,918$ variables and $310,651$ constraints) was $120.29$ seconds for the base case and $3.34$ seconds for each subsequent updated subproblem.
The total time taken to solve the case study problem was 10,863 seconds (i.e., 3 hours, 1 minute and 3 seconds).
After 24 iterations, the algorithm converged to the optimal solution of problem $\underline{GTEP}$, with an associated cost of 540,722 M\texteuro.
Based on this solution, we then determined a feasible solution to the original GTEP problem, which has a cost of 540,832 M\texteuro, only 0.02\% higher than the optimal value of problem $\underline{GTEP}$. This highlights the effectiveness of the proposed procedure in producing near-optimal solutions.
%Table \ref{tab:CPU_time} shows the size of the problems to be solved and the CPU time required by solution of the master problem, the solution of the base case, and the solution of each update subproblem. As can be noted, starting from the advanced basis of the base case and its factorisation, the computation time of subproblems is drastically reduced: while the solution of the base case takes 120.291 seconds, each of the updated subproblems is solved in only 3.342 seconds on average. Since the number of subproblems is $|\mathcal{Y}| \cdot |\mathcal{W}| = 108$, the total time required to solve all subproblems at the last iteration is 477.885 seconds.
\color{black}

\begin{comment}
 \begin{table}[ht!]
    \centering
    \begin{tabular}{lrrr}
    \hline 
                                      & \# Constraints & \# Decision Variables & CPU time (seconds) \\
        \hline 
         Master Problem               &          5,441 &                 4,795 &             10.018 \\
         First Subproblem (Base Case) &        318,651 &               306,918 &            120.291 \\
         Updated Subproblem           &        318,651 &               306,918 &              3.342 \\   
         \hline
    \end{tabular}
    \caption{Size and solution time (seconds) of master problem and subproblems at the last iteration of Benders algorithm.}
    \label{tab:CPU_time}
\end{table} 
\end{comment}

\subsection{The optimal capacity expansion plan}
%Analysis of the optimal solution
\label{sec:Results}
%\bigskip
In this section we present the optimal capacity expansion plan determined by the stochastic programming model (SPM plan) to achieve the required policy objectives (i.e., decommissioning of coal- and oil-fired power plants, minimum penetration of renewables and CO$_2$ mitigation) at minimum investment, operating and penalty costs. 
We then highlight the differences between the SPM plan and the optimal plan determined by the corresponding mean value problem (MVP plan), where each uncertain parameter is assigned the expected value of its realisations in the scenarios representing uncertainty in the stochastic programming model. 
%The numerical results show that in some scenarios with high prices of gas and CO$_2$ emissions, it is cheaper not to meet the demand, even with the penalty cost of not supplying electricity, than to buy the fuels and pay the CO$_2$ emission costs for thermoelectric production.
%In particular, we will observe that the expansion plan determined by MVP results in an energy system that is not adequate to meet the electricity demand in some scenarios.

The SPM plan is shown in Table  \ref{tab:yearly_cap_stoc}.
\begin{table}[ht!]
\caption{Expansion plan determined by the stochastic programming model.}
    \centering
%    \begin{tabular}{l|rrrrrrrrr}
\begin{tabular} { |p{0.6 cm}| R{0.72 cm}| R{0.74 cm} | R{0.89 cm} | R{0.72 cm} |R{0.82cm} | R{1.02 cm}  | R{1.08 cm} | R{1.01 cm} | R{0.91 cm} |}
\hline 
%\multirow{2}{*}{Year} & Solar & Wind  & CCGT &  GT  & Hydro & Li-Batt & Na-Batt & PtG         & \multirow{2}{*}{Connections} \\
%                      & [GW]  & [GW]  & [GW] & [GW] & [GW]  & [GW]    & [GW]    & [GW$_{th}$] &                              \\
                 Year & Solar &  Wind & CCGT &  GT  & Hydro & Li-Batt & Na-Batt &     PtG     & Links \\
                      &  [GW] &  [GW] & [GW] & [GW] &  [GW] &  [GW]   &  [GW]   & [GW$_{th}$] &             \\
\hline 
                 2023 & 20.62 & 14.36 &      &      &       &         &         &             &             \\
\hdashline
                 2024 &  8.25 &  5.46 & 0.63 & 0.40 &       &         &         &             & J-N1 J-N2   \\ 
\hdashline
                 2025 &  7.84 &  5.74 & 3.14 & 0.56 &  2.98 &         &         &             & L-R7 L-R9  \\  
\hdashline
                 2026 &  4.54 &  3.16 &      &      &       &         &         &             &             \\ 
\hdashline
                 2027 &  0.99 &  2.57 &      &      &       &         &         &             &             \\ 
\hdashline
                 2028 &  0.81 &  2.10 &      &      &       &         &          &            &             \\ 
\hdashline
                 2029 &  1.98 &  5.84 &      &      &       &         &          &            &             \\ 
\hdashline
                 2030 &  4.50 & 11.67 & 0.29 & 0.12 &       &    0.29 &          &       0.29 & L-N3        \\ 
\hdashline
                 2031 &  0.72 &  1.17 &      &      &       &         &          &            &             \\  
\hdashline
                 2032 &  0.59 &       &      &      &       &         &          &            &             \\ 
\hdashline
                 2033 &       &  1.75 &      &      &       &         &     0.35 &            &             \\ 
\hdashline
                 2034 &  0.70 &  2.09 &      &      &       &         &          &       0.14 &             \\ 
\hdashline
                 2035 &  1.70 &  5.06 &      &      &       &    2.04 &          &            &             \\ 
\hdashline
                 2036 &  0.41 &  1.92 &      &      &       &         &          &            &             \\ 
\hdashline
                 2037 &  0.64 &  1.22 &      &      &       &         &          &            &             \\ 
\hdashline
                 2038 &  0.06 &       &      &      &       &         &          &            &             \\ 
\hdashline
                 2039 &  0.59 &  2.27 &      &      &       &         &          &            &             \\ 
\hdashline
                 2040 &  1.17 &  3.14 &      &      &       &         &          &            &             \\
\hline
                Total & 56.11 & 69.52 & 4.06 & 1.08 &  2.98 &    2.33 &     0.35 &       0.43 &             \\
\hline
\end{tabular}
\label{tab:yearly_cap_stoc}
\end{table}
In 2023, the plan foresees significant investment in solar and wind power plants, accounting for 36.8\% and 20.7\% respectively of total investment in these two technologies over the planning period.
In 2024, the plan foresees four different types of intervention:
(i) the construction of the two branches of the GALSI project, i.e. J-N1, between North Africa and Sardinia, and J-N2, between Sardinia and the Italian mainland, which will allow more North African gas, the cheapest, to be imported;
(ii) the construction of all the candidate gas-fired power plants in ITsar zone, namely the CCGT19 plant and the two GT9 plants, replacing the existing coal and oil-fired thermal power plants;
(iii) the closure of all power plants to be decommissioned within the first two years of the planning period; 
(iv) further investment in solar and wind power plants, representing 14.7\% and 7.8\%, respectively, of the total investment in these technologies. 
The interventions foreseen in 2025 are: (i) further investments in solar and wind power plants (14\% and 8.3\% of total investment respectively); (ii) construction of new thermal power capacity in the continental zones, namely one plant for each of the following clusters: GT1, CCGT3 and CCGT13 in ITn zone; GT7 in ITcn zone; GT3, CCGT7 and CCGT14 in ITcs zone; GT4, CCGT15 and CCGT16 in ITs zone; CCGT17 in ITcal zone; (iii) the realisation of all the candidate projects for new pumped storage capacity, except for the NPSH1 project in zone ITcs; (iv) the reinforcement of the electricity connection between Northern Italy and Northern Europe through the construction of the L-R7 and L-R9 projects.
From 2026 to 2029, only additional investment in new solar and wind power capacity is planned, amounting to 14.8\% and 19.7\% of total investment in the respective technologies over the planning period.
%
% al 2030 sono imposti particolari vincoli di decarbinizzazione che spiegano gli interventi previsti dalla soluzione ottima? 
In 2030, the plan foresees 
(i) investments in solar and wind energy, amounting to 8.8\% and 16.8\% of the total investments in the respective technologies; 
(ii) the installation of new thermal power capacity in the ITsic zone (CCGT18 and GT8); 
(iii) 0.29 GW new capacity of lithium-ion batteries; 
(iv) 0.29 GW of PtG capacity; 
(v) the construction of the new power line between Sicily and Tunisia (L-N3).  
From 2031 to 2040, investments in solar and wind power account for 11.7\% and 26.7\% of total investment in the respective technologies. Investment decisions in flexibility sources include 0.35 GW new capacity of sodium-ion batteries in 2033, 0.14 GW of PtG capacity in 2034, and 2.04 GW new capacity of lithium-ion batteries in 2035.
\begin{table}[ht!]
\caption{Power generation capacity [GW] by technology and zone resulting from the SPM plan at the end of 2040.}
    \centering
    \normalsize
    \begin{tabular}{l|r|rrrrrrrr}
    \hline 
         Technology   & TOT    & ITn    & ITcn   & ITcs   & ITs    & ITcal  & ITsic  & ITsar \\
        \hline 
Solar                 &	77.17 & 40.00	& 10.00 & 15.00 &  5.22 &  2.00	&  2.95 & 2.00 \\
Wind                  &	80.71 &  8.00	&  8.00 & 13.95 & 24.91 & 10.00	& 12.34 & 3.50 \\
Hydro                 &	28.84 & 18.91	&  1.00 &  3.60 &  0.80 &  2.06	&  1.17 & 1.31 \\
CCGT                  &	37.50 & 19.56	&  1.38 &  5.17 &  4.17 &  3.68	&  2.91 & 0.63 \\
GT	                  &  2.21 &  0.47	&  0.08 &  0.49 &  0.22 &  0.11	&  0.64 & 0.20 \\
OIL                   &	 0.02 &  0.00	&  0.00 &  0.02 &  0.00 &  0.00	&  0.00 & 0.00 \\
%COAL                  &  0.000 &  0.000	&  0.000 &  0.000 &  0.000 &  0.000	&  0.000 & 0.000 \\
%Lithium-ion batteries &  2.326 &  0.600 &  0.056 &  0.335 &  0.600 &  0.000 &  0.299 & 0.437 \\ 
%Sodium-ion batteries  &  0.347 &  0.000 &  0.000 &  0.000 &  0.347 &  0.000 &  0.000 & 0.000 \\ 
    \hline
    \end{tabular}
        \label{tab:tot_cap_stoc}
\end{table}
The power generation capacity by technology and zone resulting from SPM plan at the end of 2040 is shown in Table \ref{tab:tot_cap_stoc}: the generation mix consists of 70\% non-programmable renewables, 12.5\% hydro and 17.5\% thermal, thus presenting a substantially different structure from the initial one in 2022, shown in Table \ref{tab:Installed_cap}.
Wind power is the main generation technology, with 69.50 GW of new capacity installed, mainly in zones ITcs, ITs, ITcal and ITsic, where wind power production coefficients are high, as shown in Figure \ref{fig:rep_days}. In the northern regions, investment is mainly in solar technology, which increases by 56.09 GW.

\begin{table}
\caption{Investment costs [M\texteuro] of the SPM plan and the MVP plan.}
    \centering
    \normalsize
    \begin{tabular}{l|rr}
    \hline
                             & SPM plan & MVP plan \\
    \hline
%     Decommissioning          &      201 &      201 \\
    Solar                    &   42,229 &   42,180 \\
    Wind                     &   62,701 &   58,103 \\
    Thermal                  &    3,488 &    3,488 \\
    Hydro                    &    2,335 &    1,959 \\
    Batteries                &    1,181 &      597 \\
    PtG                      &       79 &        0 \\
    Power Transmission Lines &      779 &      858 \\  
    Gas Pipelines            &    2,222 &    2,222 \\
    \hline
    Total investment costs   &  115,013 &  109,407 \\
%     Total investment costs   &  115,214 &  109,608 \\
        \hline
    \end{tabular}
    \label{tab:tot_costs}
\end{table}

Compared to the SPM plan, the MVP plan invests less in renewable energy and in flexibility resources for the shift of energy in time: indeed, it suggests to install 
(i) 3.43 GW of new solar capacity and 4.19 GW of new wind capacity in 2030, instead of 4.50 GW and 11.67 GW, respectively, as in the SPM plan; 
(ii) less new pumped storage hydro capacity, as the candidate plant NPSH4 in the ITsic zone is not to be built; 
(iii) less electric storage capacity, since in 2035 investment in lithium-ion batteries is only 1.03 GW and no sodium-ion batteries are to be installed; 
(iv) no PtG capacity.  
The MVP plan, on the other hand, provides for more flexibility in terms of the spatial transfer of energy: the construction of the L-R7 and L-R9 lines (reinforcing the connection between ITcn and ITcs) is brought forward to 2024, and the L-R3 line (the construction of which is not included in the SPM plan) is to be built in 2035.
In 2040, the energy mix resulting from the MVP plan consists of 68.7\% non-programmable renewable generation, 13\% hydro and 18.3\% thermoelectric.
The investment costs of the two plans are detailed in Table \ref{tab:tot_costs}.

%\subsection{The optimal operation of the integrated electricity and gas systems}
%\label{sec:Results-operation}
For each plan and scenario, Table \ref{tab:oper_costs} shows the total value over the planning period of each term $m$ of the second-stage costs, 
\begin{comment}
i.e. the values
\begin{equation}
\label{eq:obj}
\sum_{y \in \mathcal{Y}} \sum_{c \in \mathcal{C}^y} \psi_c  \sum_{t=1}^{24} \tau^m_{t,c,y}, \quad \quad 1 \leq m \leq 10,  \nonumber 
\end{equation}
\end{comment}
We note that in both plans, in all scenarios, the total second-stage costs are relevant compared to the investment cost: 
in the SPM plan, the total second-stage costs vary between 3.5 times (in the L-A scenario) and 4 times (in the H-C scenario) the investment cost; 
in the MVP plan, the total second-stage costs vary between 3.75 times (in the L-A scenario) and 7 times (in the H-C scenario) the investment cost. 
This confirms the need for detailed modeling of short-term operation so that an accurate estimate of these costs can be taken into account when determining expansion plans. 
We also note that by far the highest operating cost is for thermoelectric production.
In the SPM plan, the cost of thermoelectric production in the H-C scenario is 1.14 times the cost in the L-A scenario. 
Similarly, in the MVP plan, the cost of thermoelectric production in the H-C scenario is about 15\% higher than in the L-A scenario: 
however, in the H-C scenario, minimising the second-stage costs leads to a curtailment of electricity demand. 
In fact, the energy system resulting from the MVP plan has less renewable capacity, so the residual load to be met by programmable power plants is greater; in the MVP plan there is also less electricity storage capacity and PtG technology is not installed. 
Therefore, in the H-B, L-C and H-C scenarios, the high fuel and CO$_2$ emission costs make it cheaper to pay the penalty for the energy not supplied than to produce with gas-fired thermal power plants: demand curtailment also avoids the expensive start-up of thermal power plants and their production for several consecutive hours (forced by the minimum up-time and minimum power constraints).  
On the other hand, the system resulting from the SPM plan meets electricity demand at all hours of the planning period in all scenarios: the higher share of renewable capacity reduces the residual load to be covered by programmable generation, and the higher storage capacity allows more effective management of imbalances between generation and demand.

\begin{table}[ht!]
\caption{Operating costs and penalties [M\texteuro] over the planning period in the six scenarios of fuel and CO$_2$ prices for the SPM plan (above) and the MVP plan (below).}
\centering
\begin{tabular}{c|l|rr|rr|rr}
\hline
    & Plan                       &     SPM &     SPM &     SPM &     SPM &     SPM &     SPM \\ 
$m$ & Scenarios                  &     L-A &     H-A &     L-B &     H-B &     L-C &     H-C  \\ 
\hline
  1 & Hydro power production     &  31,034 &  31,067 &  31,045 &  31,063 &  31,062 &  31,081 \\
  2 & Start-up of thermal plants &   1,625 &   1,677 &   1,607 &   1,658 &   1,665 &   1,722 \\
  3 & Thermal power production   & 360,717 & 384,665 & 368,233 & 392,329 & 388,513 & 412,509 \\
  4 & Electricity from batteries &     954 &     970 &     953 &     952 &     947 &     964 \\
  5 & Overgeneration             &   7,093 &   7,089 &   7,089 &   7,089 &   7,089 &   7,092 \\  
  6 & Not supplied electricity   &       0 &       0 &       0 &       0 &       0 &       0 \\
  7 & Not supplied reserve       &       0 &       0 &       0 &       0 &       0 &       0 \\
  8 & Gas supply                 &     315 &     373 &     315 &     373 &     314 &     374 \\
  9 & PtG plants                 &       0 &       0 &       0 &       0 &       0 &       0 \\
 10 & Gas demand curtailment     &       0 &       0 &       0 &       0 &       0 &       0 \\  
\hline
    & Total second-stage costs   & 401,739 & 425,841 & 409,241 & 433,465 & 429,592 & 453,830 \\  
\hline
\\
\hline
    & Plan                       &     MVP &     MVP &     MVP &     MVP &     MVP &     MVP \\ 
$m$ & Scenarios                  &     L-A &     H-A &     L-B &     H-B &     L-C &     H-C  \\ 
\hline
  1 & Hydro power production     &  31,049 &  31,069 &  31,052 &  31,052 &  31,063 &  31,059 \\
  2 & Start-up of thermal plants &   1,690 &   1,668 &   1,677 &   1,639 &   1,732 &   1,776 \\
  3 & Thermal power production   & 375,917 & 400,496 & 383,781 & 408,334 & 403,418 & 427,934 \\
  4 & Electricity from batteries &     553 &     558 &     552 &     558 &     555 &     564 \\
  5 & Overgeneration             &   3,256 &   3,284 &   3,189 &   3,293 &   3,336 &   3,388 \\  
  6 & Not supplied electricity   &       0 &       0 &       0 & 292,464 & 197,454 & 304,685 \\
  7 & Not supplied reserve       &       0 &       0 &       0 &       0 &       0 &       0 \\
  8 & Gas supply                 &     315 &     374 &     315 &     368 &     310 &     367 \\
  9 & PtG plants                 &       0 &       0 &       0 &       0 &       0 &       0 \\
 10 & Gas demand curtailment     &       0 &       0 &       0 &       0 &       0 &       0 \\  
\hline
    & Total second-stage costs   & 412,781 & 437,451 & 420,568 & 737,707 & 637,869 & 769,774 \\ 
\hline    
\end{tabular}
\label{tab:oper_costs}
\end{table}

\begin{comment}
I due sistemi risultano in due modi completamente diversi di rispondere alla domanda negli scenari che rappresentano l'incertezza. (Scambiare tabella 19 e 20) La tabella 20 mostra i risultati di tale confronto. 

Tabelle 21 e 22: Per mostrare comne il sistema si comporta a regime per soddisfare la domanda annuale nel 2040 (pari a 399.99 TWh), riportiamo la Tabella 21, che mostra le produzioni annuali delle diverse fonti con l'obiettivo di mostrare l'elevata quota di generazione rinnovabile non-programmabile (che copre la quota principale della domanda), quota termica molto bassa e importo molto basso e compensato dall'export. Si nota un maggiore contributo delle risorse di flessibilita' nella soluzione stocastica. (Sottolineare differenze con modello mean value). 
\end{comment}
With a penalty of 3,000 $\frac{\text{\texteuro}}{MWh}$ for energy not delivered \cite{ItalyMarketReform2020}, 
the Value of Stochastic Solution (VSS) is $138,133$ M\texteuro, where $VSS$ is the difference between the total expected cost of the MVP plan ($678,965$ M\texteuro, including decommissioning costs, investment costs and expected operating costs and penalties over the six scenarios) and the value of the objective function ($540,832$ M\texteuro) of the stochastic programming model:
a reduction of 20.34\% in total expected costs is achieved by implementing the investment decisions determined by the SPM plan instead of those determined by the MVP plan.

\begin{figure} [ht!]
    \centering
    \includegraphics[width=0.75\linewidth]{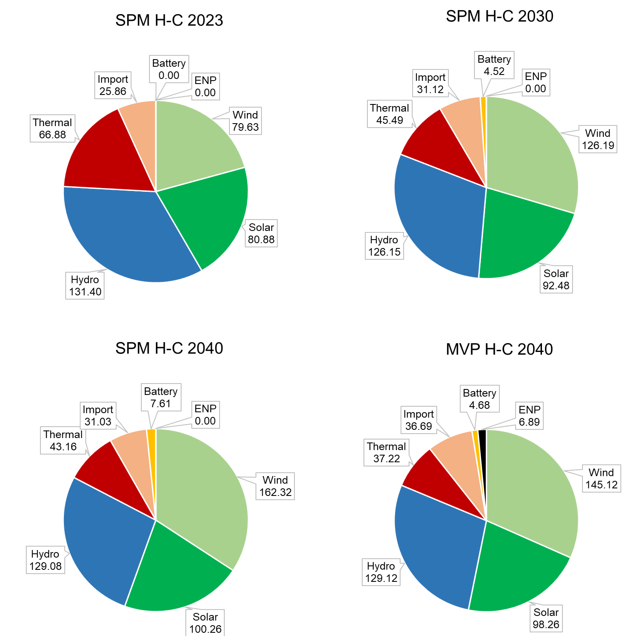}
    \caption{Electricity production [TWh] by power sources in scenario H-C in years 2023, 2030 and 2040 for the SPM plan and in year 2040 for the MVP plan.}
    \label{fig:el_balance}
\end{figure}

To illustrate the evolution of the electricity system resulting from the SPM plan, Figure \ref{fig:el_balance} shows, for scenario H-C, the annual generation per source in the first year of the planning period (year 2023, upper left panel), at the end of the first subperiod (year 2030, upper right panel) and at the end of the second subperiod (year 2040, lower left panel); for comparison, it also shows the generation per source in 2040 of the energy system resulting from the MVP plan (lower right panel).
In the SPM plan, the annual electricity production from non-programmable renewables increases from 160.51 TWh in 2023 to 262.58 TWh in 2040, and the annual energy from batteries, which is zero in 2023, is 4.52 TWh in 2030 and reaches 7.61 TWh in 2040; the annual production of thermoelectric power decreases, while the annual production of hydroelectric power decreases only slightly. In the MVP plan, there is 6.89 TWh of undelivered energy in 2040. Due to the increased interconnection of the electricity system with foreign countries, imported energy increases over the planning period. Details on electricity imports and exports with foreign areas are shown in Table \ref{tab:imp_exp}: in 2040, the MVP plan foresees more imports and less exports than the SPM plan.

%electricity imported from Northern Europe represents more than 90\% of the total imports and exceeds the exports to Northern Europe. In contrast, the exports to Montenegro, Greece and Tunisia exceed the corresponding imports, which are negligible. Overall, in the SPM plan, total exports for the Italian system exceed total imports by 5.79 TWh in 2023, 6.82 TWh in 2030 and 7.07 TWh in 2040. 

\begin{table}
    \caption{Electricity imports and exports [TWh] with foreign zones in the H-C scenario for the SPM plan (years 2023, 2030 and 2040) and the MVP plan (year 2040).}
    \centering
    \begin{tabular}{l |cccc | c | cccc | c}
    \hline
 & \multicolumn{5}{c}{Import from foreign zones}& \multicolumn{5}{| c}{Export to foreign zones}\\
          &UE&  ME&  GR&  TN&  Total&  UE&  ME&  GR&  TN& Total
\\
\hline
          SPM H-C 2023&25.04&  0.57&  0.26&  0.00&  25.86&  25.01&  3.68&  2.97&  0.00& 31.65
\\
          SPM H-C 2030&29.51&  0.69&  0.38&  0.54&  31.12&  26.89&  3.76&  3.13&  4.16& 37.95
\\
          SPM H-C 2040&28.61&  0.96&  0.62&  0.84&  31.03&  26.88&  3.79&  3.20&  4.23& 38.10
\\
          MVP H-C 2040&33.73&  1.24&  0.77&  0.95&  36.69&  24.68&  3.12&  2.23&  3.79& 33.83
\\
\hline
    \end{tabular}
    \label{tab:imp_exp}
\end{table}

\begin{figure}[ht!]
    \centering
    \includegraphics[width=1\linewidth]{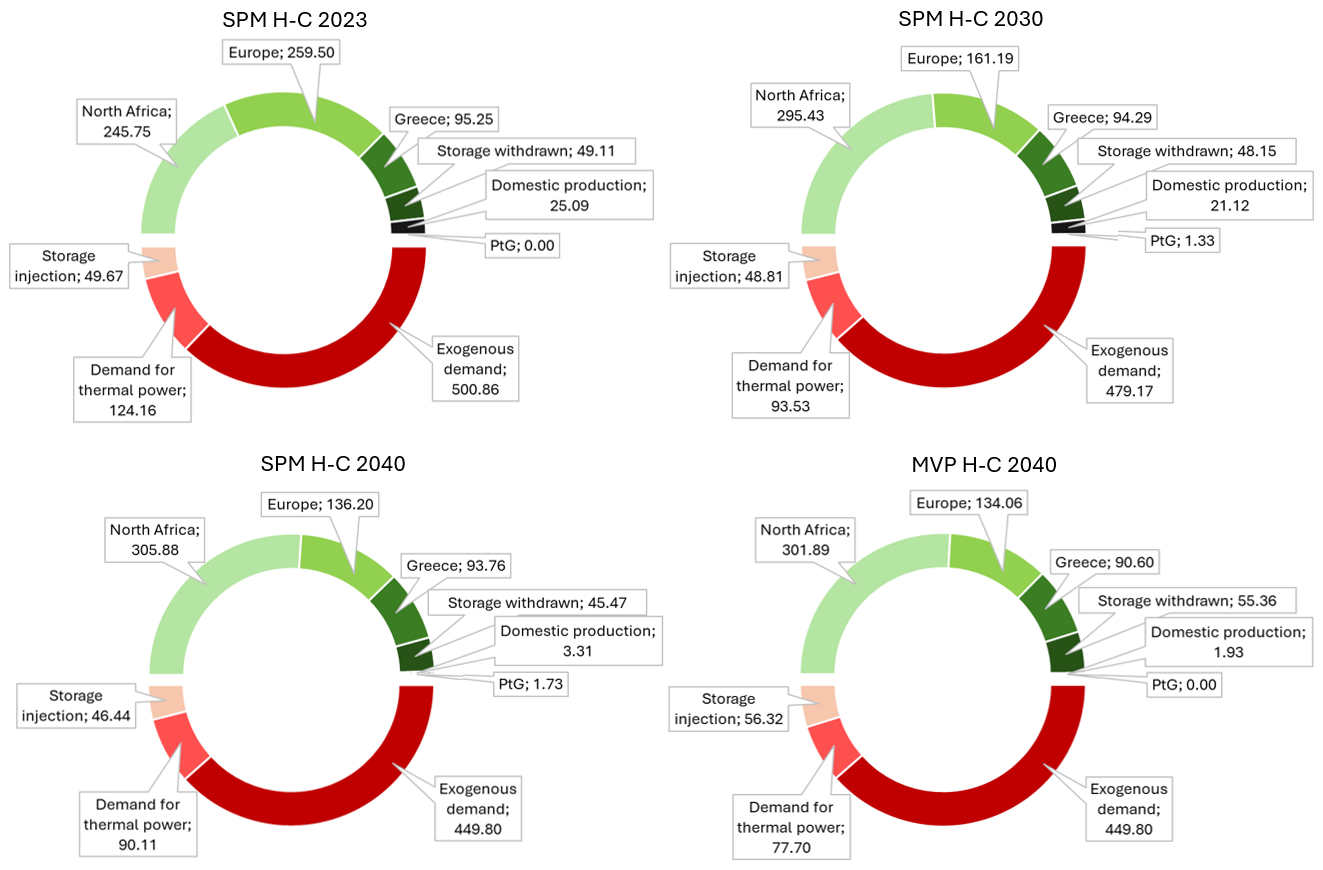}
    \caption{Annual gas supply [TWh$_{th}$] in green and annual gas consumption [TWh$_{th}$] in red in scenario H-C for the SPM plan (years 2023, 2030, 2040) and the MVP plan (year 2040).}
    \label{fig:gas_balance}
\end{figure}

To illustrate the evolution of the gas system, Figure \ref{fig:gas_balance} shows, for the H-C scenario, the annual gas supply in green and the annual gas consumption in red for the SPM plan in 2023 (upper left panel), 2030 (upper right panel) and 2040 (lower left panel); the lower right panel shows the gas supply and consumption in 2040 for the energy system resulting from the MVP plan.
In the SPM plan, gas demand decreases from 625.02 TWh$_{th}$ in 2023 to 539.92 TWh$_{th}$ in 2040, due to both a reduction in demand from the industrial and tertiary sectors and a drastic reduction in gas consumption by thermal power plants (from 124.16 TWh$_{th}$ in 2023 to 90.11 TWh$_{th}$ in 2040) as a result of the increasing penetration of renewables. In 2040, the gas consumption by thermal power plants in the MVP plan is even lower (77.70 TWh$_{th}$) due to electricity demand curtailment. 
In 2023, due to high costs (see Table \ref{tab:Fuel_Price}), domestic production is 25.09 TWh$_{th}$ and most of the gas is imported: 259.50 TWh$_{th}$ from Europe, 245.75 TWh$_{th}$ from North Africa and 95.25 TWh$_{th}$ from Greece. The realisation of the GALSI project in 2024 allows an increase in cheaper gas imports: in 2040, gas imports from North Africa increase to 305.88 TWh$_{th}$, while gas imports from Europe decrease to 136.20 TWh$_{th}$ and domestic production decreases to 3.31 TWh$_{th}$.
In the SPM plan, PtG technology is available from 2030 and produces 1.33 TWh$_{th}$ and 1.73 TWh$_{th}$ synthetic gas in 2030 and 2040 respectively.

\begin{comment}
Qui la rappresentazione utilizzata e' un grafico a barre perche' l'interesse e' sia sul lato della generazione sia sul lato della domanda.
In particolare, sul lato della generazione era interessante sottolineare il contributo delle varie fonti di approvvigionamento del gas.
Sul lato della domanda, era invece interessante mostrare l'entita' della domanda di gas per generazione termoelettrica e la sua variabilita' negli scenari.

Viceversa, nello scenario di alti prezzi, la convenienza nel realizzare il curtailment elettrico piuttosto che attivare nelle ore di picco le centrali termoelettriche si traduce in una riduzione netta della domanda di gas per generazione termoelettrica, richiedendo un ricorso minore alle varie fonti di approvvigionamento di gas.
\end{comment}

\color{black}

\section{Conclusions}
\label{sec:Conclusions}
In this paper, the problem of planning the transition towards decarbonised energy systems is addressed by proposing a two-stage multi-period stochastic programming model to define capacity expansion decisions for facilities both in the power and the gas sector, under different scenarios for future fuel and CO$_2$ prices. 
Specifically, our proposed model co-optimizes investment and operational decisions, taking into account bidirectional flows between the electricity and the natural gas systems.

Two key aspects make the proposed model suitable for assisting qualified operators in anticipative planning: (i) the very detailed representation of short-term operation, and (ii) the inclusion of uncertainty in CO$_2$ and fossil fuel prices. Despite the increased complexity caused by the introduction into the GTEP problem of the detailed modeling of short-term operation and the uncertainty on prices, both aspects need to be taken into account to reliably plan the transition to decarbonized energy systems.
In particular, a detailed consideration of the operational problem is required to simulate the actual operation of the system in the presence of high levels of non-programmable generation, so that the least-cost generation mix that best meets demand and manages the intermittency of renewable generation can be determined. Only by introducing in the GTEP problem an accurate short-term modelling, qualified operators can properly address the challenges associated with integrating large shares of non-programmable renewable.
In terms of stochasticity, taking into account the uncertainty about future prices can significantly improve the quality of the expansion plans. Indeed, performing the anticipative planning with a deterministic model may not always result in generation systems adequate in all scenarios. 

To show the effectiveness of the proposed model, a case study related to the Italian system is considered, where the joint expansion of the Italian electricity and gas system has to be planned under uncertainty of CO$_2$ and fuel prices to meet challenging policy goals set by the European Commission and the Italian regulator.
The numerical results show the relevance of the operating cost on the total system costs, demonstrating the need to represent short-term operation with the highest accuracy when planning the transition to low-carbon energy systems. 
Furthermore, by comparing the expansion plans resulting from the stochastic and the mean-value problem in the different price scenarios, we show how taking account of price stochasticity allows adequate systems to be designed. In particular, the Italian energy system resulting from the stochastic analysis, in contrast to the mean value model solution, can easily cope with stressed load conditions, as the increased renewable generation reduces the magnitude of load peaks, while the increased capacity for electrical storage and gasification allows the effective management of imbalances between generation and demand.

\bigskip

\bibliography{reference}

\newpage
\appendix
\setcounter{table}{0}
\renewcommand{\thetable}{A.\arabic{table}}
\section{Notation}
\label{Notation}
This appendix provides the full notation of the two-stage stochastic programming model introduced in Section \ref{sec:Model}: Table \ref{tab:sets} shows the sets, while in Table \ref{tab:notation} the model parameters and the decision variables are reported.
\begin{table} [ht!]
\caption{Sets.}
    \centering
%    \begin{longtable}{p{2.2 cm} p{11.7 cm}}
\normalsize
    \begin{tabular}{p{0.7 cm} p{11.6 cm}}
    \hline
    Set & Description \\
    \hline
         $\mathcal{Y}$   & Set of years. \\
         $\mathcal{C}^y$ & Set of representative days in year $y \in \mathcal{Y}$. \\ 
         $\mathcal{T}$      & Set of hours of the day $\mathcal{T} = \{ 1,\ldots,24 \}$. \\
         $\mathcal{W}$      & Set of scenarios for fuel and CO$_2$ prices. \\
         $\mathcal{Z}$      & Set of zones of the power system. \\ 
         $\mathcal{Z}_a$    & Set of zones of the power system belonging to area $a \in \mathcal{A}$. \\  
         $\mathcal{A}$      & Set of areas (e.g., countries). \\
         $\mathcal{N}$      & Set of zones of the gas system. \\
         $\mathcal{L}_E$    & Set of existing power transmission lines. \\   
         $\mathcal{L}_C$    & Set of candidate power transmission lines. \\ 
         $BS^{\text{P}}_z$  & Set of power transmission lines entering zone $z \in \mathcal{Z}$ of the power system. \\
         $FS^{\text{P}}_z$  & Set of power transmission lines leaving zone $z \in \mathcal{Z}$  of the power system. \\ 
         $\mathcal{J}_E$    & Set of existing gas pipelines. \\      
         $\mathcal{J}_C$    & Set of candidate pipelines. \\           
         $BS^{\text{G}}_n$  & Set of gas pipelines entering zone $n \in \mathcal{N}$ of the gas system. \\
         $FS^{\text{G}}_n$        & Set of gas pipelines leaving zone $n \in \mathcal{N}$ of the gas system. \\
         $\mathcal{H}_E$          & Set of existing hydropower plants. \\ 
         $\mathcal{H}_C$          & Set of candidate hydropower plants. \\  
%$\mathcal{H}_z \subset \mathcal{H}_E \cup \mathcal{H}_C$
         $\mathcal{H}_z$          & Set of hydropower plants located in zone $z \in \mathcal{Z}$. \\
         $\mathcal{H}_P$          & Set of programmable hydropower plants. \\
         $\mathcal{H}_{NP}$       & Set of non-programmable hydropower plants. \\
%         ${\Xi}$                  & Set of long-term storage controls ${\Xi} = \{1,\ldots,\overline{\xi} \}$. \\
         $\mathcal{K}_z$          & Set of clusters of thermal power plants located in power zone $z \in \mathcal{Z}$. \\
         $\mathcal{K}_z^\text{G}$ & Set of clusters of gas-fired thermal power plants located in power zone $z \in \mathcal{Z}$. \\
         $\mathcal{K}_n$          & Set of clusters of thermal power plants located in gas zone $n \in \mathcal{N}$. \\
         $\mathcal{B}_z$          & Set of battery technologies located in zone $z \in \mathcal{Z}$. \\
         $\mathcal{G}_z$          & Set of PtG technologies located in power zone $z \in \mathcal{Z}$. \\    
         $\mathcal{G}_n$          & Set of PtG technologies located in gas zone $n \in \mathcal{N}$. \\   
         $\mathcal{F}$            & Set of fuels (except gas). \\  
%         $c(d) \in \mathcal{C}^y$ & Representative day associated with calendar day $d \in \mathcal{D}^y, y \in \mathcal{Y}$ \\
%         $f(k) \in \mathcal{F}$ & Fuel employed by thermal power plant $k \in \mathcal{K}_z, z \in \mathcal{Z}$ \\
         \hline
    \end{tabular}
    \label{tab:sets}
%\end{longtable}
\end{table}

%(\texteuro/MW)
%\begin{table}
\begin{longtable}[ht!] {p{1.4 cm} p{1.4 cm} p{9.1 cm}}
%\small % Set the font size for the whole table
%    \centering
%    \begin{tabular}{p{1.5 cm} p{2.5 cm} p{8 cm}}
    \hline
    \multicolumn{3}{l}{Model parameters} \\
    \hline
    Parameter & UM & Description \\
    \hline
         $IC^\text{S}_{z,y}$   & [$\frac{\text{\texteuro}}{\text{MW}}$] & Investment cost of solar power capacity in zone $z \in \mathcal{Z}$ in year $y \in \mathcal{Y}$. \\
         $S_{z,0}$             & [$MW$]             & Solar power capacity in zone $z \in \mathcal{Z}$ at the beginning of the planning period. \\
         $\underline{S}_{z,y}$ & [$MW$]             & Lower bound to solar capacity in zone $z \in \mathcal{Z}$ in year $y \in \mathcal{Y}$. \\  
         $\overline{S}_{z,y}$  & [$MW$]             & Upper bound to solar capacity in zone $z \in \mathcal{Z}$ in year $y \in \mathcal{Y}$. \\ 
         $\mu_{z,t,c}$       & [$\frac{MWh}{MW}$] & Solar power production coefficient in zone $z \in \mathcal{Z}$ in hour $t \in \mathcal{T}$ of representative day $c \in \mathcal{C}^y$, $y \in \mathcal{Y}$. \\ 
    \hdashline
         $IC^\text{W}_{z,y}$   & [$\frac{\text{\texteuro}}{MW}$] & Investment cost of wind power capacity in zone $z \in \mathcal{Z}$ in year $y \in \mathcal{Y}$. \\
         $W_{z,0}$             & [$MW$]             & Wind power capacity in zone $z \in \mathcal{Z}$ at the beginning of the planning period. \\
         $\underline{W}_{z,y}$ & [$MW$]             & Lower bound to wind capacity in zone $z \in \mathcal{Z}$ in year $y \in \mathcal{Y}$. \\  
         $\overline{W}_{z,y}$  & [$MW$]             & Upper bound to wind capacity in zone $z \in \mathcal{Z}$ in year $y \in \mathcal{Y}$. \\ 
         $\rho_{z,t,c}$      & [$\frac{MWh}{MW}$] & Wind power production coefficient in zone $z \in \mathcal{Z}$ in hour $t \in \mathcal{T}$ of representative day $c \in \mathcal{C}^y$, $y \in \mathcal{Y}$. \\ 
    \hdashline 
         $IC_k^\text{K}$       & [$\frac{\text{\texteuro}}{MW}$]   & Investment cost for thermal power plants of type $k \in \mathcal{K}_z, z \in \mathcal{Z}$. \\
         $DC_k^\text{K}$       & [$\frac{\text{\texteuro}}{MW}$]   & Decommissioning cost for thermal power plants of type $k \in \mathcal{K}_z, z \in \mathcal{Z}$. \\
         $C_k^\text{SU}$       & [$\text{\texteuro}$]              & Start-up cost for thermal power plants of type $k \in \mathcal{K}_z, z \in \mathcal{Z}$. \\
         $\underline{P}_k$     & [$MW$]               & Minimum power output for thermal power plants of type $k \in \mathcal{K}_z, z \in \mathcal{Z}$. \\
         $\overline{P}_k$      & [$MW$]               & Capacity of thermal power plants of type $k \in \mathcal{K}_z, z \in \mathcal{Z}$. \\
         $MUT_k$               & [$h$]                & Minimum uptime for thermal power plants of type $k \in \mathcal{K}_z, z \in \mathcal{Z}$. \\ 
         ${MDT}_k$             & [$h$]                & Minimum downtime for thermal power plants of type $k \in \mathcal{K}_z, z \in \mathcal{Z}$. \\
         $N_{k,0}$             & [$-$]                & Number of thermal power plants of type $k \in \mathcal{K}_z, z \in \mathcal{Z}$ at the beginning of the planning period. \\
         $\underline{N}_{k,y}$ & [$-$]                & Lower bound to number of thermal power plants of type $k \in \mathcal{K}_z, z \in \mathcal{Z}$ in year $y \in \mathcal{Y}$. \\
         $\overline{N}_{k,y}$  & [$-$]                & Upper bound to number of thermal power plants of type $k \in \mathcal{K}_z, z \in \mathcal{Z}$ in year $y \in \mathcal{Y}$. \\
         $HR_{k}$              & [$\frac{Gcal}{MWh}$] & Heat rate for thermal power plants of type $k \in \mathcal{K}_z, z \in \mathcal{Z}$. \\
         ${CO_2}_{k}$          & [$\frac{ton}{MWh}$]           & CO$_2$ emission rate for thermal power plants of type $k \in \mathcal{K}_z, z \in \mathcal{Z}$. \\ 
         \hdashline
         $IC_h^\text{H}$       & [$\frac{\text{\texteuro}}{MW}$]         & Investment cost for hydropower plant $h \in \mathcal{H}_C$. \\
         $C_h^\text{H}$             & [$\text{\texteuro}$]         & Operating cost of hydropower plant $h \in \mathcal{H}$. \\
         $\overline{H}_h^\text{IN}$ & [$MW$]    & Maximum hourly pumping power of hydropower plant $h \in \mathcal{H}$. \\
         $\overline{H}_h^\text{OUT}$ & [$MW$]   & Maximum hourly power output of hydropower plant $h \in \mathcal{H}$. \\ 
         $\overline{H}_h^\text{SPILL}$ & [$MW$] & Maximum hourly spillage for hydropower plant $h \in \mathcal{H}$. \\
         $\lambda_h^\text{IN}$ & [$-$]          & Energy loss coefficient for pumping of hydropower plant $h \in \mathcal{H}$. \\ 
         $\lambda_h^\text{OUT}$ & [$-$]        & Energy loss coefficient for power production of hydropower plant $h \in \mathcal{H}$. \\
         $EPR_h$ & [$h$] & Energy to power ratio of hydropower plant $h \in \mathcal{H}$. \\
         $F_{h,t,c,y}$ & [$MW$]                 & Natural inflow of hydropower plant $h \in \mathcal{H}$ in hour $t \in \mathcal{T}$ of representative day $c \in \mathcal{C}^y$ of year $y \in \mathcal{Y}$. \\ 
         $H_{h,0}^\text{LT}$ & [$MWh$]          & Energy in the reservoir of programmable hydropower plant $h \in \mathcal{H}_P$ at the beginning of each year. \\ 
         \hdashline
         $IC_{b,y}^\text{B}$ & [$\frac{\text{\texteuro}}{MW}$]     & Investment cost of battery technology $b \in \mathcal{B}_z, z \in \mathcal{Z}$ in year $y \in \mathcal{Y}$. \\
         $C_b^\text{B}$ & [$\text{\texteuro}$] & Operating cost of battery technology $b \in \mathcal{B}_z, z \in \mathcal{Z}$. \\
         $B_{b,0}^\text{CAP}$ & [$MW$] & Existing capacity of battery technology $b \in \mathcal{B}_z, z \in \mathcal{Z}$ at the beginning of the planning period. \\       
         $\overline{B}_b^\text{CAP}$ & [$MW$] & Upper bound to installed capacity of battery technology $b \in \mathcal{B}_z, z \in \mathcal{Z}$. \\
         $\lambda_b^\text{IN}$ & [$-$] & Energy loss coefficient for charge of battery technology $b \in \mathcal{B}_z, z \in \mathcal{Z}$. \\
         $\lambda_h^\text{OUT}$ & [$-$] & Energy loss coefficient for discharge of battery technology $b \in \mathcal{B}_z, z \in \mathcal{Z}$. \\
         $EPR_b$ & [$h$] & Energy to power ratio of battery technology $b \in \mathcal{B}_z, z \in \mathcal{Z}$. \\
         $B_{b,0,c,y}$ & [$MWh$] & Energy level of battery technology $b \in \mathcal{B}_z, z \in \mathcal{Z}$ at the beginning of representative day $c \in \mathcal{C}$ of year $y \in \mathcal{Y}$. \\ 
         \hdashline
         $IC_g^\text{PtG}$ & [$\frac{\text{\texteuro}}{MW}$] & Investment cost of PtG technology $g \in \mathcal{G}_z, z \in \mathcal{Z}$. \\
         $C_g^\text{PtG}$ & [$\frac{\text{\texteuro}}{MWh_{th}}$] & Operating cost of PtG technology $g \in \mathcal{G}_z, z \in \mathcal{Z}$. \\
         $PtG_{g,0}^\text{CAP}$ & [$MW$] & Existing capacity of PtG technology $g \in \mathcal{G}_z, z \in \mathcal{Z}$. \\       
         $\overline{PtG}_g^\text{CAP}$ & [$MW$] & Upper bound to installed capacity of PtG technology $g \in \mathcal{G}_z, z \in \mathcal{Z}$. \\
         $\eta_g$ & [$-$] & Efficiency of PtG technology $g \in \mathcal{G}_z, z \in \mathcal{Z}$. \\ 
         \hdashline
         $IC_l^\text{L}$ & [$\text{\texteuro}$] & Investment cost for transmission line $l \in \mathcal{L}_E$. \\
         $\underline{F}_l^\text{L}$ & [$MW$] & Lower transmission limit for line $l \in \mathcal{L}$. \\
         $\overline{F}_l^\text{L}$ & [$MW$] & Upper transmission limit for line $l \in \mathcal{L}$. \\ 
         \hdashline 
         $C^\text{ENP}$ & [$\frac{\text{\texteuro}}{MWh}$] & Penalty cost for electricity not provided. \\
         $C^\text{RNP}$ & [$\frac{\text{\texteuro}}{MWh}$] & Penalty cost for reserve not provided. \\
         $C^\text{OG}$ & [$\frac{\text{\texteuro}}{MWh}$] & Penalty cost for over-generation. \\
         $D_{z,t,c,y}^\text{P}$ & [$MWh$] & Electricity demand in power zone $z \in \mathcal{Z}$ in hour $t \in \mathcal{T}$ of representative day $c \in \mathcal{C}$ of year $y \in \mathcal{Y}$. \\
         $R_{z,t,c,y}^\text{P}$ & [$MWh$] & Reserve demand in power zone $z \in \mathcal{Z}$ in hour $t \in \mathcal{T}$ of representative day $c \in \mathcal{C}$ of year $y \in \mathcal{Y}$. \\ 
         \hdashline
         $IC_j^\text{J}$ & [$\text{\texteuro}$] & Investment cost for pipeline $j \in \mathcal{J}_E$. \\
         $\underline{F}_j^\text{J}$ & [$MW$] & Lower transportation limit for pipeline $j \in \mathcal{J}_E$. \\
         $\overline{F}_j^\text{J}$ & [$MW$] & Upper transportation limit for pipeline $j \in \mathcal{J}_E$. \\ 
         \hdashline 
         $C^\text{GC}$ & [$\frac{\text{\texteuro}}{MWh_{th}}$] & Penalty cost for gas curtailment. \\
         $\underline{G}_n$ & [$MW_{th}$] & Minimum hourly gas supply in gas zone $n \in \mathcal{N}$. \\
         $\overline{G}_n$ & [$MW_{th}$]  & Maximum hourly gas supply in gas zone $n \in \mathcal{N}$. \\
         $\overline{G}_n^\text{IN}$ & [$MW_{th}$]  & Maximum hourly gas injection in storage of gas zone $n \in \mathcal{N}$. \\
         $\overline{G}_n^\text{OUT}$ & [$MW_{th}$]  & Maximum hourly gas withdraw from storage of gas zone $n \in \mathcal{N}$. \\
         $\overline{G}_n^\text{LT}$ & [$MWh_{th}$]  & Maximum gas stored in gas zone $n \in \mathcal{N}$. \\
         ${G}_{n,0}^\text{LT}$ & [$MWh_{th}$]  & Gas stored in gas zone $n \in \mathcal{N}$ at the beginning of each year. \\         
         $D^\text{G}_{n,t,c,y}$ & [$MWh_{th}$] & Gas demand in gas zone $n \in \mathcal{N}$ in hour $t \in \mathcal{T}$ of representative day $c \in \mathcal{C}$ of year $y \in \mathcal{Y}$. \\  
         \hdashline
         $M$ & [$-$] & Number of days between two checks of the long-term storage. \\
         $\phi_{a,y}$ & [$-$] & Minimum percentage of annual electric load supplied by non-programmable renewable generation in area $a \in \mathcal{A}$ in year $y \in \mathcal{Y}$. \\
         ${\overline{CO}_2}_{a,y}$ & [$ton$] & Upper bound to CO$_2$ emissions in area $a \in \mathcal{A}$ in year $y \in \mathcal{Y}$. \\ 
         $r$ & [$-$] & Annual discount rate. \\
         $y_0$ & [$-$] & Reference year to which investment costs are discounted. \\
         $\psi_c$ & [$-$] & Weight of representative day $c \in \mathcal{C}^y, y \in \mathcal{Y}$. \\
%         $\xi$ & [$-$] & index for checks of the long-term storage. \\
         \hdashline
         $pb_w$ & [$-$] & Probability of scenario $w \in \mathcal{W}$. \\
         $C_{k,y,w}^\text{M}$ & [$\frac{\text{\texteuro}}{MWh}$] & Production cost for thermal power plants of type $k \in \mathcal{K}_z, z \in \mathcal{Z}$ in year $y \in \mathcal{Y}$ in scenario $w \in \mathcal{W}$. \\
         $C^\text{G}_{n,y,w}$ & [$\frac{\text{\texteuro}}{MWh_{th}}$] & Gas supply cost in gas zone $n \in \mathcal{N}$ in year $y \in \mathcal{Y}$ in scenario $w \in \mathcal{W}$. \\
         $Pr^{\text{FUEL}}_{f,y,w}$ & [$\frac{\text{\texteuro}}{Gcal}$] & Price of fuel $f \in \mathcal{F}$ in year $y \in \mathcal{Y}$ in scenario $w \in \mathcal{W}$. \\
         $Pr^{\text{CO}_2}_{y,w}$ & [$\frac{\text{\texteuro}}{ton}$] & CO$_2$ price in year $y \in \mathcal{Y}$ in scenario $w \in \mathcal{W}$. \\
         $\gamma_{k,0,c,y,w}$ & [$-$] & Number of thermal power plants of type $k \in \mathcal{K}_z, z \in \mathcal{Z}$, online at the beginning of representative day $c \in \mathcal{C}^y$ in year $y \in \mathcal{Y}$ in scenario $w \in \mathcal{W}$. \\
         \hline     
         \hline
    \multicolumn{3}{l}{First-stage decision variables defined for each year $y \in \mathcal{Y}$} \\
    \hline
    Variable & UM & Description \\
    \hline
    $\delta_{l,y}$ & [$-$] & Binary variable: $\delta_{l,y}=1$, if candidate power transmission line $l \in \mathcal{L}_C$ is to be built in year $y$; $\delta_{l,y}=0$, otherwise.\\
    $\theta_{l,y}$ & [$-$] & Binary variable: $\theta_{l,y}=1$, if candidate power transmission line $l \in \mathcal{L}_C$ is available in year $y$; $\theta_{l,y}=0$, otherwise. \\
    $\delta_{j,y}$ & [$-$] & Binary variable: $\delta_{j,y}=1$, if candidate gas pipeline $j \in \mathcal{J}_C$ is to be built in year $y$; $\delta_{j,y}=0$, otherwise. \\
    $\theta_{j,y}$ & [$-$] & Binary variable: $\theta_{j,y}=1$, if candidate pipeline $j \in \mathcal{J}_C$ is available in year $y$; $\theta_{j,y}=0$, otherwise. \\
    $\delta_{h,y}$ & [$-$] & Binary variable: $\delta_{h,y}=1$, if candidate hydropower plant $h \in \mathcal{H}_C$ is to be built in year $y$; $\delta_{h,y}=0$, otherwise. \\
    $\theta_{h,y}$ & [$-$] & Binary variable: $\theta_{h,y}=1$, if candidate hydropower plant $h \in \mathcal{H}_C$ is available in year $y$; $\theta_{h,y}=0$, otherwise. \\
    $N_{k,y}^{-}$ & [$-$] & Integer number of thermal power plants of type $k \in \mathcal{K}_z, z \in \mathcal{Z}$, to be decommissioned. \\
    $N_{k,y}^{+}$ & [$-$] & Integer number of thermal power plants of type $k \in \mathcal{K}_z, z \in \mathcal{Z}$, to be built. \\
    $N_{k,y}$ & [$-$] & Integer number of thermal power plants of type $k \in \mathcal{K}_z, z \in \mathcal{Z}$, available for production. \\    
    $S_{z,y}$ & [$MW$] & New solar power capacity to be built in zone $z \in \mathcal{Z}$. \\
    $W_{z,y}$ & [$MW$] & New wind power capacity to be built in zone $z \in \mathcal{Z}$. \\
    $B_{b,y}^\text{CAP}$ & [$MW$] & New capacity to be built for battery technology $b \in \mathcal{B}_z, z \in \mathcal{Z}$. \\
    $PtG_{g,y}^\text{CAP}$ & [$MW_{th}$] & New capacity to be built for PtG technology $g \in \mathcal{G}_z, z \in \mathcal{Z}$. \\ 
    $RES_{z,t,c,y}$ & [$MW$] & Non-programmable renewable generation in zone $z \in \mathcal{Z}$ in hour $t, 1 \leq t \leq 24$, of representative day $c \in \mathcal{C}^y$. \\
    \hline
    \hline
    \multicolumn{3}{l}{Second-stage decision variables defined for hour $t, 1 \leq t \leq 24$, of representative} \\
    \multicolumn{3}{l}{day $c \in \mathcal{C}^y$ of year $y \in \mathcal{Y}$ under scenario $w \in \mathcal{W}$} \\
    \hline
    Variable & UM & Description \\
    \hline
    $F_{l,t,c,y,w}^\text{L}$ & [$MW$] & Power flow on transmission line $l \in \mathcal{L}$.  \\    
    $F_{j,t,c,y,w}^\text{J}$ & [$MW$] & Gas flow through pipeline $j \in \mathcal{J}$. \\
    $H_{h,t,c,y,w}^{\text{IN}}$ & [$MW$] & Pumping power of hydropower plant $h \in \mathcal{H}$.  \\
    $H_{h,t,c,y,w}^{\text{OUT}}$ & [$MW$] & Power production of hydropower plant $h \in \mathcal{H}$. \\
    $H_{h,t,c,y,w}^{\text{SPILL}}$ & [$MW$] & Spillage for hydropower plant $h \in \mathcal{H}$.  \\
    $\alpha_{k,t,c,y,w}$ & [$-$] & Integer number of thermal power plants of type $k \in \mathcal{K}_z, z \in \mathcal{Z}$ started up. \\
    $\beta_{k,t,c,y,w}$ & [$-$] & Integer number of thermal power plants of type $k \in \mathcal{K}_z, z \in \mathcal{Z}$ shut down. \\
    $\gamma_{k,t,c,y,w}$ & [$-$] & Integer number of thermal power plants of type $k \in \mathcal{K}_z, z \in \mathcal{Z}$ online. \\
    $p_{k,t,c,y,w}$ & [$MW$] & Power production over the minimum production for thermal power plants of type $k \in \mathcal{K}_z, z \in \mathcal{Z}$. \\
    $B_{b,t,c,y,w}$ & [$MWh$] & Energy level of battery technology $b \in \mathcal{B}_z, z \in \mathcal{Z}$. \\
    $B_{b,t,c,y,w}^\text{IN}$ & [$MW$] & Charge of battery technology $b \in \mathcal{B}_z, z \in \mathcal{Z}$. \\
    $B_{b,t,c,y,w}^\text{OUT}$ & [$MW$] & Discharge of battery technology $b \in \mathcal{B}_z, z \in \mathcal{Z}$. \\
    $G_{g,t,c,y,w}^{\text{PtG}}$ & [$MW_{th}$] & Production of synthetic gas from PtG technology $g \in \mathcal{G}_z, z \in \mathcal{Z}$. \\ 
    $G_{n,t,c,y,w}$ & [$MWh_{th}$] & Gas supply in gas zone $n \in \mathcal{N}$. \\  
    $G^\text{IN}_{n,t,c,y,w}$ & [$MW_{th}$] & Gas injected in storage of gas zone $n \in \mathcal{N}$. \\
    $G^\text{OUT}_{n,t,c,y,w}$ & [$MW_{th}$] & Gas withdrawn from storage in gas zone $n \in \mathcal{N}$. \\  
    $OG_{z,t,c,y,w}$ & [$MWh$] & Overgeneration in power zone $z \in \mathcal{Z}$. \\
    $E^\text{NP}_{z,t,c,y,w}$ & [$MWh$] & Electricity not provided in power zone $z \in \mathcal{Z}$. \\   
    $R^\text{NP}_{z,t,c,y,w}$ & [$MWh$] & Reserve not provided in power zone $z \in \mathcal{Z}$. \\
    $G^\text{CURT}_{n,t,c,y,w}$ & [$MWh_{th}$] & Gas curtailment in gas zone $n \in \mathcal{N}$. \\  
    \hline
    \hline
    \multicolumn{3}{l}{Second-stage decision variables defined for each long-term storage control $\xi$,} \\
    \multicolumn{3}{l}{with $1 \leq \xi \leq \overline{\xi}$ and $\overline{\xi}=\lfloor \frac{365}{M} \rfloor$, in year $y \in \mathcal{Y}$ under scenario $w \in \mathcal{W}$} \\
    \hline 
    Variable & UM & Description \\
    \hline
    $H^{\text{LT}}_{h,\xi M,y,w}$ & [$MWh$] & Energy stored at hydropower plant $h \in \mathcal{H}$ at the end of day $\xi M$. \\    
    $G^{\text{LT}}_{n,\xi M,y,w}$ & [$MWh_{th}$] & Gas stored in gas zone $n \in \mathcal{N}$ at the end of day $\xi M$. \\  
         \hline
%    \end{tabular}
\caption{Model parameters and decision variables. UM = Unit of Measure}
\label{tab:notation}
\end{longtable}

\end{document}